\documentclass[10pt]{extarticle}
\usepackage{amsfonts}
\usepackage{amsthm}
\usepackage{amsmath}
\usepackage{amssymb,amsbsy,amsopn}
\usepackage{mathtools}
\usepackage{mathrsfs}
\usepackage{graphicx,color,float,subcaption,caption}
\usepackage{mathptmx}
\usepackage{bbm}
\usepackage{bm}
\usepackage[bookmarks=true,bookmarksnumbered=true,bookmarkstype=toc]{hyperref}
\usepackage{dsfont}
\usepackage{leftidx}
\usepackage{empheq}
\usepackage{enumerate}
\usepackage{chngcntr}
\usepackage{comment}

\usepackage{tikz}
\usepackage{pgfplots}
\pgfplotsset{every axis/.append style={
		axis x line=middle,    
		axis y line=middle,    
		axis line style={<->}, 
	},
	compat=1.12,
	cmhplot/.style={color=black,mark=none,line width=1pt,<->},
	soldot/.style={color=black,only marks,mark=*},
	holdot/.style={color=black,fill=white,only marks,mark=*},
}

\theoremstyle{definition} 

\newtheorem{theorem}{Theorem}

\newtheorem{proposition}[theorem]{Proposition}
\newtheorem{example}[theorem]{Example}

\numberwithin{definition}{section}
\numberwithin{theorem}{section}

\makeatletter
\@addtoreset{equation}{section}
\makeatother

\title{Numerical aspects of shot noise representation of infinitely divisible laws and related processes\footnotetext[0]{The authors are grateful to the anonymous referee for their valuable feedback that helped improve the quality of this manuscript.}}
\author{\sc Sida Yuan\footnote{Email address: sida.yuan@cba.com.au. Postal address: 5-7 Central Avenue, Sydney, NSW 2015, Australia. Sida Yuan is currently employed by the Commonwealth Bank of Australia. The opinions expressed by the authors are solely their own and do not reflect those of the Commonwealth Bank of Australia and its related entities.} \, and Reiichiro Kawai\footnote{Email address: raykawai@g.ecc.u-tokyo.ac.jp. Postal Address: Graduate School of Arts and Sciences / Mathematics and Informatics Center, The University of Tokyo, Japan, and School of Mathematics and Statistics, The University of Sydney, Australia. This work was partially supported by JSPS Grants-in-Aid for Scientific Research 20K22301 and 21K03347.}}
\date{}

\begin{document}

\maketitle

\begin{abstract}
	\noindent The ever-growing appearance of infinitely divisible laws and related processes in various areas, such as physics, mathematical biology, finance and economics, has fuelled an increasing demand for numerical methods of sampling and sample path generation. In this survey, we review shot noise representation with a view towards sampling infinitely divisible laws and generating sample paths of related processes. In contrast to many conventional methods, the shot noise approach remains practical even in the multidimensional setting. We provide a brief introduction to shot noise representations of infinitely divisible laws and related processes, and discuss the truncation of such series representations towards the simulation of infinitely divisible random vectors, L\'evy processes, infinitely divisible processes and fields and L\'evy-driven stochastic differential equations. Essential notions and results towards practical implementation are outlined, and summaries of simulation recipes are provided throughout along with numerical illustrations. Some future research directions are highlighted.
	
	\vspace{0.3em}
	\noindent \textit{Keywords:} Infinitely divisible laws; L\'evy processes; shot noise representation; infinitely divisible processes, Monte Carlo methods.
	
	\noindent \textit{2020 Mathematics Subject Classifications:} 60E07, 60G55, 60G51, 65C05, 65C30.
\end{abstract}

\tableofcontents

\section{Introduction}
Infinitely divisible laws have long been investigated actively in the literature due to their fascinatingly rich structure. The corresponding class of stochastic processes is the class of L\'evy processes, and moreover, infinitely divisible laws also bear deep relations to infinitely divisible processes and fields and L\'evy-driven stochastic differential equations. Over the last half-century, such stochastic processes have grown to become widely popular across many domains. The attractiveness of such stochastic processes may be attributed to two reasons. Firstly, they are able to capture jump discontinuities with great versatility. Through the theory of these stochastic processes, one can construct stochastic processes with flexible jump structures under mild technical conditions. Secondly, stochastic processes related to infinitely divisible laws can capture dynamics beyond Gaussianity. For example, the subclass of such stochastic processes with heavy-tailed marginal laws offers an easy solution for the modelling of heavy-tailed dynamics. To mention only a handful of applications, these stochastic processes have appeared in physical modelling for diffusion and transport \cite{carnaffan2017cusping,stanislavsky2009transport}, degradation modelling \cite{almasary2017approximate,kahle2016degradation} and various applications in finance and insurance \cite{carr2002fine,tankov2003financial,peters2015advances}.

As stochastic processes relating to infinitely divisible laws find increasing applications in the literature, there is a clear demand for methods of sampling infinitely divisible laws and generating sample paths of related processes. One potential solution is the approximation of sample paths based on classical deterministic time discretisation. However, a major drawback is that many application contexts require simulation techniques which specify individual jumps \cite{kahle2016degradation,peters2015advances}. For example, in insurance, claims can be modelled as downward jumps, thus preserving jumps is necessary to observing the ruin time. This necessity for preserving all or some of the jumps of the stochastic process serves as an additional hurdle in our quest for appropriate numerical schemes, which eliminates such conventional methods based upon increments from consideration. Therefore, simulation of infinitely divisible laws and related stochastic processes by generating individual jumps seems not only ideal, but perhaps necessary in many practical scenarios.

In physics, the term \textit{shot noise} is used to describe noise resulting from the discreteness of charge carriers. In electrical circuits, shot noise manifests as the sporadic fluctuations of current, particularly when the current is low \cite{beenakker1992suppression,reznikov1998quantum}. In optics, shot noise manifests as the fluctuations of the number of photons detected, most apparent in a low light environment \cite{beenakker1999photon,echternach2013photon}. Such discrete noise is modelled by the shot noise process, in which the arrival times of the \textit{shots} or \textit{jumps} follow a Poisson process. More precisely, if $X_t$ denotes the system's state at time $t$, then it is most commonly expressed through one of the two following representations:

\begin{itemize}
	\setlength{\parskip}{0cm}
	\setlength{\itemsep}{0cm}
	\item where $\mu(d\mathbf z,ds)$ is a marked Poisson random measure which has a weight at $(\mathbf z,s)$ if there is a jump at time $s\ge0$ with size $\mathbf z\in\mathbb R^d\backslash\{0\}$, we write
	\[
		X_t = \int_0^t\int_{\mathbb R^d\backslash\{0\}} H(t,s,\mathbf z)\,\mu(d\mathbf z,ds);
	\]
	
	\item where $\{\Gamma_k\}_{k\in\mathbb N}$ is the sequence of standard Poisson arrival times independent of the iid marks $\{Z_k\}_{k\in\mathbb N}$ corresponding to $\mu(d\mathbf z,ds)$, we write
	\[
		X_t = \sum_{k=1}^{+\infty} H(t,\Gamma_k,Z_k).
	\]
\end{itemize}

Here, the kernel $H(t,s,\mathbf z)$ of the shot noise process describes the level, at an observation time $t$, of the shot that occurred at a previous time $s$. In physical applications, it often makes sense that the influence of the shot on the physical system decays with the passage of time. In those settings, the magnitude of the kernel is nonincreasing in the observation time $t$, and often taken as exponential decay. The shot noise phenomenon has been investigated in a wide variety of applications, for example, in metallic conductors \cite{beenakker1992suppression}, quantum systems \cite{beenakker2003quantum,reznikov1998quantum} and optics \cite{beenakker1999photon,echternach2013photon,sears2012photon,wilt2013photon}. The modelling of shot noise has also been extended mathematically, for example, to Cox processes \cite{brix1999generalized}, capturing long-range dependence \cite{bremaud2002power} and nonlinearity \cite{eliazar2005nonlinear}. Shot noise processes have been shown to relate to other important stochastic processes in the asymptotic regime, such as the fractional Brownian motion \cite{kluppelberg2004fractional} and nonstationary Gaussian processes \cite{pang2019nonstationary}. Crucially, these shot noise processes have a profound connection to infinitely divisible laws \cite{lane1984central}. By considering the shot noise kernel as a cumulative integral of a L\'evy measure and the underlying Poisson arrival times as random steps over the domain of the kernel rather than over time, one obtains a shot noise series representation of the infinitely divisible law characterised solely by a L\'evy measure. This straightforwardly leads to shot noise representations for L\'evy processes and infinitely divisible processes, which are in some sense analogous to the well-known Karhunen-Lo\`eve expansion for Gaussian processes.

Paralleling advancements in the study of shot noise processes \cite{daley1971definition,gilchrist1975shot,hawkes1974cluster,rice1977generalized,westcott1976existence}, shot noise representations of infinitely divisible laws and processes have been investigated as early as the 1970s.
Ferguson and Klass \cite{ferguson1972representation} led the seminal effort in establishing an initial method of representing independent increment processes without Gaussian components as random series. In response, Kallenberg \cite{kallenberg1974series} investigated their convergence properties, while Resnick \cite{resnick1976extremal} demonstrated their derivation via the L\'evy-It\^o decomposition of such stochastic processes. Their theoretical development and appearance have since gradually expanded, for instance, \cite{giraitis1991shot,samorodnitsky1991construction,vervaat1979stochastic,walker2000miscellanea}. In particular, Rosi\'nski established general necessary and sufficient conditions for almost sure convergence of shot noise series to infinitely divisible random vectors without Gaussian components \cite{rosinski1990series}. Shot noise representation was used in \cite{rosinski1989path} to study path properties of L\'evy-driven stochastic integrals. Almost sure uniform convergence for shot noise representations of L\'evy processes was investigated in the general setting in \cite{rosinski2001series}. Since then, shot noise representation has been at the forefront of the study of a variety of relevant stochastic processes, such as the stable process \cite{davydov2012convergence,lemke2015fully}, its generalisations via tempering \cite{bianchi2011tempered,houdre2007layered,rosinski2007tempering}, fractional stable motions \cite{carnaffan2019analytic,houdre2006fractional,kawai2016higher,leguevel2012ferguson} and L\'evy processes of type G \cite{wiktorsson2002simulation}.

Perhaps, the greatest gift of shot noise representation to the study of stochastic processes in the age of unprecedented computational power is its elegant and simple solution to sample path generation. In the case of L\'evy and infinitely divisible processes, by truncating the series representation to a finite sum, one obtains a straightforward approximation of the stochastic process from which approximate sample paths can be generated. Among the few sample path generation schemes \cite{rosinski2008simulation}, this has been the go-to numerical method for contemporary applications involving L\'evy processes, \cite{kahle2016degradation,kawai2006importance,lecourtois2018some,peters2015advances} to name a few more examples. Moreover, truncation of shot noise can be naturally extended to the multidimensional setting, including L\'evy copulas \cite{grothe2013vine,tankov2016levy}. As one would expect, the widespread use of the numerical technique demands analyses of the associated truncation error. To give some examples of particular stochastic processes, error analyses have been performed for the stable process \cite{bondesson1982simulation}, gamma process \cite{kahle2016degradation}, tempered stable process \cite{imai2011finite}, higher order fractional stable motion \cite{kawai2016higher} and L\'evy-driven CARMA processes \cite{kawai2017sample}.
More general treatments of error analysis have been studied, for example, in terms of moments \cite{imai2013numerical} and Gaussian approximation \cite{asmussen2001approximations,cohen2007gaussian}. As such, using numerical schemes based on shot noise representations carries the advantages of tractable generalisability to multidimensional settings and error analysis. Through the present survey, we hope to clearly establish the practical value of numerical methods for simulating infinitely divisible laws and related processes based on shot noise representations. In doing so, we demonstrate the viability of using jump models in applied contexts, and encourage further enrichment of the technique.

This survey aims to summarise shot noise representation with a view towards sampling infinitely divisible laws and generating sample paths of related processes.
We review some preliminary notions of infinitely divisible laws and related processes in Section \ref{section preliminaries}, and offer several important examples of infinitely divisible laws. Section \ref{section series representation} outlines shot noise representations of infinitely divisible laws via the L\'evy-It\^o decomposition. Examples with various infinitely divisible laws are provided. We describe the approximation of infinitely divisible laws via truncation of shot noise representation in Section \ref{section truncation}, along with important results on the error. In Section \ref{section truncation for processes}, we discuss the truncation scheme for L\'evy processes, infinitely divisible processes and fields and L\'evy-driven stochastic differential equations. Various examples of error analysis and numerical illustrations are presented, along with summaries of simulation recipes. Section \ref{section computing expectations} briefly visits some practical numerical topics for computing expectations via shot noise representations, such as various variance reduction methods and using quasi-Monte Carlo methods. Finally, we summarise our discourse in Section \ref{section conclusion} along with brief suggestions for future research directions.

\section{Preliminaries}\label{section preliminaries}
We begin by reviewing some preliminaries of the theory of infinitely divisible laws and related processes. Some well-known examples of particular interested to the literature are given.

We introduce some notations which will be used throughout. In what follows, we will be working under the probability space $(\Omega,\mathscr{F},\mathbb P)$. We denote the Borel $\sigma$-algebra over a space $S$ as $\mathcal B(S)$. Let $\mathbb N\coloneqq\{1,2,\cdots\}$ and $\mathbb N_0 \coloneqq\{0,1,2,\cdots\}$. Denote $\langle\cdot,\cdot\rangle$ as the inner product, $\|\cdot\|$ as the Euclidean norm on $\mathbb R^d$ for any $d\in\mathbb N$ and $\mathbb R^d_0\coloneqq\mathbb R^d\backslash\{0\}$. The Dirac delta measure concentrated at $\mathbf x\in\mathbb R^d$ is denoted by $\delta_{\mathbf x}$. We denote the positive part of functions as $(f(x))_+ \coloneqq0\vee f(x)$. Let $\stackrel{\mathscr L}{=}$ represent equality in law and $\mathscr L(\cdot)$ the law of a random vector. We denote the indicator function of a set $A$ as $\mathbbm1_A(\cdot)$ and sometimes as $\mathbbm1(\cdot\in A)$. We refer to the uniform distribution over $(0,1)$ and the exponential distribution with unit rate as the standard uniform and exponential distributions, respectively. The sequence $\{\Gamma_k\}_{k\in\mathbb N}$ will be used throughout to denote the arrival times of the standard Poisson process.

\subsection{Infinitely divisible laws and related processes}
Perhaps the most familiar definition of infinite divisibility is as follows. A law $F$ is said to be \textit{infinitely divisible} if for every $n\in\mathbb N$, there exists a sequence of iid random vectors $\{X_{k,n}\}_{k=1,\dots,n}$ such that $F = \mathscr L(\sum\nolimits_{k=1}^n X_{k,n})$. Immediately from this definition, we see that if $\varphi$ is the characteristic function of a random vector $X$, then $X$ is infinitely divisible if and only if there exists a characteristic function $\varphi_n$ for every $n\in\mathbb N$ such that $(\varphi_n)^n \equiv \varphi$. This criterion is often useful in determining the infinite divisibility of a law given that we know its characteristic function. Furthermore, the characteristic function of infinitely divisible laws can provide even more insight through the following celebrated result on their characterisation.

\begin{theorem}[L\'evy-Khintchine representation]\label{theorem LK representation}
	A probability law $F$ is infinitely divisible if and only if there exists a triple $(\mathbf{a},S,\nu)$, where $\mathbf{ a}\in\mathbb R^d$, $S\in\mathbb R^{d\times d}$ is a symmetric nonnegative-definite matrix and $\nu(d\mathbf z)$ is a measure on $\mathbb R^d_0$ such that
	\begin{equation}\label{Levy measure condition}
		\int_{\mathbb R^d_0}(1\wedge\|\mathbf z\|^2)\, \nu(d\mathbf z)<+\infty,
	\end{equation}
	and the characteristic function of $F$ is given by
	\begin{equation}\label{Levy-Khintchine formula}
		\varphi(\bm\theta) = \exp\left[i\langle\bm\theta,a\rangle - \frac{1}{2}\langle\bm\theta,S\bm\theta\rangle + \int_{\mathbb R^d_0}\left(e^{i\langle\bm\theta,\mathbf z\rangle} - 1 - i\langle\bm\theta,\mathbf z\rangle \mathbbm1_{(0,1]}(\|\mathbf z\|)\right)\, \nu(d\mathbf z)\right],\quad\bm\theta\in\mathbb R^d.
	\end{equation}
	Moreover, if it exists, the triple $(a, S,\nu)$ is unique. 
\end{theorem}

A rigorous proof of the Theorem \ref{theorem LK representation} can be found in \cite[Section 8]{sato1999levy}. We call the measure $\nu(d\mathbf z)$ which satisfies the integrability condition \eqref{Levy measure condition} a \textit{L\'evy measure}. A stochastic process $\{X_t:t\ge0\}$ in $\mathbb R^d$ with $X_0 = 0$ a.s. is a \textit{L\'evy process} if
\begin{enumerate}[(i)]
	\setlength{\parskip}{0cm}
	\setlength{\itemsep}{0cm}
	\item it has stationary and independent increments,
	\item it is continuous in probability, that is, for any $\epsilon>0$ and $t\ge0$, it holds that $\lim_{\Delta\to0}\mathbb P(\|X_{t+\Delta}-X_t\|>\epsilon)=0$, and
	\item $t\mapsto X_t$ is c\`adl\`ag.
\end{enumerate}
 
While the general infinitely divisible law and L\'evy process may contain a Gaussian component, our focus is restricted to the absence of such components, that is, when $S=0$ in the L\'evy-Khintchine formula \eqref{Levy-Khintchine formula}.
Among the most elementary examples of L\'evy processes are the Poisson process and its generalisation, the compound Poisson process. The (compound) Poisson process is a L\'evy process with a (compound) Poisson marginal law. A random vector $X$ is distributed under a compound Poisson law if and only if it can be expressed as the random sum $X = \sum_{k=1}^N Y_k$, where $N$ is a Poisson random variable with rate $\lambda>0$ and $\{Y_k\}_{k\in\mathbb N}$ is a sequence of iid random vectors with law $\rho(d\mathbf z)$ independent of $N$, with the characteristic function
\begin{equation}\label{cf of compound Poisson law}
	\varphi_X(\bm\theta) = \exp\left[\int_{\mathbb R_0^d}\left(e^{i\langle\bm\theta,\mathbf z\rangle} - 1\right)\lambda\,\rho(d\mathbf z)\right],\quad\bm\theta\in\mathbb R^d.
\end{equation}

The connection between L\'evy processes and infinitely divisible laws runs deep. Specifically, the relationship is a correspondence. Where $\{X_t:t\ge0\}$ is a L\'evy process, it holds that $\varphi_{X_t} \equiv (\varphi_{X_1})^t$, so every increment of a L\'evy process is infinitely divisible. Conversely, every infinitely divisible law admits the existence of a L\'evy process with a matching marginal law.
To interpret the L\'evy-Khintchine representation in the context of L\'evy processes, note that the characteristic function $\varphi$ \eqref{Levy-Khintchine formula} corresponds to a convolution of a Brownian motion and many Poisson processes of different jump sizes and intensities governed by the L\'evy measure $\nu(d\mathbf z)$. The vector $\mathbf a$ corresponds to a linear drift while the matrix $S$ corresponds to the diffusion matrix of the Brownian motion. With respect to the Poisson component, the L\'evy measure $\nu(B)$, for any $B\in\mathcal B(\mathbb R^d_0)$, corresponds to the expected number of jump sizes in $B$ in a unit time interval.
The term $i\langle\bm\theta,\mathbf z\rangle$ in the integrand is a centring term for the convolved Poisson processes with small jump sizes, which ensures that the integral exists in the case where $\nu(d\mathbf z)$ is an infinite L\'evy measure with a heavy density near the origin. Specifically, since $\|e^{i\langle\bm\theta,\mathbf z\rangle}-1\| \sim |\langle\bm\theta,\mathbf z\rangle|$ as $\|\mathbf z\|\to0$, if linear functions are not $\nu$-integrable near the origin, then the compensation terms are necessary and cannot be integrated separately in general. An exception is the so-called subordinator where the L\'evy measure has positive support and finite first moment around the origin. The subordinator forms an important subclass of L\'evy processes, which can be employed to express random clocks. In this survey, the stable process (L\'evy measure \eqref{stable Levy measure}) with stability $\alpha\in(0,1)$ without negative jumps, tempered stable process (L\'evy measure \eqref{tempered stable Levy measure}) with stability $\alpha\in(0,1)$ without negative jumps and the gamma process (L\'evy measure \eqref{gamma Levy measure}) are examples of subordinators.

To summarise, a L\'evy process is generally comprised of a diffusion component, which is Gaussian, and a jump component, which can be decomposed into compound and compensated Poisson components for large and small jumps, respectively. For a comprehensive review of the theory of L\'evy processes, we refer the reader to \cite{bertoin1998levy,sato1999levy}. An important takeaway relevant to our discussion later is as follows: an infinitely divisible random vector $X$ characterised by the triplet $(0,0,\nu)$ can be understood as the (possibly infinite) sum of the jumps of a L\'evy process characterised by the same L\'evy-Khintchine triple over the unit time interval. This allows us to make sense of the concept of \textit{jumps} in the context of an infinitely divisible random vector.

A related concept to infinitely divisible laws and L\'evy processes is the notion of \textit{infinitely divisible processes} \cite{rajput1989spectral,rosinski1989path,talagrand1993regularity}. A stochastic process is infinitely divisible if its finite dimensional distributions are infinitely divisible. Naturally, all L\'evy processes are infinitely divisible. However, the class of infinitely divisible processes is more general, as it includes the following. A stochastic process $\{X_t:t\ge0\}$ in $\mathbb R^d$ is a \textit{stochastic integral process} if it can be represented in the stochastic integral form as
\begin{equation}\label{stochastic integral process}
	\{X_t:t\ge0\} \stackrel{\mathscr L}{=} \left\{\int_S f(t,s)\,\Lambda(ds) : t\ge0\right\},
\end{equation}
where $f$ is a suitable deterministic function and $\Lambda(ds)$ is an independently scattered infinitely divisible measure on a suitable space. 
We refer the reader to Theorems 4.11 and 5.2 in \cite{rajput1989spectral} for details regarding the stochastic integral form \eqref{stochastic integral process}. An important example of such a random measure is that generated by the increments of an additive process $\{Z_t:t\in S\}$, where $S$ is a possibly unbounded interval. Moreover, by imposing stationary increments, the infinitely divisible measure corresponds to a L\'evy process, which is most relevant to our discussion.
Stochastic integral processes driven by L\'evy processes of infinity jump activity are of interest in Section \ref{section Levy-driven stochastic integrals}. We exclude the case of finite jump activity from the discussion, as approximations are not required in that setting. We mention that stochastic integral processes of the form \eqref{stochastic integral process} where the integrator is not independently scattered are also of interest in the literature, such as in the case of a cluster compound Poisson random measure \cite{samorodnitsky1996class}.

\subsection{Examples of infinitely divisible laws without Gaussian components}\label{section ID examples}
We provide some examples of infinitely divisible laws, and equivalently, of L\'evy processes. In particular, our definitions will be provided at the level of the L\'evy measure, as this not only exemplifies the immense utility of the L\'evy-Khintchine formula \eqref{Levy-Khintchine formula}, but also sets up for their later use as demonstrations of shot noise representations in Section \ref{section series representation}.

We begin by noting that a compound Poisson law is an infinitely divisible law without Gaussian components where the L\'evy measure $\nu(d\mathbf z)$ is finite. For example, if $\nu \equiv \delta_1$, then we have the standard Poisson distribution.
Of course, we obtain the compound Poisson process and standard Poisson process when we consider the L\'evy-Khintchine triple in the context of L\'evy processes. The compound Poisson distribution and process have found use in a variety of applications \cite{kahle2016degradation,peters2015advances}.

We now define the stable law according to its L\'evy measure \cite[Section 14]{sato1999levy}.	A \textit{stable law with stability $\alpha\in(0,2)$ and scale $a>0$} is an infinitely divisible law without Gaussian components and with its L\'evy measure given by
\begin{equation}\label{stable Levy measure}
	\nu(B) = \int_{S^{d-1}}\int_{0+}^{+\infty}\mathbbm1_B(r\bm\xi)\frac{a}{r^{\alpha+1}}\, dr\, \sigma(d\bm\xi),\quad B\in\mathcal B(\mathbb R^d_0),
\end{equation}
where $S^{d-1}$ is the unit sphere of $\mathbb R^d$ and $\sigma(d\bm\xi)$ a probability measure on $S^{d-1}$.
The stable law, which contains and generalises the Cauchy ($\alpha=1$) and Gaussian ($\alpha\to2$) laws, has received widespread attention in the literature due to its usefulness as a heavy-tailed distribution. Similarly, the corresponding stable process has also appeared in wide-ranging applications \cite{peters2015advances,stanislavsky2009transport}, to name a couple. See \cite{samorodnitsky1994stable} for a comprehensive review of its properties. The simulation of stable processes have been thoroughly studied in the literature, for example, in \cite{janicki1993simulation,samorodnitsky1994stable}. It is important to note that the L\'evy measure \eqref{stable Levy measure} is expressed in polar form, that is, the integrand $a/r^{\alpha+1}$ captures the expected number of jumps with magnitude $r$ over any unit time interval, while $\sigma(C)$ captures the proportion of jumps in the set of directions $C\subseteq S^{d-1}$. The shot noise representation for the stable law is provided in Example \ref{example stable series representation} and error analysis of its truncation is given in Example \ref{example stable error analysis}.

As a heavy-tailed distribution, the stable law does not have second-order moments for $\alpha\in(0,2)$, and no first-order moment for $\alpha\in(0,1]$. One can construct a law which preserves all moments and yet resembles many of the properties of the stable law by truncating the L\'evy density \cite{mantegna1994stochastic}. Another elegant approach is through the exponential tempering \cite{koponen1995analytic,rosinski2007tempering} of the L\'evy measure \eqref{stable Levy measure}, of which we summarise in the following. Suppose we apply a polar decomposition on the stable L\'evy measure $\nu_S(d\mathbf z)$ to obtain $\nu_S(dr,d\bm\xi) = h(dr,\bm\xi)\sigma(d\bm\xi)$, where $\{h(\cdot,\bm\xi)\}_{\bm\xi\in S^{d-1}}$ is an appropriate family of L\'evy measures defined on $(0,+\infty)$ and $\sigma(d\bm\xi)$ remains a probability measure on $S^{d-1}$. Then, the infinitely divisible law without Gaussian components with L\'evy measure $\nu(dr,d\bm\xi) = q(r,\bm\xi)\nu_S(dr,d\bm\xi)$ is called a \textit{tempered law}, where $r\mapsto q(r,\bm\xi)$ is completely monotone and $\lim_{r\to+\infty}q(r,\bm\xi) = 0$ for every $\bm\xi\in S^{d-1}$. By the complete monotonicity, the tempering function can be represented as $q(r,\bm\xi) = \int_{0+}^{+\infty} e^{-rs}\,Q(ds,\bm\xi)$,
where $\{Q(ds,\bm\xi)\}_{\bm\xi\in S^{d-1}}$ is a family of finite Borel measures on $(0,+\infty)$. Tempering the stable L\'evy measure \eqref{stable Levy measure} with the function $q(r,\bm\xi)$, we obtain the L\'evy measure
\begin{equation}\label{tempered stable Levy measure}
	\nu(B) = \int_{S^{d-1}}\int_{0+}^{+\infty} \mathbbm1_B(r\bm\xi) \frac{1}{r^{\alpha+1}}q(r,\bm\xi)\,dr\,\sigma(d\bm\xi) = \int_{\mathbb R_0^d}\int_{0+}^{+\infty} \mathbbm1_B(u\mathbf v)\frac{e^{-u}}{u^{\alpha+1}}\,du\,\rho(d\mathbf v),\quad B\in\mathcal B(\mathbb R_0^d),
\end{equation}
where the measure $\rho(d\mathbf v)$ satisfies
\begin{equation}\label{tempered stable law rho equation}
	\rho(A) = \int_{S^{d-1}}\int_{0+}^{+\infty} \mathbbm1_A(\bm\xi/s)s^\alpha\,Q(ds,\bm\xi)\,\sigma(d\bm\xi),\quad A\in\mathcal B(\mathbb R_0^d).
\end{equation}
We refer the reader to \cite{rosinski2007tempering} for the proofs for the formulas \eqref{tempered stable Levy measure} and \eqref{tempered stable law rho equation}.
Thus, we call an infinitely divisible law without Gaussian components with L\'evy measure of the form \eqref{tempered stable Levy measure} a \textit{tempered stable} law, where $\alpha\in(0,2)$ is the stability parameter and the measure $\rho(d\mathbf v)$ satisfies $\int_{\mathbb R_0^d}\|\mathbf v\|^\alpha\,\rho(d\mathbf v)<+\infty$.
A fascinating attribute of the tempered stable process is that it behaves like a stable process in short time and a Brownian motion in long time. The tempered stable process has found applications in financial modelling \cite{carr2002fine,lecourtois2018some} and statistical mechanics \cite{carnaffan2017cusping,carnaffan2019analytic}. Shot noise representations of the tempered stable law are presented in Example \ref{example tempered stable law} and comparison of truncation errors for the representations are provided in Example \ref{example Rosinski's series comparison}.

In the one-dimensional setting, the \textit{CGMY} law is related to the tempered stable law, which is associated with the CGMY process introduced in \cite{carr2002fine} as a model for asset returns, governed by the L\'evy measure
\begin{equation}\label{CGMY Levy measure}
	\nu(dz) = C\frac{e^{-G|z|}\mathbbm1_{(-\infty,0)}(z) + e^{-M|z|}\mathbbm1_{(0,+\infty)}(z)}{|z|^{1+Y}}\,dz,
\end{equation}
where $C>0$, $G\ge0$, $M\ge0$ and $Y<2$. Clearly, the CGMY L\'evy measure \eqref{CGMY Levy measure} with $G,M>0$ and $Y\in(0,2)$ can be expressed in terms of the tempered stable L\'evy measure \eqref{tempered stable Levy measure} with $\rho(dv) = C(G^Y\delta_{\{-1/G\}}(dv) + M^Y\delta_{\{1/M\}}(dv))$. We mention here that shot noise representations for the tempered stable law (Example \ref{example tempered stable law}) do not cover CGMY laws with $Y\le0$.

An important case of the CGMY law is when $Y=0$. A \textit{gamma law with shape $a>0$ and scale $\beta>0$} is an infinitely divisible law without Gaussian components and with the L\'evy measure
\begin{equation}\label{gamma Levy measure}
	\nu(dz) = \frac{a}{z}e^{-\beta z}\, dz,
\end{equation}
defined on $(0,+\infty)$. The gamma process has seen applications in various areas including degradation modelling \cite{kahle2016degradation} and statistical mechanics \cite{carnaffan2017cusping}. We emphasize the distinctness of properties between the gamma and tempered stable laws, despite both being cases of the CGMY law with finite moments of all polynomial orders. For example, we will later see immense differences in their shot noise representation (Examples \ref{example gamma law} and \ref{example tempered stable law}), and the inapplicability of Gaussian approximation for the truncation error in the case of the gamma law (Section \ref{section Gaussian approximation}).

Following a similar idea of tempering the stable L\'evy measure \eqref{stable Levy measure}, another direction of generalising the stable law is through the layered stable law \cite{houdre2007layered}.
An infinitely divisible law without Gaussian components is a \textit{layered stable law} if its L\'evy measure satisfies
\begin{equation}\label{layed stable law Levy measure}
	\nu(B) = \int_{S^{d-1}}\int_0^\infty \mathbbm1_B(r\bm\xi)q(r,\bm\xi)\,dr\,\sigma(d\bm\xi),\quad B\in\mathscr{B}(\mathbb R_0^d),
\end{equation}
where $\sigma(d\bm\xi)$ is a probability measure on the unit sphere $S^{d-1}$ of $\mathbb R^d$, and $q:(0,\infty)\times S^{d-1}\to(0,\infty)$ is a locally integrable function such that
\[
	q(r,\bm\xi) \sim \begin{dcases}
		c_1(\bm\xi)r^{-\alpha-1},\quad& \text{as }r\to0,\\
		c_2(\bm\xi)r^{-\beta-1}, \quad& \text{as }r\to\infty,
	\end{dcases}
\]
for almost every $\bm\xi\in S^{d-1}$, where $c_1$ and $c_2$ are $\sigma$-integrable positive functions on $S^{d-1}$, where $(\alpha,\beta)\in(0,2)\times(0,\infty)$ are the inner and outer stability indices, respectively.
Similarly to the tempered stable process, the layered stable process exhibits transient behaviour across different time scales. Specifically, it behaves as a $\alpha$-stable process in short time, and as a $\beta$-stable process in long time. When $\beta>2$, the long time behaviour resembles a Brownian motion. Shot noise representations of the layered stable law are presented in Example \ref{example layered stable law}.

\section{Shot noise representation of infinitely divisible laws}\label{section series representation}
We work towards shot noise series representations of infinitely divisible laws without Gaussian components. That is, we seek to represent such laws using series with summands dependent on Poisson arrival times and random shot markings. We build up our understanding of series representations from the least general case to the most, with the most general formulation being the generalised shot noise method \cite{rosinski1990series} (Theorem \ref{theorem generalised shot noise method}). We motivate series representations of infinitely divisible laws via the L\'evy-It\^o decomposition \cite{resnick1976extremal,rosinski2001series}, and demonstrate that such representations are possible precisely because the laws can be summarised entirely by the jumps of the corresponding L\'evy process.

\subsection{The case of finite L\'evy measure}\label{section finite Levy measure}
As a motivating example, we begin by offering a simple shot noise representation for the infinitely divisible random variable characterised by a finite and absolutely continuous L\'evy measure with bounded positive support \cite{ferguson1972representation}. Recall that the case of a finite L\'evy measure corresponds to a compound Poisson distribution. Such a random variable can be thought of as the position of the corresponding compound Poisson process at unit time. Let $\{\Gamma_{k}\}_{k\in\mathbb N}$ be the arrival times of a standard Poisson process.
Consider a nonnegative intensity function $h$ such that $\int_0^T h(s)\,ds<+\infty$ for a fixed truncation time $T>0$, and define $H_I(t) \coloneqq \inf\{u\in[0,T]:\int_0^u h(s)\,ds<t\}$ as the generalised inverse of the corresponding mean value function.
Then, the random variable
\begin{equation}\label{series representation of ID laws}
	\sum_{k=1}^{+\infty} H_I(\Gamma_{k})\mathbbm1(\Gamma_{k}\in[0,\textstyle{\int}_0^Th(s)\,ds])
\end{equation}
is well-defined and infinitely divisible with L\'evy measure $h(z)dz$ defined on $(0,T]$.
This result can be verified by checking the characteristic function of the series, which is easily computed by conditioning on the number of jumps of the standard Poisson process over $[0,\int_0^T h(s)\,ds]$.

With the shot noise representation \eqref{series representation of ID laws}, we have represented every univariate infinitely divisible law with a finite and absolutely continuous L\'evy measure defined on $(0,T]$ as a series based on Poisson arrival times. Specifically, if $X$ is an infinitely divisible random variable in $\mathbb R$ without Gaussian components and with a finite L\'evy measure $\nu(dz) = h(z)\, dz$ whose density has support over a positive bounded interval $(0,T]$, then $X$ is equal in distribution to the series \eqref{series representation of ID laws}. So indeed, by setting $H_I(t)$ as the cumulative integral of the L\'evy measure $\nu(dz)$ and interpreting the Poisson arrival times as random steps over the domain of $H_I(t)$, a transformed state space, rather than over time, we obtain a shot noise representation of the infinitely divisible law without Gaussian components. This idea naturally leads the more general inverse L\'evy measure method discussed later in Section \ref{section inverse Levy measure method}. The subscript for the kernel $H_I$ is used to distinguish between the kernel of the inverse L\'evy measure method and that of the generalised shot noise method in Section \ref{section generalised shot noise method}.


\subsection{L\'evy-It\^o decomposition}\label{section Levy-Ito decomposition}
We would like to introduce the L\'evy-It\^o decomposition as a means to work towards a general approach for shot noise series representation, which extends to infinitely divisible random vectors with infinite L\'evy measures. Recall that we can interpret infinitely divisible random variates as the unit-time increment $X_1 -X_0$ of its corresponding L\'evy process. This interpretation will prove to be fruitful for more general shot noise representations for infinitely divisible laws through the L\'evy-It\^o decomposition for L\'evy processes, as first demonstrated by Resnick \cite{resnick1976extremal}.

For a fixed $T>0$, consider a L\'evy process $\{X_t:t\in[0,T]\}$ in $\mathbb R^d$ without Gaussian components and with a general L\'evy measure $\nu(d\mathbf z)$ satisfying \eqref{Levy measure condition}. Let us define $X_{0-} \coloneqq X_0$ and a measure $\mu(d\mathbf z,ds)$ on $\mathbb R_0^d\times[0,+\infty)$ such that
\[
	\mu(B,J) = |\{s\in J : X_s - X_{s-}\in B\}| = \int_J \mathbbm1_B(X_s - X_{s-})\, ds,\quad (B,J)\in\mathcal B(\mathbb R_0^d)\times\mathcal B([0,+\infty)).
\]
That is, $\mu(B,J)$ is a random variable that counts the number of jumps with sizes in $B$ during the time interval $J$.
Define $S_{t,n}\coloneqq\{s\in[0,t] : |X_s - X_{s-}|>1/n\}$ to be the set of points at which jumps with sizes greater than $1/n$ occur in $[0,t]$. Since $X$ is almost surely finite on $[0,t]$, we must have that $|S_{t,n}|<+\infty$ almost surely, so the set of all jumps $S_t\coloneqq\cup_{n\in\mathbb N} S_{t,n}$ over $[0,t]$ must be almost surely at most countable. Thus, we can express
\[
	\mu(B,[0,t]) = \sum_{s\in S_t}\delta_{(X_s - X_{s-})}(B)\quad \text{a.s.},
\]
for every $B\in\mathcal B(\mathbb R_0^d)$. Clearly, $\{\mu(\cdot,[0,t])\}_{t\ge0}$ forms a family of random counting measures. We state a couple of results from \cite{applebaum2009levy}. Firstly, if $\{N_t:t\ge0\}$ is a L\'evy process that is nondecreasing and $N_t - N_{t-}$ takes values in $\{0,1\}$ for every $t\ge0$, then it is a Poisson process. This naturally leads to the following result. Let $B\in\mathcal B(\mathbb R^d_0)$ such that $\nu(B)<+\infty$. Then, $\{\mu(B,[0,t]):t\ge0\}$ is a Poisson process in $\mathbb R^d_0$ with intensity $\nu(B)$.
Additionally, we can include the case of infinite L\'evy measures by replacing the assumption $\nu(B)<+\infty$ with $B$ not including zero in its closure, as this avoids the accumulation of infinitely many small jumps. The upshot is that $\mu(d\mathbf z,ds)$ is a Poisson random measure on $\mathbb R_0^d\times[0,+\infty)$ with intensity measure $(\nu\times m)(d\mathbf z,ds)$, where $m(ds)$ is the Lebesgue measure on $\mathbb R$. Note that this verifies the interpretation of $\nu(B)$ as the expected number of jumps with sizes in $B\in\mathcal B(\mathbb R^d_0)$ over the unit interval, as
\[
	\nu(B) = \nu(B)m([0,1]) = \mathbb E[\mu(B,[0,1])].
\]
Moreover, the stated result on counting jump discontinuities shows that every L\'evy process admits a Poisson random measure on $\mathbb R^d_0\times[0,+\infty)$ with intensity measure $(\nu\times m)(d\mathbf z,ds)$.
With this understanding of the connection between Poisson random measures and L\'evy processes, it is reasonable now to present the following case of the L\'evy-It\^o decomposition.

\begin{theorem}[L\'evy-It\^o decomposition]\label{theorem LI decomposition}
	Let $\{X_t:t\ge0\}$ be a L\'evy process in $\mathbb R^d$ without Gaussian components with L\'evy measure $\nu(d\mathbf z)$. Where $\mu(d\mathbf z,ds)$ is the Poisson random measure on $\mathbb R^d_0\times [0,+\infty)$ with intensity measure $(\nu\times m)(d\mathbf z,ds)$ associated with $X$, it holds that
	\begin{equation}\label{LI decomposition}
		X_t = \int_0^t \int_{\|{\bf z}\|\in (0,1]}{\mathbf z}(\mu (d{\mathbf z},ds)-\nu (d{\mathbf z})ds) + \int_0^t\int_{\|\mathbf z\|>1}\mathbf z\, \mu(d\mathbf z,ds),\quad t\ge0.
	\end{equation}
\end{theorem}

We refer the reader to \cite[Chapter 4]{sato1999levy} for details. 


One must take care when working with the double integral term in the L\'evy-It\^o decomposition; the linearity of the integral cannot be applied if the integrand $\mathbf z$ is not integrable with respect to the relevant measures. This highlights the importance of the compensation term $v(d\mathbf z)ds$.
In the case that $\mathbf z$ is integrable, linearity of the integral can be applied to \eqref{LI decomposition} to obtain
\[
	X_t = \int_0^t\int_{\mathbb R_0^d}\mathbf z\,\mu(d\mathbf z,ds) - \int_0^t\int_{\|\mathbf z\|\in(0,1]}\mathbf z\,\nu(d\mathbf z)ds = \int_{\mathbb R_0^d}\mathbf z\, \mu(d\mathbf z,[0,t]) - t\int_{\|\mathbf z\|\in(0,1]}\mathbf z\, \nu(d\mathbf z).
\]
In the case where $\mathbf z$ is not integrable, define
\[
	X_{t}^{(n)} \coloneqq \int_0^t \int_{\|{\bf z}\|\in[1/n,1]}{\mathbf z}(\mu (d{\mathbf z},ds)-\nu (d{\mathbf z})ds) + \int_{\|\mathbf z\|>1}\mathbf z\ \mu(d\mathbf z,[0,t]).
\]
Then, we are guaranteed that $\mathbf z$ is integrable for every $n\in\mathbb N$ and hence we can apply linearity of the integral to obtain
\[
	X_{t}^{(n)} = \int_{\|\mathbf z\|\ge1/n}\mathbf z\ \mu(d\mathbf z,[0,t]) - t\int_{\|\mathbf z\|\in[1/n,1]}\mathbf z\, \nu(d\mathbf z).
\]
It is clear that $X_{t}^{(n)}\to X_t$ pointwise as $n\to+\infty$. Thus, we can derive a series representation for L\'evy processes by finding an appropriate representation of $\mu(d\mathbf z,ds)$ as a random sum, which is possible in principle as it is a random counting measure. For the moment, suppose that $\mu(d\mathbf z,ds) = \sum_{k=1}^{+\infty}\delta_{(J_k,T_k)}(d\mathbf z,ds)$ almost surely, where $\{J_k\}_{k\in\mathbb N}$ and $\{T_k\}_{k\in\mathbb N}$ are suitable independent sequences of random vectors in $\mathbb R^d$ and $[0,T]$, respectively. Then, we have that
\[
	X_{t}^{(n)} = \sum_{k=1}^{+\infty}\int_{\|\mathbf z\|\ge\frac{1}{n}}\mathbf z\, \delta_{(J_k,T_k)}(d\mathbf z,[0,t]) - t\int_{\|\mathbf z\|\in[1/n,1]}\mathbf z\, \nu(d\mathbf z) = \sum_{k=1}^{+\infty} J_k\mathbbm1_{[1/n,+\infty)}(\|J_k\|)\mathbbm1_{[0,t]}(T_k) - t\int_{\|\mathbf z\|\in[1/n,1]}\mathbf z\, \nu(d\mathbf z) \quad \text{a.s.}
\]
Finally, by letting $n\to+\infty$, we obtain the shot noise series representation \cite[Section 1]{rosinski2001series}
\begin{equation}\label{series rep via Levy-Ito}
	X_t = \sum_{k=1}^{+\infty} (J_k\mathbbm1_{[0,t]}(T_k) - tc_k) \quad \text{a.s.,}
\end{equation}
where $\{c_k\}_{k\in\mathbb N}$ is a sequence of suitable compensation vectors in $\mathbb R^d$ (sometimes referred to as \textit{centres} if the expectation of each term is used) which guarantee the convergence of the series. In the case where the L\'evy process is a subordinator, convergence without compensation vectors is ensured, as its L\'evy measure has finite first moment about the origin.
Where $X_1$ is our infinitely divisible random vector of interest, its shot noise representation is given by $\sum_{k=1}^{+\infty} (J_k - c_k)$ in law.

As demonstrated, the crux of this series representation based upon the L\'evy-It\^o decomposition is the representation of the Poisson random measure associated with a L\'evy process as a random sum of Dirac delta measures. That is to say, this approach has reduced the problem of deriving the shot noise representation of L\'evy processes in general settings to finding expressions of $\mu(d\mathbf z,ds)$ as random sums. We explore one such method in the following.

\subsection{Inverse L\'evy measure method}\label{section inverse Levy measure method}
We would like to generalise the shot noise representation in Section \ref{section finite Levy measure} to infinitely divisible random variables with infinite L\'evy measures with support extending to negative real numbers. In Section \ref{section finite Levy measure}, we assumed $\int_0^T h(s)\, ds<+\infty$ so that the generalised inverse of the mean value function $H_I$ is well-defined. However, a simple extension to infinite L\'evy measures is to instead define the kernel $H_{I}(r) \coloneqq \inf\{u\in(0,+\infty) : \nu((u,+\infty]) < r\}$ to \textit{run down from infinity} instead.
This kernel is well-defined even in the case of infinite L\'evy measures since, by definition, the tail of the L\'evy measure is finite.
Let $\{X_t:t\in[0,1]\}$ be a L\'evy process in $(0,+\infty)$ without Gaussian components with L\'evy measure $\nu(dz)$. Unlike the setting of Section \ref{section finite Levy measure}, we do not assume the L\'evy measure $\nu(dz)$ to admit a density on $(0,+\infty)$. As mentioned in the previous section,
our goal is to represent the Poisson random measure $\mu(dz,ds)$ with intensity measure $\nu\times m$ associated with the L\'evy process $X$ as a random sum.
We state the following result from \cite[Proposition 2.1]{rosinski2001series}, which will shortly prove to be essential.

\begin{proposition}\label{prop transformation of Poisson random measure}
	Let $M(ds)$ and $N(dz)$ be Poisson random measures on Borel spaces $S$ and $T$ with intensity measures $\mu$ and $\nu$, respectively. If there exists a measurable function $h:S\to T$ such that $\nu=\mu\circ h^{-1}$ on $\mathcal B(T)$, then $N(dz)$ is equal in law to $M\circ h^{-1}(dz)$. If, in addition, the Poisson random measure $N(dz)$ is defined on a probability space which admits the existence of a standard uniform random variable independent of $N(dz)$, and there exists a sequence of random elements $\{Y_k\}_{k\in\mathbb N}$ on $S$ such that
	\[
		M(ds)=\sum_{k=1}^{+\infty}\delta_{Y_k}(ds),
	\]
	then there exists a sequence of random elements $\{X_k\}_{k\in\mathbb N}$ defined on a common probability space as $N(dz)$ that is identical in law to $\{Y_k\}_{k\in\mathbb N}$ and
	\[
		N(dz) = \sum_{k=1}^{+\infty} \delta_{h(X_k)}(dz) \quad \text{a.s.}
	\]
\end{proposition}

With this, we will now present the \textit{inverse L\'evy measure} method for computing a shot noise representation in the one-dimensional setting. Let $\{\Gamma_k\}_{k\in\mathbb N}$ be a sequence of standard Poisson arrival times. Then, its corresponding Poisson random measure with intensity measure $m(dz)$ can be expressed by $\sum\nolimits_{k=1}^{+\infty} \delta_{\Gamma_k}(dz)$. Define the marked Poisson random measure
\[
	M(dz,ds)\coloneqq\sum_{k=1}^{+\infty} \delta_{(\Gamma_k,T_k)} (dz,ds),
\]
where $\{T_k\}_{k\in\mathbb N}$ are iid standard uniform random variables independent of $\{\Gamma_k\}_{k\in\mathbb N}$. Then, it holds that $M(dz,ds)$ is a Poisson random measure on $(0,+\infty)\times[0,1]$ with intensity measure $m\times m$. Define $H_I^*:(0,+\infty)\times[0,1]\to(0,+\infty)\times[0,1]$ by $(y,u)\to(H_I(y),u)$, then $(m\times m)\circ (H_I^*)^{-1} = \nu\times m$. By Proposition \ref{prop transformation of Poisson random measure}, it holds that
\begin{equation}\label{Poisson random measure as series}
	\mu(dz,ds)\stackrel{\mathscr L}{=}\sum_{k=1}^{+\infty} \delta_{(H_I(\Gamma_k),T_k)}(dz,ds).
\end{equation}
Substituting this representation for $\mu(dz,ds)$ into the L\'evy-It\^o decomposition of $X_t$ and following the derivation of \eqref{series rep via Levy-Ito}, we obtain a shot noise representation for $X$ via the inverse L\'evy measure method as \cite{ferguson1972representation,rosinski2001series}
\[
	\{X_t:t\in[0,1]\} \stackrel{\mathscr L}{=} \left\lbrace\sum_{k=1}^{+\infty}( H_I(\Gamma_k)\mathbbm{1}_{[0,t]}(T_k) - tc_k) : t\in[0,1]\right\rbrace,
\]
where $\{c_k\}_{k\in\mathbb N}$ is a sequence of suitable centres. Where the unit-time marginal $X_1$ is the infinitely divisible random variable of our interest, its shot noise representation is given by $\sum_{k=1}^{+\infty}( H_I(\Gamma_k) - c_k)$ in law. This is almost identical to the series \eqref{series representation of ID laws} except for the appearance of compensating constants, which should be expected given the possibility of heavy intensity of the L\'evy measure about the origin.

We wish to extend the inverse L\'evy measure method to the multidimensional setting with jumps in any direction. We will loosely present LePage's approach \cite{lepage1980multidimensional}. Let $\{X_t:t\in[0,1]\}$ now be a L\'evy process in $\mathbb R^d$ without Gaussian components with L\'evy measure $\nu(d\mathbf z)$. We seek a representation of the associated Poisson random measure $\mu(d\mathbf z,ds)$ on $\mathbb R^d_0\times[0,1]$ with intensity measure $\nu\times m$ as a series like \eqref{Poisson random measure as series}. Consider a radial disintegration of $\nu(d\mathbf z)$ given by
\begin{equation}\label{radial decomposition of Levy measure}
	\nu(B) = \int_{S^{d-1}}\int_{0+}^{+\infty}\mathbbm1_B (r\bm\xi)\, h(dr,\bm\xi)\, \sigma(d\bm\xi),\qquad B\in\mathcal B(\mathbb R_0^d),
\end{equation}
where $\sigma(d\bm\xi)$ is some probability measure on the unit sphere $S^{d-1}$ of $\mathbb R^d$ and $\{h(\cdot,\bm\xi)\}_{\bm\xi\in S^{d-1}}$ is a measurable family of L\'evy measures on $(0,+\infty)$. The idea behind this disintegration is to decompose the intensities of jumps sizes by magnitude and direction; specifically, the intensity of jumps of magnitude $r>0$ in the direction $\bm\xi\in S^{d-1}$ is represented by $h(dr,\bm\xi)$, while the proportion of all jumps with directions in the set $C\subseteq S^{d-1}$ is given by $\sigma(C)$. We define the generalised inverse of the tail of $h((u,+\infty),\bm\xi)$ in the direction $\bm\xi$ as
\begin{equation}\label{inverse Levy measure kernel}
	H_I(r,\bm\xi)\coloneqq\inf\{u\in(0,+\infty): h((u,+\infty),\bm\xi)<r\}.
\end{equation}
Let $\{U_k\}_{k\in\mathbb N}$ be a sequence of iid random vectors with law $\sigma(d\bm\xi)$ on $S^{d-1}$, independent of $\{\Gamma_k\}_{k\in\mathbb N}$ and $\{T_k\}_{k\in\mathbb N}$. Similarly to before, we define the marked Poisson process $\sum\nolimits_{k=1}^{+\infty} \delta_{(\Gamma_k,U_k,T_k)}(dr,d\bm\xi,ds)$ on $(0,+\infty)\times S^{d-1}\times[0,1]$ with intensity measure $m\times \sigma\times m$. Define $H_I^*:(0,+\infty)\times S^{d-1}\times[0,1]\to \mathbb R^d_0\times[0,1]$ such that $(r,\bm\xi,s) \mapsto (H_I(r,\bm\xi)\bm\xi,s)$. Due to the radial representation \eqref{radial decomposition of Levy measure} of the L\'evy measure $v(d\mathbf z)$, we have that $(m\times\sigma\times m)\circ (H_I^*)^{-1} = \nu\times m$. By Proposition \ref{prop transformation of Poisson random measure}, it holds that \cite[Section 3]{rosinski2001series}
\[
	\mu(d\mathbf z,ds)\stackrel{\mathscr L}{=}\sum_{k=1}^{+\infty} \delta_{( H_I(\Gamma_k,U_k)U_k,T_k)}(d\mathbf z,ds).
\]
Substituting the expression for $\mu(d\mathbf z,ds)$ into the L\'evy-It\^o decomposition, we obtain
\begin{equation}\label{inverse Levy measure method series over [0,1]}
	\{X_t:t\in[0,1]\} \stackrel{\mathscr L}{=} \left\lbrace\sum_{k=1}^{+\infty} (H_I(\Gamma_k,U_k)U_k \mathbbm1_{[0,t]}(T_k) - tc_k):t\in[0,1]\right\rbrace.
\end{equation}
Where the unit-time marginal $X_1$ is our infinitely divisible random vector of interest, its shot noise representation is given by $\sum_{k=1}^{+\infty} (H_I(\Gamma_k,U_k)U_k - c_k)$ in law. By comparing the shot noise representations between infinitely divisible laws and their corresponding L\'evy measures, we see a correspondence in which the series pertaining to the latter is merely the former but with uniform scattering of summands. Intuitively, the uniform scattering is necessary to preserve the stationarity of increments. An alternative derivation of the inverse L\'evy measure method in the multidimensional case can be obtained via the generalised shot noise method of Section \ref{section generalised shot noise method} in the following.

We provide the shot noise representation of the stable law obtained via the inverse L\'evy measure method \cite{lepage1981convergence}, which is crucial from a theoretical perspective \cite{davydov2012convergence,janicki1993simulation,samorodnitsky1991construction,samorodnitsky1994stable} as well as for practical use \cite{kawai2016higher,kawai2017sample,lemke2015fully}, to name a few examples.

\begin{example}[Inverse L\'evy measure method for stable random vector]\label{example stable series representation}
	Let $X$ be a stable law with L\'evy measure \eqref{stable Levy measure}. Then, by the inverse L\'evy measure method, it holds that
	\begin{equation}\label{series rep of stable processes}
		X \stackrel{\mathscr L}{=} \sum_{k=1}^{+\infty}\left(\left(\frac{\alpha\Gamma_k}{a}\right)^{-1/\alpha}U_k - c_k\right),
	\end{equation}
	where $\{U_k\}_{k\in\mathbb N}$ is a sequence of iid random vectors with distribution $\sigma(d\bm\xi)$ and $\{c_k\}_{k\in\mathbb N}$ is a sequence of suitable centres given in Theorem \ref{theorem generalised shot noise method}. If $\sigma(d\bm\xi)$ is isotropic or the L\'evy measure $\nu(d\mathbf z)$ satisfies $\int_{\mathbb R_0^d}(1\wedge \|\mathbf z\|)\,\nu(d\mathbf z) < +\infty$, then the infinite series \eqref{series rep of stable processes} converges almost surely without the centres $\{c_k\}_{k\in\mathbb N}$. $\hfill\qed$
\end{example}

\subsection{Generalised shot noise method}\label{section generalised shot noise method}
We will now present the generalised shot noise method \cite{rosinski1990series,rosinski2001series}, which generalises the inverse L\'evy measure method of Section \ref{section inverse Levy measure method}.

\begin{theorem}[Generalised shot noise method]\label{theorem generalised shot noise method}
	Suppose a L\'evy measure $\nu(d\mathbf z)$ on $\mathbb R^d_0$ can be decomposed as
	\begin{equation}\label{generalised shot noise method decomposition}
		\nu(B) = \int_{0+}^{+\infty} \mathbb P(H(r,U)\in B)\, dr, \quad B\in\mathcal B(\mathbb R^d_0),
	\end{equation}
	where $U$ is a random vector in some space $\mathscr U$ and $H:(0,+\infty)\times\mathscr U\to\mathbb R^d_0$ is such that for every $u\in\mathscr U$, $r\mapsto \|H(r,u)\|$ is nonincreasing. Then, the following statements hold.
	\begin{enumerate}[(i)]
		\setlength{\parskip}{0cm}
		\setlength{\itemsep}{0cm}
		\item It holds that
		\begin{equation}\label{generalised shot noise series}
			X\stackrel{\mathscr L}{=}\sum_{k=1}^{+\infty} \left(H\left(\Gamma_{k},U_k\right)-c_k\right),
		\end{equation}
		where $X$ is an infinitely divisible random vector without Gaussian components with the L\'evy measure $\nu(d\mathbf z)$, $\{\Gamma_{k}\}_{k\in\mathbb N}$ are the arrival times of the standard Poisson process, $\{U_k\}_{k\in\mathbb N}$ are iid copies of $U$ independent of $\{\Gamma_k\}_{k\in\mathbb N}$, and $\{c_k\}_{k\in\mathbb N}$ is a sequence of suitable centres in $\mathbb R^d$. Moreover, we can take
		\[
			c_k \coloneqq \int_{k-1}^k \mathbb E\left[H(s,U)\mathbbm1_{(0,1]}(\|H(s,U)\|)\right]\,ds,\quad k\in\mathbb N.
		\]
		If additionally, the L\'evy measure $\nu(d\mathbf z)$ on $\mathbb R^d_0$ satisfies $\int_{\|\mathbf z\|>1}\|\mathbf z\|\,\nu(d\mathbf z)<+\infty$, then we can instead take $c_k = \int_{k-1}^k \mathbb E[H(s,U)]\,ds$ for every $k\in\mathbb N$.
		
		\item If $a\coloneqq \lim_{s\to+\infty}\int_{0+}^s\int_{\|\mathbf z\|\le1}\mathbf z\,\mathbb P(H(r,U)\in d\mathbf z)\,dr$ exists in $\mathbb R^d$, then it holds that $\sum_{k=1}^{+\infty} H\left(\Gamma_{k},U_k\right)$ is an infinitely divisible random vector characterised by the L\'evy-Khintchine triplet $(a,0,\nu)$.
	\end{enumerate}
\end{theorem}

A rigorous proof of the generalised shot noise method can be found in \cite{rosinski2001series}, along with the almost sure convergence of the series \eqref{generalised shot noise series}. Moreover, there exist independent random sequences $\{\Gamma_k\}_{k\in\mathbb N}$ and $\{U_k\}_{k\in\mathbb N}$ such that the equality \eqref{generalised shot noise series} holds almost surely. This result not only generalises the inverse L\'evy measure method, but as the decomposition \eqref{generalised shot noise method decomposition} of the L\'evy measure is not unique, the generalised shot noise method can be used to derive distinct shot noise representations of the same infinitely divisible random vector. In particular, Theorem \ref{theorem generalised shot noise method} can be used to derive the inverse L\'evy measure, rejection, thinning and Bondesson's methods for shot noise representation of an infinitely divisible random vector, as follows.

\begin{proposition}[Inverse L\'evy measure, rejection, thinning and Bondesson's methods] \label{prop series representations}
	Let $\nu(d\mathbf z)$ be a L\'evy measure on $\mathbb R_0^d$ such that
	\begin{equation}\label{radial decomposition}
		\nu(B) = \int_{S^{d-1}}\int_{0+}^{+\infty}\mathbbm 1_B(r\bm\xi)h(dr,\bm\xi)\sigma(d\bm\xi), \quad B\in\mathcal B(\mathbb R_0^d),
	\end{equation}
	where $\sigma(d\bm\xi)$ is a probability measure on the unit sphere $S^{d-1}$ of $\mathbb R^d$ and $\{h(\cdot,\bm\xi)\}_{\bm\xi\in S^{d-1}}$ is a measurable family of L\'evy measures on $(0,+\infty)$. Let $\{\Gamma_k\}_{k\in\mathbb N}$ be a sequence of standard Poisson arrival times, and $\{U_k\}_{k\in\mathbb N}$ a sequence of iid random vectors under $\sigma(d\bm\xi)$ independent of $\{\Gamma_k\}_{k\in\mathbb N}$. Then, an infinitely divisible random vector $X$ without Gaussian components with the L\'evy measure $\nu(d\mathbf z)$ has the following representations.
	\begin{enumerate}[(i)]
		\setlength{\parskip}{0cm}
		\setlength{\itemsep}{0cm}
		\item \textit{Inverse L\'evy measure method} \cite{ferguson1972representation,lepage1980multidimensional}: Define
		\begin{equation}\label{inverse Levy measure method kernel}
			H_I(r,\bm\xi) \coloneqq \inf\{u\in(0,+\infty) : h((u,+\infty),\bm\xi) < r\},\quad (r,\bm\xi)\in(0,+\infty)\times S^{d-1}.
		\end{equation}
		Then, it holds that
		\[
			\nu(B) = \int_{0+}^{+\infty}\mathbb P(H_I(r,U)U\in B)\,dr,\quad B\in\mathcal B(\mathbb R_0^d),
		\]
		where $U$ distributed under $\sigma(d\bm\xi)$. Hence, by Theorem \ref{theorem generalised shot noise method}, we have
		\begin{equation}\label{inverse Levy measure method series}
			X\stackrel{\mathscr L}{=}\sum_{k=1}^{+\infty} \left(H_I\left(\Gamma_k,U_k\right)U_k - c_{1,k}\right),
		\end{equation}
		where $\{c_{1,k}\}_{k\in\mathbb N}$ is a sequence of suitable centres.
		(Note that \cite{ferguson1972representation} only establishes the univariate case.)
		
		\item \textit{Rejection method} \cite{rosinski2001series}: Let $\nu_p$ be a L\'evy measure on $\mathbb R_0^d$ such that
		\[
			\nu_p(B) = \int_{S^{d-1}}\int_{0+}^{+\infty}\mathbbm1_B(r\bm\xi)h_p(dr,\bm\xi)\sigma(d\bm\xi),\quad B\in\mathcal B(\mathbb R_0^d),
		\]
		where $\{h_p(\cdot,\bm\xi)\}_{\bm\xi\in S^{d-1}}$ is a measurable family of L\'evy measures on $(0,+\infty)$ such that for every $\bm\xi\in S^{d-1}$, $h(\cdot,\bm\xi)$ is absolutely continuous with respect to $h_p(\cdot,\bm\xi)$ and satisfies $(dh/dh_p)\le1$. Define
		\[
			H_p(r,\bm\xi) \coloneqq \inf\{u\in(0,+\infty) : h_p((u,+\infty),\bm\xi)<r\},\quad (r,\bm\xi)\in(0,+\infty)\times S^{d-1}.
		\]
		Then, it holds that
		\[
			\nu(B) = \int_{0+}^{+\infty}\mathbb P\left(H_p(r,U)\mathbbm1\left(\frac{dh}{dh_p}(H_p(r,U),U) > V\right)U\in B\right)\,dr,\quad B\in\mathcal B(\mathbb R_0^d),
		\]
		where $U$ distributed under $\sigma(d\bm\xi)$ and independent of $V$, a standard uniform random variable. Hence, by Theorem \ref{theorem generalised shot noise method}, we have
		\begin{equation}\label{rejection method series}
			X\stackrel{\mathscr L}{=} \sum_{k=1}^{+\infty} \left(H_p\left(\Gamma_k,U_k\right)\mathbbm1\left(\frac{dh}{dh_p}\left(H_p\left(\Gamma_k,U_k\right),U_k\right) > V_k\right)U_k - c_{2,k}\right),
		\end{equation}
		where $\{V_k\}_{k\in\mathbb N}$ is a sequence of iid standard uniform random variables and $\{c_{2,k}\}_{k\in\mathbb N}$ is a sequence of suitable centres.
		
		\item \textit{Thinning method} \cite{rosinski1990series}: Suppose $\{F(\cdot,\bm\xi)\}_{\bm\xi\in S^{d-1}}$ is a family of probability measures on $(0,+\infty)$ such that for every $\bm\xi\in S^{d-1}$, $h(\cdot,\bm\xi)$ is absolutely continuous with respect to $F(\cdot,\bm\xi)$. Then, it holds that
		\[
			\nu(B) = \int_{0+}^{+\infty}\mathbb P\left(V\mathbbm1\left(\frac{dh}{dF}(V,U) > r\right)\in B\right)\,dr,\quad B\in\mathcal B(\mathbb R_0^d),
		\]
		where $(V,U)$ is distributed under $F(dv,\bm\xi)\times\sigma(d\bm\xi)$. Hence, by Theorem \ref{theorem generalised shot noise method}, we have
		\begin{equation}\label{thinning method series}
			X\stackrel{\mathscr L}{=} \sum_{k=1}^{+\infty}\left(V_k\mathbbm1\left(\frac{dh}{dF}(V_k,U_k)>\Gamma_k\right)U_k - c_{3,k} \right),
		\end{equation}
		where $\{(V_k,U_k)\}_{k\in\mathbb N}$ is a sequence of iid copies of $(V,U)$ and $\{c_{3,k}\}_{k\in\mathbb N}$ is a sequence of suitable centres.
		
		\item \textit{Bondesson's method} \cite{bondesson1982simulation}: Suppose that the L\'evy measure $h(\cdot,\bm\xi)$ can be decomposed as
		\[
			h(B,\bm\xi) = \int_{0+}^{+\infty} G(B/g(r,\bm\xi),\bm\xi)\,dr, \quad (B,\bm\xi)\in\mathcal B(\mathbb R_0^d)\times S^{d-1},
		\]
		where $\{G(\cdot,\bm\xi)\}_{\bm\xi\in S^{d-1}}$ is a family of probability measures on $(0,+\infty)$ and $\{g(\cdot,\bm\xi)\}_{\bm\xi\in S^{d-1}}$ is a family of nonincreasing functions from $(0,+\infty)$ to $(0,+\infty)$. Then, it holds that
		\[
			\nu(B) = \int_{0+}^{+\infty}\mathbb P(g(r,U)VU\in B)\,dr,\quad B\in\mathcal B(\mathbb R_0^d),
		\]
		where $(V,U)$ is distributed under $G(dv,\bm\xi)\times\sigma(d\bm\xi)$. Hence, by Theorem \ref{theorem generalised shot noise method}, we have
		\begin{equation}\label{Bondesson's method series}
			X \stackrel{\mathscr{L}}{=} \sum_{k=1}^{+\infty} \left(g\left(\Gamma_k,U_k\right)V_k U_k - c_{4,k}\right),
		\end{equation}
		where $\{(V_k,U_k)\}_{k\in\mathbb N}$ is a sequence of iid copies of $(V,U)$ and $\{c_{4,k}\}_{k\in\mathbb N}$ is a sequence of suitable centres.
	\end{enumerate}
\end{proposition}

Thus, we see that the generalised shot noise method allows us to choose from different shot noise representations of the same infinitely divisible law due to the nonuniqueness of the decomposition \eqref{generalised shot noise method decomposition} of the L\'evy measure. We remark that alternative series representations of L\'evy processes exist. One such example is via the Karhunen-Lo\`eve expansion \cite{hackmann2018karhunen}, leading to a Fourier-like series with infinitely divisible coefficients, which can be approximated via the truncation of shot noise representation (Section \ref{section truncation}).

\subsection{Examples}\label{section examples of series representation}
Armed with several different shot noise representation methods, we provide examples of shot noise representations of infinitely divisible laws. We begin by revisiting Example \ref{example stable series representation} of the stable law.

\begin{example}[Shot noise representations for the stable law]
	The shot noise representation \eqref{series rep of stable processes} can also be obtained via Bondesson's method, for instance, with $G\equiv\delta_{\{+1\}}$ and $h(r,\bm\xi) = a/(\alpha r^\alpha)$, as well as with the rejection method with the trivial choice $h_p(\cdot)\equiv h(\cdot)$. While the thinning method can be applied, with $F(dr,\bm\xi)=e^{-r}\,dr$ for example, it is significantly elapsed by the representation \eqref{series rep of stable processes} in terms of elegance and practicality. $\hfill\qed$
\end{example}

Next, in contrast to this lack of choice regarding shot noise representations for the stable law, we consider several shot noise representations for the tempered stable law.

\begin{example}[Shot noise representations of tempered stable law]\label{example tempered stable law}	
	We first consider shot noise representations of the tempered stable law based on the thinning, rejection and inverse L\'evy measure methods which we have seen previously in Proposition \ref{prop series representations} \cite{imai2011finite}. For every $\alpha\in(0,2)$ and $r>0$, define
	\begin{equation}\label{eq tempered stable law kernel}
		H_\alpha(r) \coloneqq \left(\frac{\alpha r}{m_{\alpha,\rho}}\right)^{-1/\alpha},
	\end{equation}
	where $m_{\alpha,\rho}\coloneqq \int_{\mathbb R_0^d}\|\mathbf v\|^\alpha\,\rho(d\mathbf v)$ and $\rho(d\mathbf v)$ is defined as in \eqref{tempered stable law rho equation}. Define the following kernels:
	\begin{align}
		\notag & H_1(r,\mathbf v) \coloneqq\inf\left\lbrace u\in(0,+\infty) : \int_u^{+\infty} m_{\alpha,\rho}\frac{e^{-s/\|\mathbf v\|}}{s^{\alpha+1}}\,ds>r\right\rbrace\frac{\mathbf v}{\|\mathbf v\|}\\
		\notag & H_2(r,w,\mathbf v) \coloneqq w\mathbbm1\left(r\le m_{\alpha,\rho}\frac{\|\mathbf v\|}{\lambda}\frac{\exp\left[-\frac{1-\lambda}{\|\mathbf v\|}w\right]}{w^{\alpha + 1}}\right)\frac{\mathbf v}{\|\mathbf v\|},\\
		\notag & H_3(r,w,\mathbf v) \coloneqq w\mathbbm1\left(r\le m_{\alpha,\rho}\Gamma(\lambda_1)\left(\frac{\|\mathbf v\|}{\lambda_2}\right)^{\lambda_1}\frac{\exp\left[-\frac{1-\lambda_2}{\|\mathbf v\|}w\right]}{w^{\alpha + \lambda_1}}\right)\frac{\mathbf v}{\|\mathbf v\|},\\
		\notag & H_4(r,u,\mathbf v) \coloneqq H_\alpha(r)\mathbbm1\left(e^{- H_\alpha(r)/\|\mathbf v\|} > u\right)\frac{\mathbf v}{\|\mathbf v\|},
	\end{align}
	where $\lambda,\lambda_2\in(0,1]$ and $\lambda_1>0$. Let $\{W_k^{(2)}\}_{k\in\mathbb N}$ be a sequence of independent exponential random variables with rate $\lambda/\|V_k\|$, and $\{W_k^{(3)}\}_{k\in\mathbb N}$ a sequence of independent gamma random variables with shape $\lambda_1$ and scale $\|V_k\|/\lambda_2$. Define the series
	\begin{align}
		&X_1 \coloneqq \sum_{k=1}^{+\infty} (H_1(\Gamma_k,V_k) - c_{1,k}), \label{inverse Levy measure method for TS law}\\
		&X_2 \coloneqq \sum_{k=1}^{+\infty} \left(H_2\left(\Gamma_k,W_k^{(2)},V_k\right) - c_{2,k}\right), \label{thinning method for TS law 1}\\
		&X_3 \coloneqq \sum_{k=1}^{+\infty} \left(H_3\left(\Gamma_k,W_k^{(3)},V_k\right) - c_{3,k}\right), \label{thinning method for TS law 2}\\
		&X_4 \coloneqq \sum_{k=1}^{+\infty} (H_4(\Gamma_k,U_k,V_k)-c_{4,k}), \label{rejection method for TS law}
	\end{align}
	where $\{c_{j,k}\}_{k\in\mathbb N}$, $j=1,2,3,4$, are sequences of suitable centres in $\mathbb R^d$. Then, for $j=1,2,3,4$, the series $X_{j}$ converges almost surely and is equal in law to the tempered stable law with parameters $\alpha$ and $\rho(d\mathbf v)$. The first two representations are derived from the thinning method and the latter two are derived from the rejection and the inverse L\'evy measure methods, respectively.
	
	Yet another shot noise representation for the tempered stable law is Rosi\'nski's representation \cite{rosinski2007tempering}, which can be verified via the generalised shot noise method (Theorem \ref{theorem generalised shot noise method}) but falls outside of the methods described in Proposition \ref{prop series representations}.
	Let $\{W_k\}_{k\in\mathbb N}$, $\{U_k\}_{k\in\mathbb N}$ and $\{V_k\}_{k\in\mathbb N}$ be mutually independent sequences of iid standard exponential, standard uniform random variables and random vectors in $\mathbb R_0^d$ with distribution $\|\mathbf v\|^\alpha\rho(d\mathbf v)/m_{\alpha,\rho}$, respectively.
	Then, where $\{\Gamma_k\}_{k\in\mathbb N}$ is a sequence of standard Poisson arrival times and $k_0$ and $z_0$ are suitable constants depending only on $\alpha$ and $\rho$, we have that the series
	\begin{equation}
		\sum_{k=1}^{+\infty} \left[\left(H_\alpha(\Gamma_k) \wedge W_k U_k^{1/\alpha}\|V_k\|\right)\frac{V_k}{\|V_k\|} - H_\alpha(k)k_0\right] + z_0\label{eq tempered stable law series representation}
	\end{equation}
	converges almost surely and is equal in law to the tempered stable law with parameters $\alpha$ and $\rho(d\mathbf v)$.
	Along with the aforementioned representations \eqref{inverse Levy measure method for TS law}--\eqref{rejection method for TS law}, the shot noise representation \eqref{eq tempered stable law series representation} shares the advantages of being explicit and exact. $\hfill\qed$
\end{example}

Note that the shot noise representations for the tempered stable law in Example \ref{example tempered stable law} with $\rho(dv) = C(G^Y\delta_{\{-1/G\}}(dv) + M^Y\delta_{\{1/M\}}(dv))$ do not cover the case of the CGMY law with $Y\le0$.
With $Y<0$, the CGMY law becomes compound Poisson. For the case when $Y=0$, which leads to the gamma law, we present its shot noise representations \cite[Section 6]{rosinski2001series} as follows.

\begin{example}[Shot noise representations of gamma law]\label{example gamma law}
	Let $X$ be a gamma law with L\'evy measure \eqref{gamma Levy measure}. Let $\{\Gamma_k\}_{k\in\mathbb N}$ be standard Poisson arrival times. Then, it holds that
	\begin{enumerate}[(i)]
		\setlength{\parskip}{0cm}
		\setlength{\itemsep}{0cm}
		\item by the inverse L\'evy measure method,
		\begin{equation}\label{inverse Levy measure method for gamma law}
			X \stackrel{\mathscr L}{=} \frac{1}{\beta}\sum_{k=1}^{+\infty} E_1^{-1}(\Gamma_k/a),
		\end{equation}
		where $E_1(x) \coloneqq \int_x^{+\infty}u^{-1}e^{-u}\,du$ denotes the exponential integral function and $E_1^{-1}$ its inverse.
		
		\item by the rejection method with $h(r) = ae^{-\beta r}/r$ and $h_p(r) = a/r(1+\beta r)$,
		\begin{equation}\label{rejection method for gamma law}
			X \stackrel{\mathscr L}{=} \frac{1}{\beta}\sum_{k=1}^{+\infty} \frac{1}{e^{\Gamma_k/a} - 1} \mathbbm1\left(\frac{e^{\Gamma_k/a}}{e^{\Gamma_k/a} - 1}\exp\left[-\frac{1}{e^{\Gamma_k/a}-1}\right] > V_k\right),
		\end{equation}
		where $\{V_k\}_{k\in\mathbb N}$ is a sequence of iid standard uniform random variables.
		
		\item by the thinning method with $F(dr) = \beta e^{-\beta r}\,dr$,
		\begin{equation}\label{thinning method for gamma law}
			X \stackrel{\mathscr L}{=} \frac{1}{\beta}\sum_{k=1}^{+\infty} V_k\mathbbm1\left(\Gamma_k V_k < a\right),
		\end{equation}
		where $\{V_k\}_{k\in\mathbb N}$ is a sequence of iid standard exponential random variables.
		
		\item by Bondesson's method with $G(du) = \beta e^{-\beta u}\,du$ and $g(r) = e^{-r/a}$,
		\begin{equation}\label{Bondesson's method for gamma law}
			X \stackrel{\mathscr L}{=} \frac{1}{\beta}\sum_{k=1}^{+\infty} e^{-\Gamma_k/a}V_k,
		\end{equation}
		where $\{V_k\}_{k\in\mathbb N}$ is a sequence of iid standard exponential random variables. $\hfill\qed$
	\end{enumerate}
\end{example}

Of the shot noise representations of the gamma law above, the easiest series to work with is perhaps the one associated with Bondesson's method \eqref{Bondesson's method for gamma law}, first appearing in \cite{bondesson1982simulation}. By contrast, the most difficult series from an implementation point of view is the one resulting from the inverse L\'evy measure method. Comparing this to the case of the stable law (Example \ref{example stable series representation}) where we saw that the inverse L\'evy measure method yields the most convenient representation, it is clear that there is an advantage to having a variety of shot noise representation methods such as those in Theorem \ref{theorem generalised shot noise method} and Proposition \ref{prop series representations}.
We mention that the inverse L\'evy measure method \eqref{inverse Levy measure method for gamma law} is employed in \cite{kotz2001laplace} to describe shot noise representations for the variance gamma law and process.
Next, we present shot noise representations of a layered stable law \cite{houdre2007layered}.

\begin{example}[Shot noise representations of a layered stable law]\label{example layered stable law}
	Let $X$ be a layered stable law with the L\'evy measure \eqref{layed stable law Levy measure} and
	\[
		q(r,\bm\xi) = r^{-\alpha-1}\mathbbm1_{(0,1]}(r) + r^{-\beta-1}\mathbbm1_{(1,\infty)}(r).
	\]
	Let $\{\Gamma_k\}_{k\in\mathbb N}$ be standard Poisson arrival times independent of a sequence $\{V_k\}_{k\in\mathbb N}$ of iid random vectors distributed under $\sigma(d\bm\xi)$. Denote $z_0 \coloneqq \mathbb E[V_1]$ and $\{b_k\}_{k\in\mathbb N}$ a sequence of suitable centring constants depending only on $\beta$. Then, it holds that
	\begin{enumerate}[(i)]
		\setlength{\parskip}{0cm}
		\setlength{\itemsep}{0cm}
		\item by the inverse L\'evy measure method,
		\[
			X \stackrel{\mathscr{L}}{=} \sum_{k=1}^{+\infty}\left[\left(\left(\beta\Gamma_k\right)^{-1/\beta}\mathbbm1_{(0,1/\beta]}(\Gamma_k) + \left(\alpha\Gamma_k + 1 - \frac\alpha\beta\right)^{-1/\alpha}\mathbbm1_{(1/\beta,+\infty)}(\Gamma_k)\right)V_k - b_k z_0\right].
		\]
		
		\item assuming $\alpha < \beta$, by the rejection method,
		\[
			X \stackrel{\mathscr{L}}{=} \sum_{k=1}^{+\infty}\left[(\alpha\Gamma_k)^{-1/\alpha}\mathbbm1\left(\frac{d\nu_{\sigma,q}^{\alpha,\beta}}{d\nu_\sigma^\alpha}\left((\alpha \Gamma_k)^{-1/\alpha}V_k\right) \ge U_k\right)V_k - b_kz_0\right],
		\]
		where $\{U_k\}_{k\in\mathbb N}$ is a sequence of iid standard uniform random variables independent of the other random sequences and
		\[
			\frac{d\nu_{\sigma,q}^{\alpha,\beta}}{d\nu_\sigma^\alpha}(\mathbf z) = \mathbbm1_{(0,1]}(\|\mathbf z\|) + \|\mathbf z\|^{\alpha-\beta}\mathbbm1_{(1,+\infty)}(\|\mathbf z\|) \le 1.\]
		
		\item assuming $\alpha < \beta$, by the rejection method,
		\[
			X \stackrel{\mathscr{L}}{=} \sum_{k=1}^{+\infty}\left[(\beta\Gamma_k)^{-1/\beta}\mathbbm1\left(\frac{d\nu_{\sigma,q}^{\alpha,\beta}}{d\nu_\sigma^\beta}\left((\beta \Gamma_k)^{-1/\beta}V_k\right) \ge U_k\right)V_k - b_kz_0\right],
		\]
		where $\{U_k\}_{k\in\mathbb N}$ is a sequence of iid standard uniform random variables independent of the other random sequences and
		\[
			\frac{d\nu_{\sigma,q}^{\alpha,\beta}}{d\nu_\sigma^\beta}(\mathbf z) = \|\mathbf z\|^{\beta-\alpha}\mathbbm1_{(0,1]}(\|\mathbf z\|) + \mathbbm1_{(1,+\infty)}(\|\mathbf z\|) \le 1.\hfill\qed
		\]
	\end{enumerate}
\end{example}

These examples of shot noise representations are only a handful of the applications of shot noise methods in the literature, and the derivation of shot noise representations and their usage in sampling and simulation are still ongoing topics of research. We remark that shot noise series is the only known representation of infinitely divisible laws in many cases, with the multivariate stable law and its generalisations as such examples. In the following, we discuss a truncation scheme for sampling via shot noise representations and the associated error analysis.

\section{Truncation of shot noise representations}\label{section truncation}
We have seen in Section \ref{section series representation} that for shot noise representations, the summand $H(\Gamma_k,U_k)$ in \eqref{generalised shot noise series} corresponds to jumps associated with the L\'evy measure. In particular, as $\|H(\cdot,\bm\xi)\|$ is nonincreasing, the summands are expressed in the descending order of jump magnitudes. Naturally, this implies that the first finitely many (large) jumps account for significantly more variation of the infinitely divisible random vector than the remaining smaller jumps \cite{imai2010quasi}. With the validated notion that the shot noise series can be reasonably approximated by its partial sums, our shot noise representation established previously provides us with a powerful method to approximating infinitely divisible laws for sampling. Namely, we do so by truncating the series representation to a finite sum. Simulation via shot noise is well scalable to higher dimensional settings (for example, see \cite{scherer2012shot}). Simulation via a finite truncation approach in the case of infinite jump intensity is typical for shot noise processes \cite{moller2006approximate,moller2005generalised}. In what follows, we describe a particular truncation scheme in more detail and provide examples and error analysis.

Let $X$ be an infinitely divisible random vector without Gaussian components and with an infinite L\'evy measure $\nu(d\mathbf z)$. To provide intuition for approximations via shot noise representation, we describe two finite truncation schemes with the inverse L\'evy measure method. Suppose a shot noise representation for $X$ is given by
\[
	X \stackrel{\mathscr L}{=} \sum_{k=1}^{+\infty}H_I(\Gamma_k,U_k)U_k,
\]
where the random sequences are as in Proposition \ref{prop series representations} (i). While we have assumed a scenario in which the compensation constants are not required for simplicity, it should be noted that their presence does not lead to any substantial difference in the analysis. Immediately, one may approximate $X$ by only including the summands corresponding to the index set $\{1,2,\cdots,n\}$ for a fixed truncation parameter $n\in\mathbb N$ rather than the entire infinite series. This is simple to implement, and it is clear that greater accuracy can be obtained by increasing $n$. However, there are two important aspects of this deterministic truncation scheme to consider. Firstly, by fixing $n$, we are conditioning on the number of jumps of the approximation, which may be undesirable for certain applied contexts in which the number of jumps is required to remain random. Secondly, while we can see that by truncating the series, we discard all jumps below some magnitude, for the direction $\bm\xi\in S^{d-1}$ this threshold magnitude is given by $H_I(\Gamma_n,\bm\xi)$, which is random.

An alternative truncation scheme is to instead perform summation with respect to the random index set $\{k\in\mathbb N: \Gamma_k \le n\}$, where $n>0$ is the truncation parameter. We refer to this framework as the \textit{Poisson truncation} approximation, which differs with the deterministic truncation scheme described previously by allowing the index set to depend on the underlying Poisson arrival times. In this way, we no longer condition on the number of jumps, but rather for each direction $\bm\xi\in S^{d-1}$, we include all jumps with magnitudes greater or equal to the deterministic threshold $H_I(n,\bm\xi)$. This fixed threshold immediately gives an indication of the error. Equivalently, in every direction $\bm\xi\in S^{d-1}$, this truncation method exactly simulates the tail of the L\'evy measure $\nu(d\mathbf z)$ over $(H_I(n,\bm\xi),+\infty)$.

In the case that the shot noise representation does not correspond to the inverse L\'evy measure method, while the shape of the domain simulated via Poisson truncation may not necessarily be as simple, it still holds that the simulated region is deterministic with finite measure and thus can still provide an indication of the error. We see that the Poisson truncation approximation of an infinitely divisible random vector is in essence an approximation by a compound Poisson random vector.
In the form of the generalized shot noise method of Theorem \ref{theorem generalised shot noise method}, the partial Levy measure described via the Poisson truncation $\{k\in \mathbb{N}:\,\Gamma_k\le n\}$, say $\nu_n(d\mathbf z)$, is given by 
\[
	\nu_n(B) = \int_{0+}^n \mathbb{P}(H(r,U)\in B)dr,\quad B\in\mathcal{B}(\mathbb{R}_0^d),
\]
for $n\in \mathbb{N}$. Hence, the total mass that the Poisson truncation describes is $\nu_n(\mathbb R_0^d) = n$.
Hereafter, we focus on the setting of the Poisson truncation method and reserve the notation $\nu_n$ for such a truncated Levy measure via the Poisson truncation approximation.

Note that the average number of summands under Poisson truncation is $n$, since the $\Gamma_k$'s in the index set $\{k\in\mathbb N:\Gamma_k\le n\}$ corresponds to the arrival times of the standard Poisson process on $[0,n]$. However, this may be merely an upper bound for the average number of jumps, as for example, some summands may evaluate to zero in the case of the rejection and thinning methods (see \eqref{rejection method series} and \eqref{thinning method series}, respectively).
Consider the case when the infinitely divisible random vector without Gaussian components has a finite L\'evy measure $\nu(d\mathbf z)$, that is, is a compound Poisson random vector. As the total jump intensity $\nu(\mathbb R_0^d)<+\infty$ is originally finite, the shot noise representation must almost surely have a finite number of nonzero terms. So in this case, there is no need to artificially truncate the series representation \eqref{generalised shot noise series}, that is,
\begin{equation}\label{random sum for finite Levy measure}
	\sum_{k=1}^{+\infty} \left(H\left(\Gamma_{k},U_k\right)-c_k\right) \stackrel{\mathscr L}{=} \sum_{k=1}^{N} \left(H\left(V_{(k)},U_k\right)-c_k\right),
\end{equation}
where $N$ is a Poisson random variable with rate $\nu(\mathbb R_0^d)$ and $\{V_{(k)}\}_{k\in\{1,\cdots,N\}}$ is a sequence of order statistics of $N$ iid uniform random variables on $(0,\nu(\mathbb R_0^d))$.

The representation \eqref{random sum for finite Levy measure} provides an alternative random sum representation of compound Poisson random vectors, compared to
$X \stackrel{\mathscr L}{=} \sum_{k=1}^N Y_k,$
where $\{Y_k\}_{k\in\mathbb N}$ is a sequence of iid random vectors distributed under $\nu(d\mathbf z)/\nu(\mathbb R^d_0)$. While the latter representation presents the jump structure as iid random vectors, the shot noise representation \eqref{random sum for finite Levy measure} decomposes the jumps in decreasing order of contribution to variation. For the sampling of the compound Poisson law, the random sum \eqref{random sum for finite Levy measure} may be more advantageous to implement if sampling the random sequence $\{Y_k\}_{k\in\mathbb N}$ is not as computationally convenient, which is often the case in multidimensional settings.

One should also note that by the memoryless property of the exponential distribution, the two sets of Poisson arrival times $\{\Gamma_k:\Gamma_k\le n,k\in\mathbb N\}$ and $\{\Gamma_k:\Gamma_k\in (n,m],k\in\mathbb N\}$ are independent for every $n<m$. Consequently, we can consider the L\'evy measures $(\nu - \nu_n)(d\mathbf z)$ and $(\nu_n-\nu_m)(d\mathbf z)$ as corresponding to independent components of the L\'evy process with L\'evy measure $\nu(d\mathbf z)$. This fact may be useful for incrementally simulating a L\'evy process based on the domain of its L\'evy measure, as well as for the analysis of error (Section \ref{section Levy-driven stochastic integrals}).

As an example of error analysis, we can consider the mean-squared error associated with the Poisson truncation of shot noise representation. Denote $\nu(d\mathbf z)$ and $\mu(d\mathbf z,ds)$ as the L\'evy measure and Poisson random measure associated with the infinitely divisible random vector of interest, and denote $\nu_n(d\mathbf z)$ and $\mu_n(d\mathbf z,ds)$ as that of its Poisson truncation approximation. Assume the L\'evy measure $\nu(d\mathbf z)$ is isotropic. The truncation error for the shot noise representation \eqref{generalised shot noise series} is expressed as the tail series
\[
	\sum_{k=1}^{+\infty}H(\Gamma_k,U_k) - \sum_{\{k\in\mathbb N:\Gamma_k\le n\}}H(\Gamma_k,U_k) = \sum_{\{k\in\mathbb N: \Gamma_k>n\}}H(\Gamma_k,U_k).
\]
By Theorem \ref{theorem LI decomposition} and the It\^o-Wiener isometry, the mean-squared error is given by
\begin{equation}\label{mse}
	\mathbb E\left[\left\|\sum_{\{k\in\mathbb N: \Gamma_k>n\}}H(\Gamma_k,U_k)\right\|^2\right] = \mathbb E\left[\left\|\int_0^1\int_{\mathbb R_0^d}\mathbf z\,((\mu-\mu_n) - (\nu - \nu_n)\mathbbm1_{(0,1]}(\|\mathbf z\|))(d\mathbf z,ds)\right\|^2\right] = \int_{\mathbb R_0^d}\|\mathbf z\|^2\,(\nu-\nu_n)(d\mathbf z),
\end{equation}
for every $n>0$. For the case where the L\'evy measure $\nu(d\mathbf z)$ is not isotropic, application of the It\^o-Wiener isometry in \eqref{mse} can still be invoked in the case of the inverse L\'evy measure method for $n$ sufficiently large such that the support of $(\nu-\nu_n)(d\mathbf z)$ is contained in the unit ball.
As an example, we present the evaluation of the mean-squared truncation error \eqref{mse} for the isotropic stable law in the following.

\begin{example}\label{example stable error analysis}
	Consider the Poisson truncation of the shot noise representation of the stable vector \eqref{series rep of stable processes} obtained via the inverse L\'evy measure method \cite[Section 4]{bondesson1982simulation}. Suppose that the measure $\sigma(d\bm\xi)$ is isotropic, that is, $c_k\equiv0$.
	As $(\nu-\nu_n)(d\mathbf z) = \mathbbm1(\|\mathbf z\|\le(\alpha n/a)^{-1/\alpha})\nu(d\mathbf z)$, the mean-squared error \eqref{mse} of the Poison truncation approximation is given by
	\begin{equation}\label{mse of truncation of stable}
		\mathbb E\left[\left\|\sum_{\{k\in\mathbb N:\Gamma_k>n\}}\left(\frac{\alpha\Gamma_k}{a}\right)^{-1/\alpha}U_k\right\|^2\right] = \int_{\|\mathbf z\|\le (\alpha n/a)^{-1/\alpha}}\|\mathbf z\|^2\,\nu(d\mathbf z) = a^{2/\alpha}\frac{\alpha^{1-2/\alpha}}{2-\alpha}n^{1-2/\alpha},
	\end{equation}
	for every $n>0$. As mentioned previously, the mean-squared error \eqref{mse of truncation of stable} still holds if $\sigma(d\bm\xi)$ is not isotropic for sufficiently large $n>0$ satisfying $(\alpha n/a)^{-1/\alpha} < 1$. We see that in the case of the stable random vector, mean-squared convergence of the Poisson truncation approximation is very fast for $\alpha$ close to zero, so truncating the series \eqref{series rep of stable processes} to a relatively small number of terms provides a good approximation. However, for $\alpha$ close to two, convergence is much slower and thus, to achieve a level of accuracy, significantly more summands must be computed. This is intuitive; increasing $\alpha$ leads to thinner tails and hence the magnitudes of the largest jumps decrease. Consequently, the variation explained by the largest jumps becomes diluted. In the univariate case, it is suggested \cite{bondesson1982simulation} that the approximation can be improved by the inclusion of a normal random variable with variance given by \eqref{mse of truncation of stable}. This idea is the basis for Gaussian approximation of the truncation error (Section \ref{section Gaussian approximation}), which can be extended to the multivariate setting. Alternatively, investigation of the truncation error for the representation for stable laws can be found in \cite{bentkus1996bounds,bentkus2001levy} in terms of optimal bounds for the variation. 
	$\hfill \qed$
\end{example}


The idea that the first few summands of the shot noise representation \eqref{generalised shot noise series} accounts for a large amount of variation has useful applications outside of sampling (Section \ref{section computing expectations}). For instance, it has been shown in \cite{breton2010regularity} that absolute continuity of the law of the first few summands guarantees the absolute continuity of the entire shot noise series.

\subsection{Comparisons among shot noise representations}\label{section analysis of finite truncation}
We compare among the errors from Poisson truncation approximations of the various shot noise representation methods established in Proposition \ref{prop series representations}. Consider the following general result \cite{imai2013numerical}.

\begin{theorem}\label{theorem finite truncation comparison}
	Assume the setting of Proposition \ref{prop series representations}. Fix $n>0$. Let $\nu_{k,n}(d\mathbf z)$, $k=1,2,3,4$, be the L\'evy measures of the partial sums obtained from inverse L\'evy measure \eqref{inverse Levy measure method series}, rejection \eqref{rejection method series}, thinning \eqref{thinning method series} and Bondesson's \eqref{Bondesson's method series} methods over the index set $\{k\in\mathbb N:\Gamma_k\le n\}$, respectively. Then, the following hold:
	\begin{enumerate}[(i)]
		\setlength{\parskip}{0cm}
		\setlength{\itemsep}{0cm}
		\item If $\nu(\mathbb R_0^d) \ge n$, then $\nu_{k,n}(\mathbb R_0^d) = n$ for $k=1,2,3,4$.
		\item For every $x>0$ and $C\in\mathcal B(S^{d-1})$,
		\[
			\nu_{k,n}((x,+\infty)C)\le \nu_{1,n}((x,+\infty)C)\le\nu((x,+\infty)C),\quad k=2,3,4.
		\]
		\item For every $q\ge0$ such that $\int_{\|\mathbf z\|>1}\|\mathbf z\|^q\,\nu(d\mathbf z)<+\infty$,
		\[
			\int_{\mathbb R_0^d}\|\mathbf z\|^q\,\nu_{k,n}(d\mathbf z) \le\int_{\mathbb R_0^d}\|\mathbf z\|^q\,\nu_{1,n}(d\mathbf z)<+\infty,\quad k=2,3,4.
		\]
		\item For every $q\ge2$,
		\[
			\int_{\mathbb R_0^d}\|\mathbf z\|^q(\nu-\nu_{1,n})\,(d\mathbf z) \le\int_{\mathbb R_0^d}\|\mathbf z\|^q\,(\nu - \nu_{k,n})(d\mathbf z),\quad k=2,3,4.
		\]
	\end{enumerate}
\end{theorem}

Note that the above result holds for the general L\'evy measure $\nu(d\mathbf z)$.
The result (i) echoes our understanding of Poisson truncation in simulating a mass $n$ irrespective of the underlying shot noise representation, and that truncation of distinct representations corresponds to the simulation of different regions of the L\'evy measure $\nu(d\mathbf z)$.
By (ii), it is clear that under the Poisson truncation method, no method presented in Proposition \ref{prop series representations} can simulate the tail of the L\'evy measure more accurately than the inverse L\'evy measure method. Moreover, in light of (i), we deduce that the other methods simulate regions of the L\'evy measure closer to the origin. By (iii), we see that the inverse L\'evy measure method captures more variation of the jumps than the other methods. The result (iv) is also useful, as for $q=2$, the integrals correspond to the variances of discarded jumps, which can be viewed as representative of the error associated with the approximation. The smaller the variance, the better the approximation. Thus, we see that at least under the Poisson truncation scheme, the L\'evy measure method is preferred over the others. The main drawback to the inverse L\'evy measure method is that the tail of the L\'evy measure may not be conveniently invertible (for instance, see Example \ref{example gamma law}), however in such cases one may still numerically invert the L\'evy measure as required \cite{imai2013numerical}.
If the numerical inversion of the tail of the L\'evy measure is computationally expensive, then one would prefer not to use the inverse L\'evy measure method. Additionally, the error of numerical inversion may accumulate in the series.

While Theorem \ref{theorem finite truncation comparison} only addresses inverse L\'evy measure, rejection, thinning and Bondesson's methods for any L\'evy measure, it appears that the result extends beyond those methods. As an example, the following demonstrates that the result holds in the case of Rosinski's series representation for the tempered stable law \cite{imai2011finite}.

\begin{example}[Comparisons among shot noise representations of the tempered stable law]\label{example Rosinski's series comparison}
	Assume the setting of Example \ref{example tempered stable law}. Let $\nu_{k,n}(d\mathbf z)$, $k=1,2,3,4,5$, denote the L\'evy measure corresponding to the Poisson truncation of the shot noise representation by the inverse L\'evy measure method \eqref{inverse Levy measure method for TS law}, the thinning method \eqref{thinning method for TS law 1} and \eqref{thinning method for TS law 2}, the rejection method \eqref{rejection method for TS law} and Rosi\'nski's representation \eqref{eq tempered stable law series representation}, respectively.
	
	\begin{enumerate}[(i)]
		\setlength{\parskip}{0cm}
		\setlength{\itemsep}{0cm}		
		\item It holds that for every $(x,C)\in(0,+\infty)\times\mathcal B(S^{d-1})$,
		\[
			\nu_{k,n}((x,+\infty)C) \le \nu_{1,n}((x,+\infty)C), \quad k=2,3,4,5,
		\]
		and
		\begin{align*}
			\nu_{4,n}((x,+\infty)C) &= \nu\left(\left(\left(\frac{\alpha n}{m_{\alpha,\rho}}\right)^{-1/\alpha}\vee x,+\infty\right)C\right),\\
			\nu_{5,n}((x,+\infty)C) &= \left(\frac{\alpha n}{m_{\alpha,\rho}}x^\alpha \wedge 1\right)\nu((x,+\infty)C).
		\end{align*}
		
		\item It holds that for every $q\ge0$ such that $\int_{\|\mathbf z\|>1}\|\mathbf z\|^q\,\nu(d\mathbf z)<+\infty$,
		\[
			\int_{\mathbb R_0^d}\|\mathbf z\|^q\,\nu_{k,n}(d\mathbf z) \le \int_{\mathbb R_0^d}\|\mathbf z\|^q\,\nu_{1,n}(d\mathbf z) \le \int_{\mathbb R_0^d}\|\mathbf z\|^q\,\nu(d\mathbf z), \quad k=2,3,4,5.
		\]
	\end{enumerate}
	We observe that the inverse L\'evy measure method outperforms the methods of Proposition \ref{prop series representations}, as expected due to Theorem \ref{theorem finite truncation comparison}, but as well as Rosi\'nski's shot noise representation.
	We also state the L\'evy measure $\nu_{5,n}(d\mathbf z)$ as follows. Recall from Section \ref{section ID examples} that the L\'evy measure for the tempered stable law can be written in polar form as $\nu(dr,d\bm\xi) = r^{-\alpha-1}q(r,\bm\xi)dr\,\sigma(d\bm\xi)$, which is often provided as such. Then, it holds that \cite{cohen2007gaussian}
	\[
		\nu_{5,n}(dr,d\bm\xi) = h_n(r,\bm\xi)dr\,\sigma(d\bm\xi),
	\]
	where
	\[
		h_n(r,\bm\xi) =
		\begin{dcases}
			n\left(r^{-1} q(r,\bm\xi) - r^{\alpha-1}\int_r^{+\infty}s^{-\alpha-1}q(r,\bm\xi)\,ds\right),\quad& r\in(0,n^{-1/\alpha}),\\
			r^{-\alpha-1}q(r,\bm\xi),\quad& r\ge n^{-1/\alpha}.
		\end{dcases}
	\]
	From this, one can proceed with error analysis by investigating the L\'evy measure corresponding to the truncation error.
	$\hfill\qed$
\end{example}

\begin{example}[Truncation of shot noise representations of the gamma law]\label{example gamma truncation error}
	Assuming the setting of Example \ref{example gamma law}. Let $\nu_{k,n}(d\mathbf z)$, $k=1,2,3,4$, denote the L\'evy measure corresponding to the Poisson truncation of the shot noise representation by the inverse L\'evy measure \eqref{inverse Levy measure method for gamma law}, the rejection \eqref{rejection method for gamma law}, the thinning \eqref{thinning method for gamma law} and Bondesson's \eqref{Bondesson's method for gamma law} methods, respectively. We present computations relating to the error of the approximations \cite[Section 5.4]{kawai2017sample} as follows.
	\begin{enumerate}[(i)]
		\setlength{\parskip}{0cm}
		\setlength{\itemsep}{0cm}
		\item For every $x>0$, it holds that
		\[
			\int_x^{+\infty}\nu_{k,n}(dz) = \begin{cases}
				(n\wedge a E_1(\beta x)),\quad &\text{if }k=1,\\
				a E_1\left((e^{n/a}-1)^{-1}\vee\beta x\right),\quad &\text{if }k=2,\\
				\int_{\beta x}^{+\infty} (n \wedge a/u)e^{-u}\,du,\quad &\text{if }k=3,\\
				a\left(E_1(\beta x) - E_1(\beta xe^{n/a})\right),\quad &\text{if }k=4.
			\end{cases}
		\]
		
		\item It holds that
		\[
			\int_{0+}^{+\infty}z\,\nu_{k,n}(dz) = \frac{a}{\beta}\times\begin{cases}
				e^{-E_1^{-1}(n/a)}, \quad& \text{if } k=1,\\
				e^{-(e^{n/a}-1)^{-1}}, \quad& \text{if } k=2,\\
				\frac{n}{a}\gamma(2,a/n) + e^{-a/n}, \quad& \text{if } k=3,\\
				1-e^{-n/a}, \quad& \text{if } k=4.
			\end{cases}
		\]
		
		\item The variance of the discarded jumps is given by
		\[
			\sigma_{k,n}^2 \coloneqq \sigma^2 - \int_{0+}^{+\infty}z^2\,\nu_{k,n}(dz) = \sigma^2\times\begin{cases}
				\gamma\left(2,E_1^{-1}(n/a)\right), \quad&\text{if }k=1,\\
				\gamma\left(2,(e^{n/a}-1)^{-1}\right),\quad&\text{if }k=2,\\
				\gamma\left(2,a/n\right) - \frac{n}{a}\gamma(3,a/n),\quad&\text{if }k=3,\\
				e^{-2n/a}, \quad&\text{if }k=4,
			\end{cases}
		\]
		where $\sigma^2\coloneqq \int_{0+}^{+\infty}z^2\,v(dz) = a/\beta^2$ and we denote $\gamma(a,x)\coloneqq \int_{0+}^x u^{a-1}e^{-u}\,du$ as the lower incomplete gamma function for $a>0$ and $x>0$.
	\end{enumerate}
	
	The result (iii) corresponds to Theorem \ref{theorem finite truncation comparison} (iv) on the variance of the discarded jumps. Without direct comparison, we already know from Theorem \ref{theorem finite truncation comparison} that the variance associated with the inverse L\'evy measure method cannot be outperformed by the other methods examined. $\hfill\qed$
\end{example}

To mention a more recent example of Poisson truncation error analysis, it is found in \cite{massing2018simulation} that in the case of the $t$-distribution, the mean-squared truncation error for the inverse L\'evy measure and rejection methods are both bounded by $2\nu/(\pi(n-1))$, where $\nu$ is the degrees of freedom parameter.

\subsection{Gaussian approximation of small jumps}\label{section Gaussian approximation}
In the Poisson truncation approximation of an infinitely divisible random vector, as the norm of the kernel $r\mapsto \|H(r,\mathbf v)\|$ in \eqref{generalised shot noise series} is nonincreasing, we discard jumps of sizes within some neighbourhood about the origin. Under an appropriate method (for example, the inverse L\'evy measure method), suppose the truncation parameter is sufficiently large so only jumps of magnitudes less than a small $\epsilon>0$ are discarded.
Then, the variance of the discarded jumps is given by
\begin{equation}\label{variance of discarded jumps}
	\sigma^2(\epsilon) \coloneqq \int_{\|\mathbf z\|<\epsilon}\|\mathbf z\|^2\,\nu(d\mathbf z)<+\infty.
\end{equation}
A natural idea to pursue is to approximate the discarded jumps by a Gaussian random vector \cite{bondesson1982simulation}. This is particularly practical when the shot noise representation of the infinitely divisible random vector converges slowly, that is, when $r\mapsto \|H(r,\mathbf v)\|$ decreases slowly, as is the case with stable random vectors with stability $\alpha$ close to two (Example \ref{example stable error analysis}).

In what follows, we restrict ourselves to the one-dimensional setting. Specifically, let $X$ be an infinitely divisible random variable with L\'evy measure $\nu(dz)$ and $X^\epsilon$ the truncation approximation of $X$ with jump sizes of $|z|<\epsilon$ discarded. Define the error, consisting of the discarded small jumps, as 
\begin{equation}\label{error process}
	X_\epsilon \coloneqq X - X^\epsilon.
\end{equation}
We seek the conditions in which it is valid to approximate $X_\epsilon/\sigma(\epsilon)$ with a standard normal random variable, as such a Gaussian approximation cannot be applied for every infinitely divisible law, with the Poisson distribution being a simple counterexample. We examine some notions from \cite{asmussen2001approximations} on this matter. A necessary and sufficient condition for the Gaussian approximation of discarded jumps to hold is as follows.

\begin{theorem}[Gaussian approximation of discarded jumps] \label{theorem Gaussian approximation of discarded jumps}
	Let $\sigma(\epsilon)$ and $X_\epsilon$ be defined as in \eqref{variance of discarded jumps} and \eqref{error process}, respectively. Then, we have that $X_\epsilon/\sigma(\epsilon)$ converges in law to the standard normal distribution as $\epsilon\to0$ if and only if for every $c>0$,
	\begin{equation}\label{Gaussian approximation of discarded jumps}
		\sigma(c\sigma(\epsilon)\wedge\epsilon)\sim\sigma(\epsilon),
	\end{equation}
	as $\epsilon\to0$.
\end{theorem}

A sufficient yet more verifiable condition is given in \cite[Proposition 2.1]{asmussen2001approximations}, which shows that a Gaussian approximation of the discarded jumps is valid when its variation decreases at a rate slower than its upper bound. If $\sigma(\epsilon)/\epsilon\to+\infty$ as $\epsilon\to0$,
then the condition \eqref{Gaussian approximation of discarded jumps} holds. Moreover, if the L\'evy measure $\nu(dz)$ does not have any atoms in a neighbourhood about the origin, then this condition is equivalent to the condition \eqref{Gaussian approximation of discarded jumps}.

There is a profound relation between the decay of the generalised inverse of the tail of the L\'evy measure to the validity of approximating the discarded jumps by a normal random variable \cite[Proposition 2.2]{asmussen2001approximations}. Suppose the L\'evy measure $\nu(dz)$ is symmetric and infinite. By symmetry, the kernel \eqref{inverse Levy measure method kernel} for the inverse L\'evy measure method is independent of the direction and simplifies to $H_I(r) = \inf\{u\in(0,+\infty): 2\nu((u,\infty)) < r\}$ for $r>0$. If, for every $a>0$, it holds that $\lim_{t\to\infty}H_I(t+a)/H_I(t)=1$, then $X_\epsilon/\sigma(\epsilon)$ converges in law to the standard normal distribution as $\epsilon\to0$. For example, it is straightforward to show that the stable distribution satisfies this condition. Thus, the error from truncation asymptotically behaves like a normal random variable. However, not every infinitely divisible distribution enjoys the benefit of Gaussian approximation, as the gamma law is a counterexample where the small jumps are infinitely active yet not intense enough to resemble Brownian infinite variation \cite{asmussen2001approximations,bondesson1982simulation}. For Gaussian approximation of discarded jumps in the case of higher dimensions with more general methods than the inverse L\'evy measure method, we refer the reader to \cite{cohen2007gaussian}.

In the light of Theorem \ref{theorem finite truncation comparison} (iv), the error analysis of Poisson truncation can be revisited from the perspective of Gaussian approximation of the error \cite{imai2013numerical}. Suppose that the Poisson truncation error of an infinitely divisible random variable can be approximated by a normal random variable.
Denote the L\'evy measure associated with the Poisson truncation of each of the shot noise representation methods in Proposition \ref{prop series representations} as $\nu_{k,n}(dz)$, $k=1,2,3,4$. Then, for each of the methods, the normal random variable which approximates the discarded jumps has the variance
\[
	\sigma_{k,n}^2 \coloneqq\int_{\mathbb R_0}z^2\,(\nu - \nu_{k,n})(dz), \quad k=1,2,3,4.
\]
Recall Theorem \ref{theorem finite truncation comparison} (iv), in which we have seen that $\sigma_{1,n}\le\sigma_{k,n}$ holds for $k=2,3,4$, where $k=1$ corresponds to the inverse L\'evy measure method. Thus, no normal random variable used in approximating the error of the methods examined has a lower variance than that associated with the inverse L\'evy measure method. This can be interpreted as the inverse L\'evy measure method being an optimal approximation in the $L^2$ sense, in comparison to the other methods under Poisson truncation. It is worth noting that this discussion is also important in higher dimensions.

\section{Numerical schemes via truncation of shot noise representation}\label{section truncation for processes}
So far, we have focused primarily on infinitely divisible laws without Gaussian components, with glimpses of shot noise representation for L\'evy processes in Sections \ref{section Levy-Ito decomposition} and \ref{section inverse Levy measure method}. We have seen the correspondence between infinitely divisible laws and L\'evy processes reflected in their shot noise representations. In Section \ref{section truncation}, we have discussed the method of approximating infinitely divisible random vectors via the truncation of shot noise representation. In this section, we focus entirely on the technique of Poisson truncation of shot noise representations for simulating L\'evy processes (Section \ref{section simulating Levy processes}), infinitely divisible processes (Section \ref{section Levy-driven stochastic integrals}) and fields (Section \ref{section Levy-driven random fields}) and L\'evy-driven stochastic differential equations (Section \ref{section Levy-driven SDEs}). For simulation purposes, while the truncation of shot noise representation is useful for infinitely divisible laws and L\'evy processes, often times it is merely optional among a wider array of numerical methods available. However, for infinitely divisible processes and L\'evy-driven stochastic differential equations, truncation of shot noise representations offers a very effective approach due to the necessity of jump-based approximation in distilling the structure of the stochastic process in general. As a means of demonstrating the effectiveness of truncation of shot noise representation, we provide various numerical illustrations along the way.

\subsection{Simulating L\'evy processes}\label{section simulating Levy processes}
The inverse L\'evy measure method established in Section \ref{section inverse Levy measure method} presents us with a shot noise representation \eqref{inverse Levy measure method series over [0,1]} for a L\'evy process $\{X_t:t\in[0,1]\}$ in $\mathbb R^d$ without Gaussian components.
Recall that this is indeed a shot noise representation for the infinitely divisible random vector $X_1$, except the summands are scattered uniformly and the centres are subtracted in proportion over the unit time interval.

We begin by briefly describing the general shot noise method for additive processes, which generalise L\'evy processes by relaxing the property of stationary increments. Consider an additive process described by a more general L\'evy-It\^o decomposition \cite[Chapter 4]{sato1999levy}
\[
	X_t = \int_0^t\int_{\|\mathbf z\|\in(0,1]}{\mathbf z}(\mu(d\mathbf z,ds) - \nu(d\mathbf z)\lambda(s)ds) - \int_0^t\int_{\|\mathbf z\|>1}\mathbf z\,\mu(d\mathbf z,ds),\quad t\ge0,
\]
where $\lambda(\cdot)\ge0$ is a nonnegative function with support including $[0,T]$ for a fixed $T>0$ and $\mu(d\mathbf z,ds)$ is a Poisson random measure on $\mathbb R_0^d\times[0,+\infty)$ with intensity measure $\nu(d\mathbf z)\lambda(s)ds$. Suppose the L\'evy measure $\nu(d\mathbf z)$ satisfies the decomposition \eqref{generalised shot noise method decomposition}. By a change of variable, it holds that for every $(B,S)\in\mathcal B(\mathbb R_0^d)\times\mathcal B([0,T])$,
\[
	\int_S\int_B\nu(d\mathbf z)\lambda(s)ds = \int_{0+}^{+\infty}\mathbb P(H(r,U)\in B)\,dr\int_S\lambda(s)\,ds = \int_{0+}^{+\infty}\mathbb P\left(H\left(\frac{r}{\int_0^T\lambda(s)\,ds},U\right)\in B\right)\,dr \,\frac{\int_S\lambda(s)\,ds}{\int_0^T\lambda(s)\,ds}.
\]
Thus, by scattering the summands in \eqref{generalised shot noise series} according to the density $\lambda(\cdot)/\int_0^T\lambda(s)\,ds$ instead of uniformly, and replacing $\Gamma_k$ by $\Gamma_k/\int_0^T\lambda(s)\,ds$, we obtain the shot noise representation for the additive process as
\begin{equation}\label{series representation for additive processes}
	\{X_t:t\in[0,T]\} \stackrel{\mathscr L}{=} \left\{\sum_{k=1}^{+\infty}\left(H\left(\frac{\Gamma_k}{\int_0^T\lambda(s)\,ds},U_k\right)\mathbbm1_{[0,t]}(S_k) - c_k\int_0^t\lambda(s)\,ds\right):t\in[0,T]\right\},
\end{equation}
where $\{S_k\}_{k\in\mathbb N}$ is a sequence of iid random variables with density $\lambda(\cdot)/\int_0^T\lambda(s)\,ds$ and $\{c_k\}_{k\in\mathbb N}$ is a suitable sequence of centres. In particular, by considering the case where $\lambda\equiv 1$, we obtain the series representation for any L\'evy process over $[0,T]$ based on the generalised shot noise method (Theorem \ref{theorem generalised shot noise method}) as
\begin{equation}\label{generalised shot noise series for Levy process}
	\{X_t:t\in[0,T]\} \stackrel{\mathscr L}{=} \left\{\sum_{k=1}^{+\infty} \left(H\left(\frac{\Gamma_{k}}{T},U_k\right)\mathbbm1_{[0,t]}(T_k)-tc_k\right): t\in[0,T]\right\},
\end{equation}
where $\{T_k\}_{k\in\mathbb N}$ is a sequence of iid uniform random variables on $(0,T)$ independent of the other random sequences. It is known \cite[Theorem 5.1]{rosinski2001series} that the infinite series \eqref{generalised shot noise series for Levy process} converges almost surely and uniformly on $[0,T]$. Moreover, there exists a version of the random sequences in the right hand side of \eqref{generalised shot noise series for Levy process} such that the equality holds almost surely.

We also mention L\'evy processes of type G, which admit special shot noise representations related to the inverse L\'evy measure method of Section \ref{section inverse Levy measure method}. A L\'evy process $\{X_t:t\ge0\}$ in $\mathbb R$ is said to be of type G if its increments can be represented in law as $X_{t_2} - X_{t_1} \stackrel{\mathscr L}{=} V^{1/2}G$ for $t_1<t_2$, where $V$ is a nonnegative infinitely divisible random variable with L\'evy measure $\nu(dz)$ and $G$ is a standard normal random variable. With $H_I(\cdot)$ defined in \eqref{inverse Levy measure kernel} but independent of the directional argument, it holds that \cite{rosinski1991class}
\begin{equation}\label{type G series}
	\{X_t:t\in[0,T]\} \stackrel{\mathscr L}{=} \left\{\sum_{k=1}^{+\infty} G_k H_I\left(\frac{\Gamma_k}{T}\right)^{1/2}\mathbbm1_{[0,t]}(T_k):t\in[0,T]\right\},
\end{equation}
where $\{\Gamma_k\}_{k\in\mathbb N}$ is a sequence of standard Poisson arrival times, $\{G_k\}_{k\in\mathbb N}$ is a sequence of iid standard normal random variables and $\{T_k\}_{k\in\mathbb N}$ is a sequence of iid uniform random variables over $(0,T)$ such that the random sequences are mutually independent. We refer to \cite{wiktorsson2002simulation} for results on error analysis for Poisson truncation of the shot noise representation \eqref{type G series}. We also mention that the L\'evy process of type G admits the subordinated Brownian motion representation $\{X_t:t\ge0\}\stackrel{\mathscr L}{=} \{B_{V_t}:t\ge0\}$, where $\{B_t:t\ge0\}$ is a standard Brownian motion in $\mathbb R$ and $\{V_t:t\ge0\}$ is a subordinator with L\'evy measure $\nu(dz)$. Thus, an alternative simulation method to the Poisson truncation of the shot noise representation \eqref{type G series} is via the truncation of a shot noise representation for the subordinator $\{V_t:t\ge0\}$.

With the knowledge of shot noise representations of L\'evy and additive processes and Poisson truncation from Section \ref{section truncation}, we have a powerful and general method for simulating their sample paths. All general as well as specific Poisson truncation error analyses for infinitely divisible laws from Section \ref{section truncation} carry over to the case of L\'evy and additive processes. In particular, Gaussian approximation of discarded jumps holds for the case of L\'evy and additive processes by including a Brownian motion (and its deterministically time-changed version), instead of a normal random vector, scaled by the variance \eqref{variance of discarded jumps} of the small jumps \cite{asmussen2001approximations}. The conditions for the one-dimensional setting as well as their multidimensional generalisations (Section \ref{section Gaussian approximation}) apply in this setting.

\subsubsection{Simulation recipes}
The implementation of Poisson truncation of shot noise representation is attractively clear-cut. Hereafter, we focus on L\'evy processes, as the generalisation to additive processes is rather trivial in light of \eqref{series representation for additive processes} and \eqref{generalised shot noise series for Levy process}. Due to the uniform scattering of jumps along the time interval $[0,T]$, it is little more than simply a matter of independently sampling the random sequences appearing in the series \eqref{generalised shot noise series for Levy process} and evaluating the shot noise kernel $H$. To generate the sequence of Poisson arrival times $\{\Gamma_k\}_{k\in\mathbb N}$, one may take advantage of the fact that the interarrival times are iid standard exponential random variables $\{E_k\}_{k\in\mathbb N}$. Based on the Poisson truncation of the shot noise representation \eqref{generalised shot noise series for Levy process}, we describe the simulation recipe for generating approximate sample paths of the L\'evy process in $\mathbb R^d$ over the time interval $[0,T]$ as follows.

\begin{enumerate}[\bf Step 1.]
	\setlength{\parskip}{0cm}
	\setlength{\itemsep}{0cm}
	\item Generate a standard exponential random variable $E_1$. If $E_1\le nT$, then assign $\Gamma_1 \leftarrow E_1$. Otherwise, return the degenerate zero process as the approximate sample path and terminate the algorithm.
	
	\item While $\Gamma_k \le nT$, generate a standard exponential random variable $E_{k+1}$ and assign $\Gamma_{k+1} \leftarrow \Gamma_k + E_{k+1}$. Denote this sequence as $\{\Gamma_k\}_{k\in\{1,\cdots,N\}}$, where $N$ satisfies $\Gamma_N\le nT < \Gamma_{N+1}$.
	
	\item Generate iid uniform random variables $\{T_k\}_{k\in\{1,\cdots,N\}}$ on $(0,T)$.
	
	\item Generate iid random vectors $\{U_k\}_{k\in\{1,\cdots,N\}}$ with the distribution specified by the shot noise representation \eqref{generalised shot noise series for Levy process}.
	
	\item For every $k\in\{1,\cdots,N\}$, compute the jump $J_k \leftarrow H(\Gamma_k/T,U_k)$.
	
	\item Sort $\{(J_k,T_k)\}_{k\in\{1,\cdots,N\}}$ to obtain $\{(J_{(k)},T_{(k)})\}_{k\in\{1,\cdots,N\}}$ such that $T_{(1)}<\cdots< T_{(N)}$.
	
	\item Assign $X_0 \leftarrow 0$ and $T_{(0)} \leftarrow 0$. For every $k\in\{1,\cdots,N\}$, assign $X_k \leftarrow X_{k-1} + J_{(k)} - (T_{(k)} - T_{(k-1)})\sum_{j=1}^N c_j$.
	
	\item Return $\{X_0,X_1,\cdots,X_N\}$ as the position of the approximate sample path at the sample times $\{T_{(0)},T_{(1)},\cdots,T_{(N)}\}$.
\end{enumerate}

In the case where the centres $\{c_k\}_{k\in\mathbb N}$ are not necessary for the almost sure convergence of the shot noise series \eqref{generalised shot noise series for Levy process}, the assignments in step 7 simplify to $X_k \leftarrow X_{k-1} + J_{(k)}$. We can alternatively sample the Poisson arrival times in the spirit of the compound Poisson representation \eqref{random sum for finite Levy measure}, as follows.

\begin{enumerate}[\bf Step 1.]
	\setlength{\parskip}{0cm}
	\setlength{\itemsep}{0cm}
	\item Generate a Poisson random variable $N$ with rate $nT$. If $N=0$, then return the degenerate zero process as the approximate sample path and terminate the algorithm.
	
	\item Generate $N$ iid uniform random variables $\{V_{k}\}_{k\in\{1,\dots,N\}}$ on $(0,nT)$ and sort in ascending order to obtain $\{V_{(k)}\}_{k\in\{1,\dots,N\}}$.
	
	\item Generate iid uniform random variables $\{T_k\}_{k\in\{1,\cdots,N\}}$ on $(0,T)$.
	
	\item Generate iid random vectors $\{U_k\}_{k\in\{1,\cdots,N\}}$ with the distribution specified by the shot noise representation \eqref{generalised shot noise series for Levy process}.
	
	\item For every $k\in\{1,\cdots,N\}$, compute the jump $J_k \leftarrow H(V_{(k)}/T,U_k)$.
	
	\item Sort $\{(J_k,T_k)\}_{k\in\{1,\cdots,N\}}$ to obtain $\{(J_{(k)},T_{(k)})\}_{k\in\{1,\cdots,N\}}$ such that $T_{(1)}<\cdots< T_{(N)}$.
	
	\item Assign $X_0 \leftarrow 0$ and $T_{(0)} \leftarrow 0$. For every $k\in\{1,\cdots,N\}$, assign $X_k \leftarrow X_{k-1} + J_{(k)} - (T_{(k)} - T_{(k-1)})\sum_{j=1}^N c_j$.
	
	\item Return $\{X_0,X_1,\cdots,X_N\}$ as the position of the approximate sample path at the sample times $\{T_{(0)},T_{(1)},\cdots,T_{(N)}\}$.
\end{enumerate}

Note that in the above simulation recipes, we work under the truncation $\{k\in\mathbb N : \Gamma_k\le nT\}$ instead of $\{k\in\mathbb N : \Gamma_k\le n\}$ to simulate a mass $n$ of the L\'evy measure per unit time (recall Theorem \ref{theorem finite truncation comparison} (i)).
While the two simulation recipes offered only differ in the method for sampling Poisson arrival times, we highlight some contrasting features which may be important in practice.
The second simulation recipe, based on the conditional uniformity of Poisson arrival times, is simpler to implement. In particular, the lack of a need for a while loop makes the this recipe more natural to implement from a functional programming perspective. In contrast, the first recipe based on successive exponential sampling carries the advantage of easy modification for adaptive truncation. That is, one may continually sample summands until some criterion, which may dynamically update, is reached. The second recipe based on conditional uniformity of Poisson arrivals does not share this advantage, as resampling the Poisson random variable for the number of jumps with a larger rate does not guarantee a greater number of jumps. Additionally, previous jump timings cannot be reused and would require complete resampling from scratch, which would become a source of inefficiency. Hence, in some scenarios, there are grounds to prefer one recipe over the other.

\subsubsection{Numerical illustrations}
\begin{figure}[h]\centering
	\begin{tabular}{cc}
		\includegraphics[width=90mm]{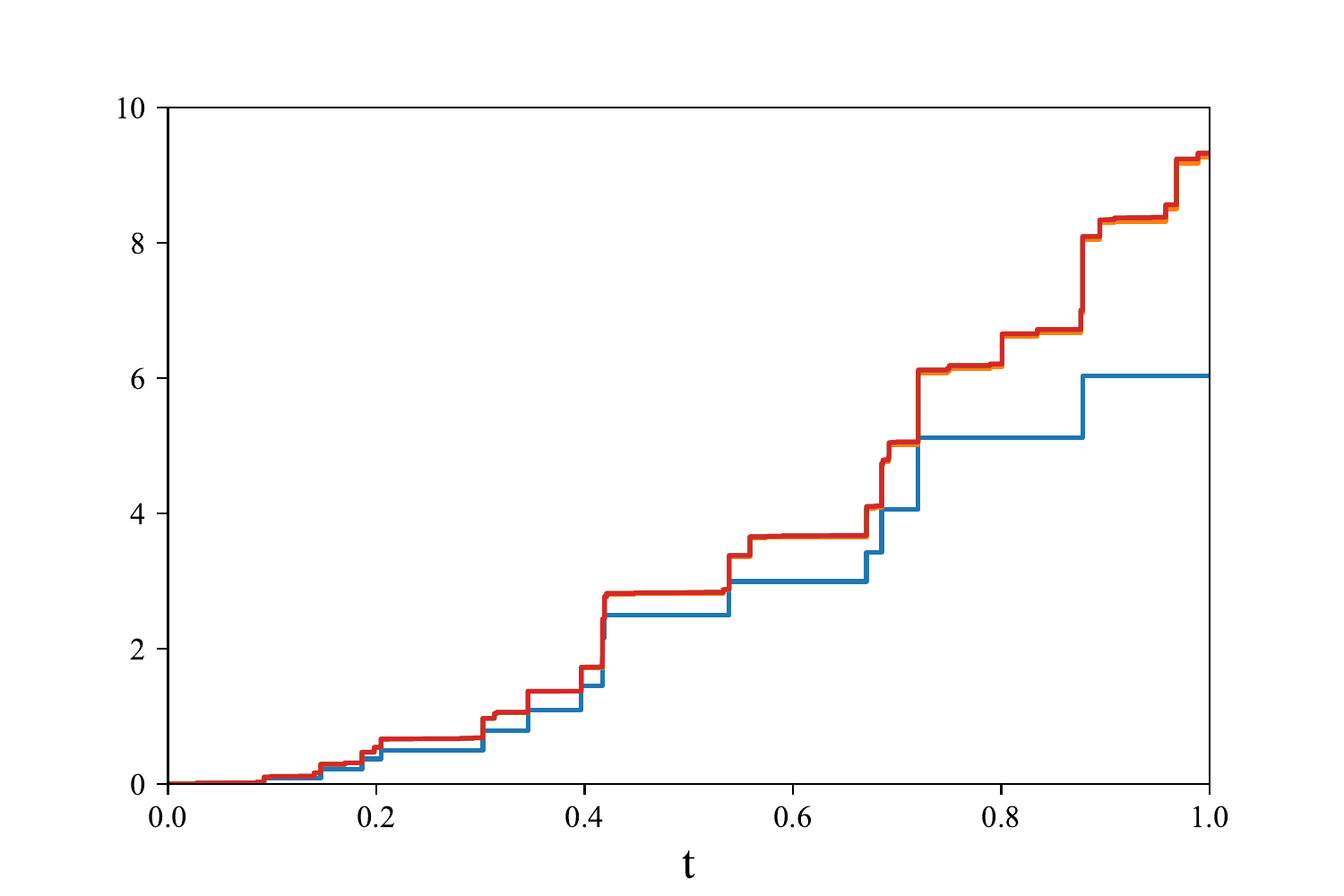} & \includegraphics[width=90mm]{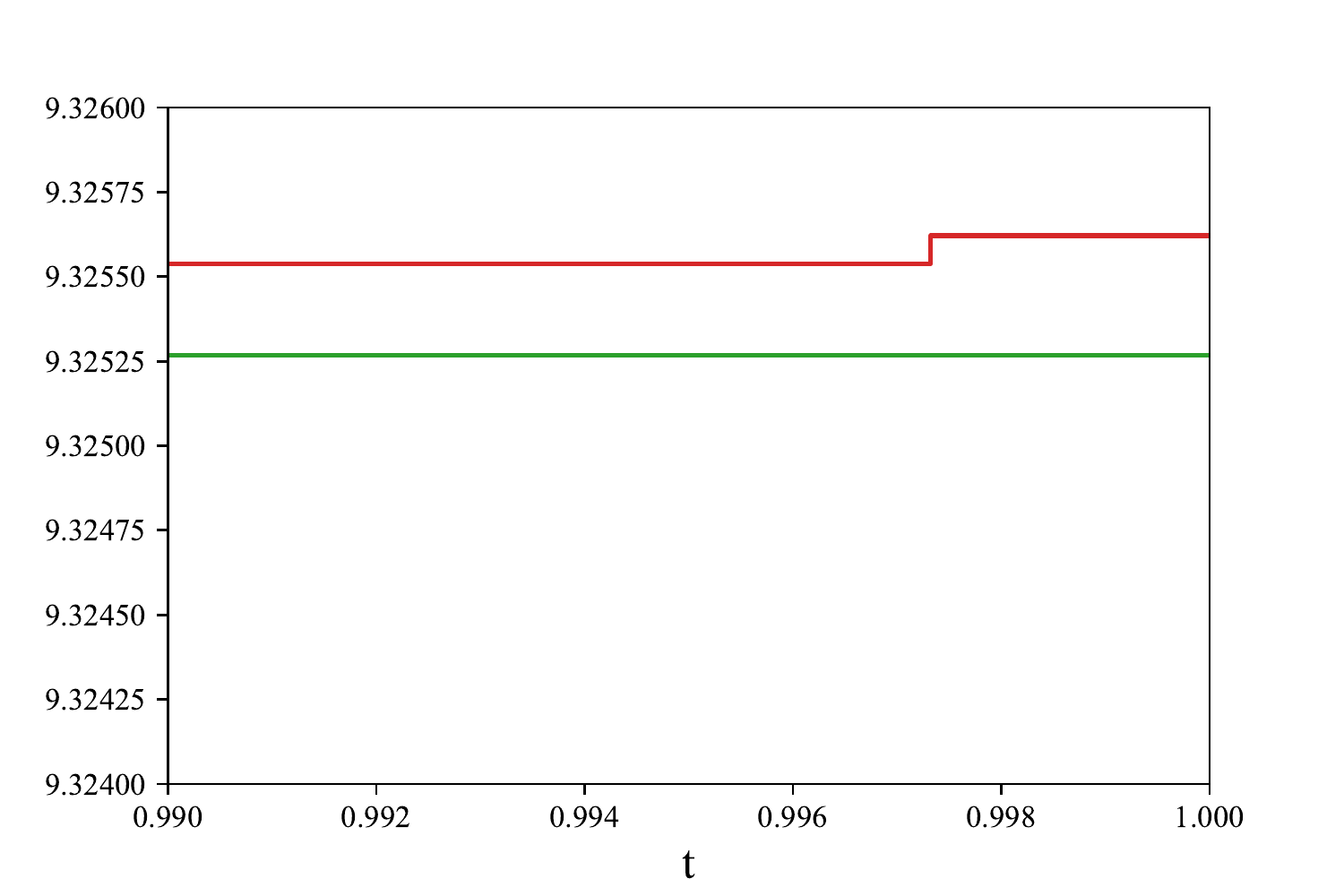} \\
		(a) & (b)
	\end{tabular}
	\caption{Sample paths of the gamma process with $a=10$ and $\beta=1$ generated by truncation of Bondesson's shot noise representation \eqref{Bondesson's method for gamma law} over a common probability space. The truncation parameters are $n=10,50,100,1000$. As the summands of the representation \eqref{Bondesson's method for gamma law} are positive, the sample paths are increasing with $n$. The estimates for the unit-time marginal $X_1$ (that is, the rightmost levels) are 6.02889, 9.26361, 9.32527 and 9.32562 (to 5 decimal places), respectively. As the latter sample paths are difficult to distinguish in plot (a), the plot (b) is zoomed near the terminal time in so as to visualise the slight gap between the sample paths for $n=100$ (lower) and $n=1000$ (higher).}
	\label{fig gamma process}
\end{figure}

In what follows, we provide numerical examples of approximate sample paths for some representative L\'evy processes of both theoretical and practical interests. We first demonstrate the truncation scheme for the gamma process in Figure \ref{fig gamma process} based on Bondesson's method described in Example \ref{example gamma law} (iv), with different truncation parameters to illustrate the convergence and error of the Poisson truncation method.
The estimates for the unit-time marginal in Figure \ref{fig gamma process} suggest that the approximate sample paths based on the truncation of Bondesson's shot noise representation for the gamma process \eqref{Bondesson's method for gamma law} converge very fast in $n$.
This has already been verified by Example \ref{example gamma truncation error} (the $k=4$ case), in which the mean and variance of the truncation error converge to zero exponentially fast. Thus, despite the inapplicability of Gaussian approximation of the discarded jumps for the gamma process (Section \ref{section Gaussian approximation}), the exponential mean-squared convergence of this method makes such accuracy considerations superfluous.

\begin{figure}[h]\centering
	\begin{tabular}{ccc}
		\includegraphics[width=61mm]{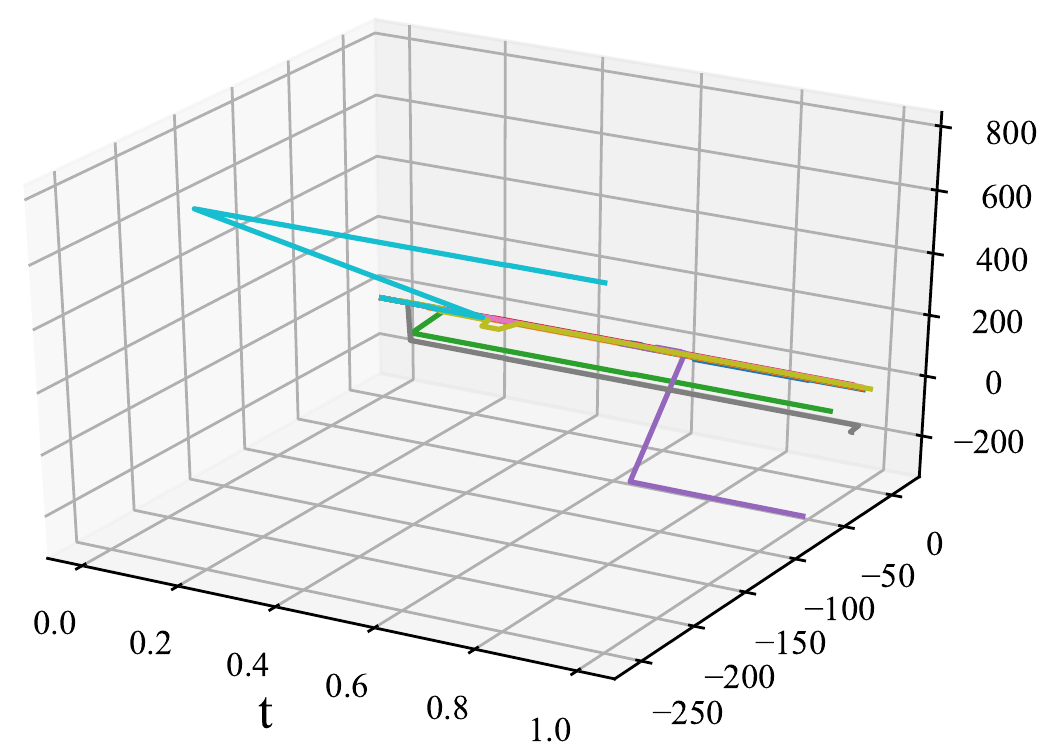} & \includegraphics[width=61mm]{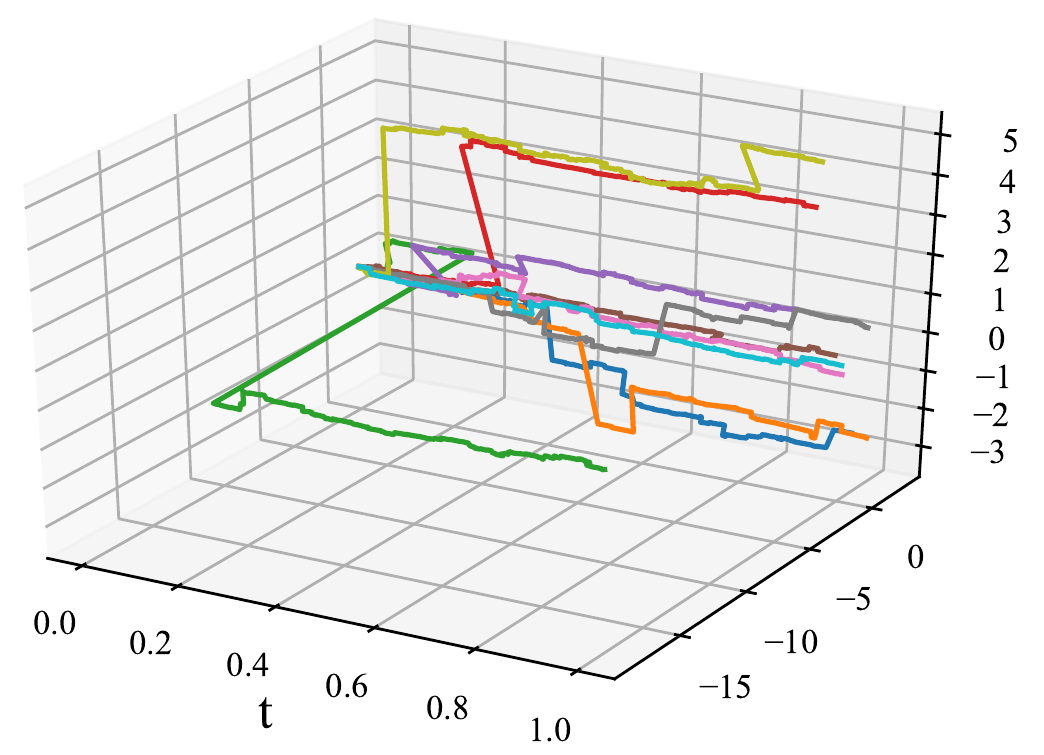} & \includegraphics[width=61mm]{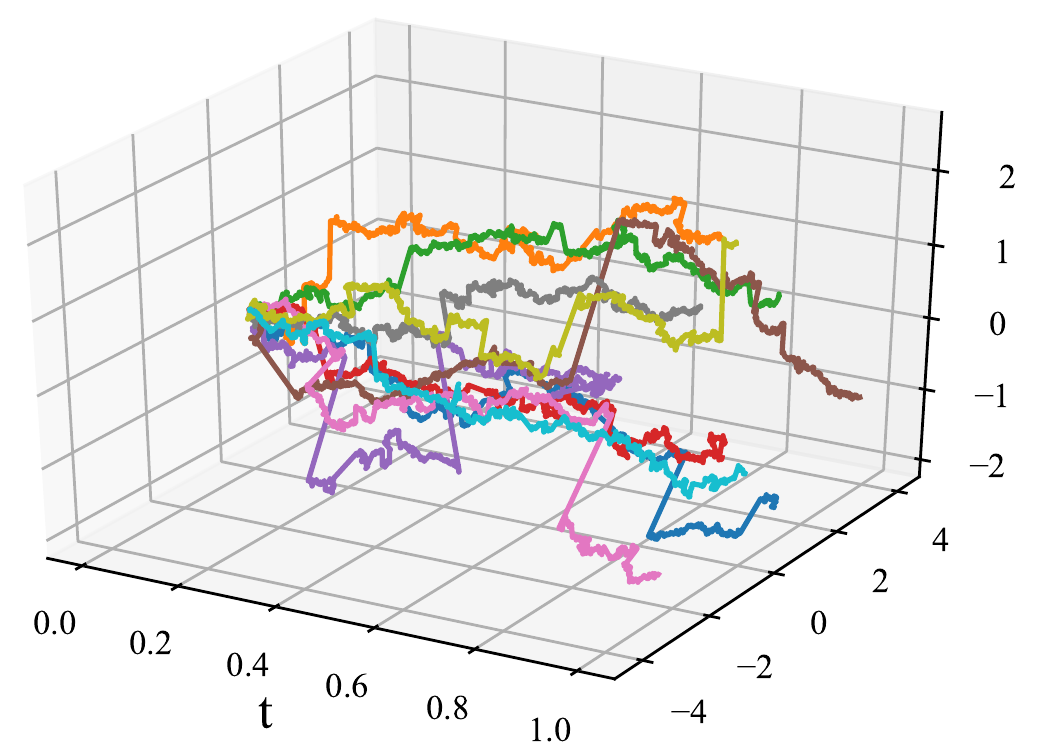} \\
		(a) $\alpha=0.5$ & (b) $\alpha=1$ & (c) $\alpha=1.5$
	\end{tabular}
	\caption{Examples of approximate sample paths of the 2-dimensional stable process with the isotropic L\'evy measure based on the truncation of the shot noise representation \eqref{series rep of stable processes}. The other parameters are $a=1$, $n=10^5$ and $T=1$. The mean-squared errors \eqref{mse of truncation of stable} are (a) $5.33\times10^{-15}$, (b) $10^{-5}$ and (c) $3.76\times10^{-2}$ (3 significant figures). The smallest jump magnitudes included in the truncations are (a) $4\times10^{-10}$, (b) $10^{-5}$ and (c) $3.54\times10^{-4}$ (3 significant figures). Each plot contains 10 iid sample paths.}
	\label{fig 2D stable paths}
\end{figure}

Recall that the shot noise representation of the stable law is provided in Example \ref{example stable series representation}, and the discussion of the mean-squared error of its truncation is given in Example \ref{example stable error analysis}. We provide sample paths of the 2-dimensional stable process with the isotropic L\'evy measure in Figure \ref{fig 2D stable paths}. From our numerical illustration, we see that lower values of the stability parameter correspond to larger jumps (Figure \ref{fig 2D stable paths} (a)), while higher values of the stability parameter correspond to smaller jumps and resemble closer to a Brownian motion (Figure \ref{fig 2D stable paths} (c)). In contrast to the case of the gamma process, the mean-squared convergence remains slower than exponential. As such, Gaussian approximation of the discarded jumps (Section \ref{section Gaussian approximation}) can be practical for enhancing the accuracy of the simulation in this case.

We provide sample paths of the tempered stable process with the isotropic L\'evy measure based on Rosi\'nski's series representation \eqref{eq tempered stable law series representation} in Figure \ref{fig 2D tempered stable paths} below. Comparing with Figure \ref{fig 2D stable paths}, we see that the jumps of the tempered stable sample paths tend to be smaller, which is expected from the exponential tempering of the L\'evy measure. This is most prominent for $\alpha=0.5$, where the stable sample paths (Figure \ref{fig 2D stable paths} (a)) see jumps with magnitudes easily exceeding $100$, while the tempered stable counterparts (Figure \ref{fig 2D tempered stable paths} (a)) do not observe jumps with magnitudes greater than one, due to the random truncation $W_kU_k^{1/\alpha}\|V_k\|$ in every summand of \eqref{eq tempered stable law series representation}.

\begin{figure}[h]\centering
	\begin{tabular}{ccc}
		\includegraphics[width=61mm]{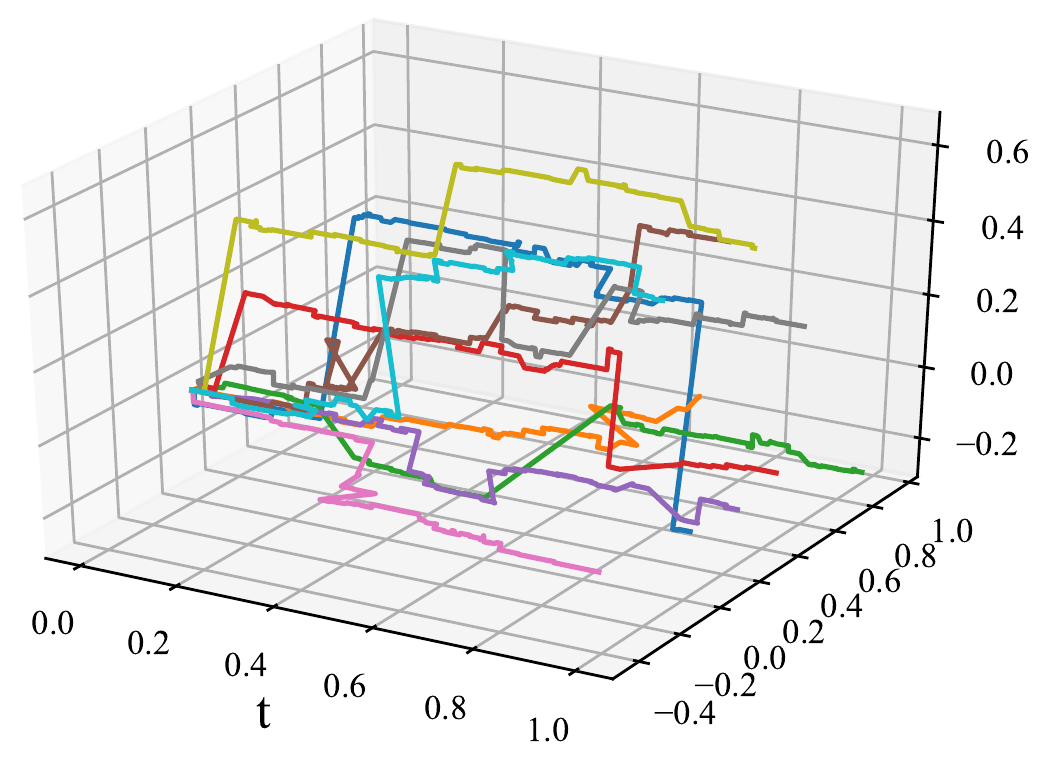} & \includegraphics[width=61mm]{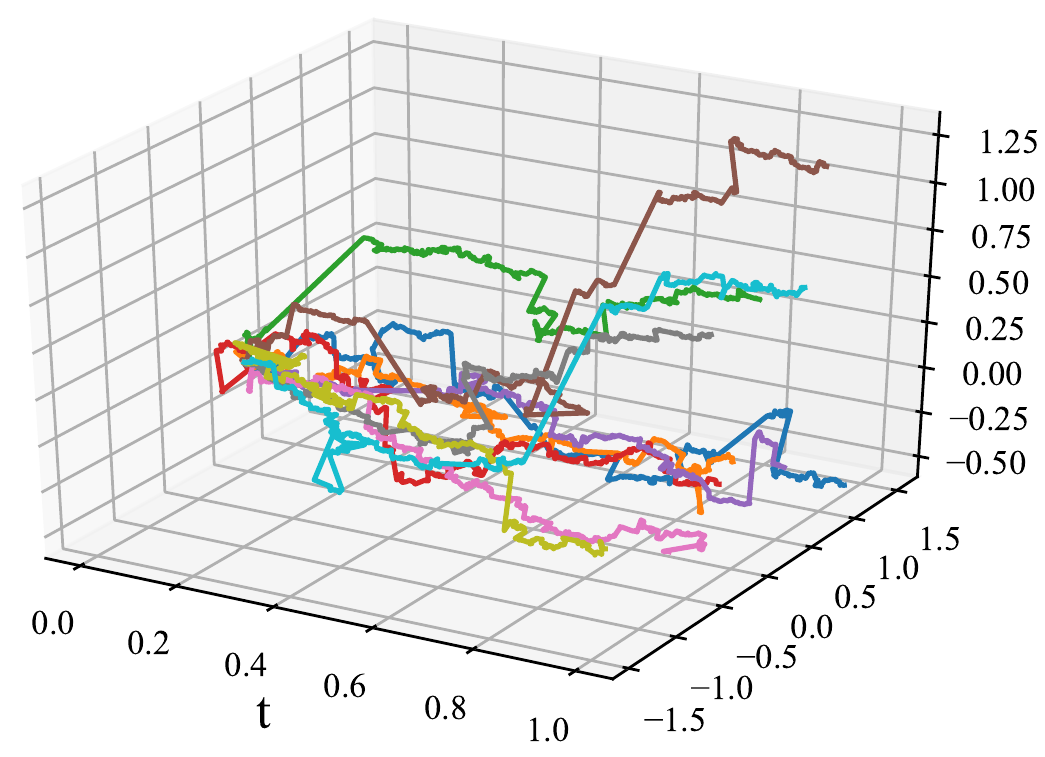} & \includegraphics[width=61mm]{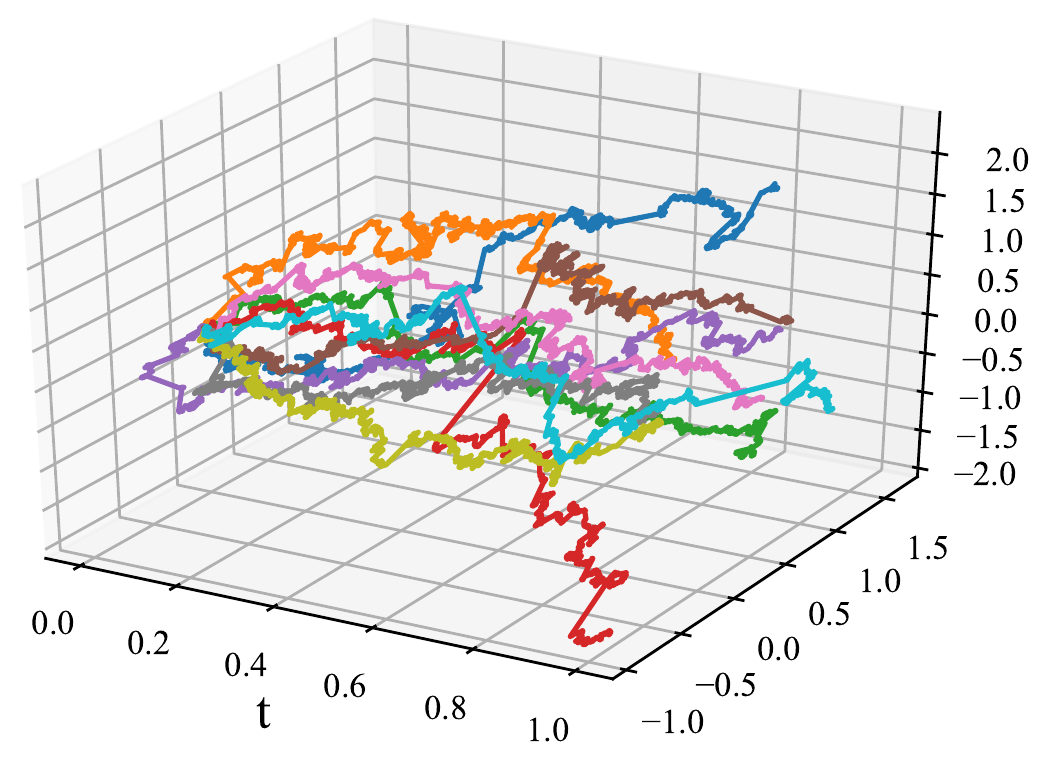} \\
		(a) $\alpha=0.5$ & (b) $\alpha=1$ & (c) $\alpha=1.5$
	\end{tabular}
	\caption{Examples of approximate sample paths of the 2-dimensional tempered stable process with the isotropic L\'evy measure based on the truncation of Rosi\'nski's series representation \eqref{eq tempered stable law series representation}. The other parameters are $a=1$, $n=10^5$ and $T=1$. Each plot contains 10 iid sample paths.}
	\label{fig 2D tempered stable paths}
\end{figure}

We mention here that for the stable process, allowing the stability parameter $\alpha:[0,+\infty) \to (0,2)$ to vary with time leads to the multistable process. Shot noise representations for the multistable processes can be found in \cite{lecourtois2018some,leguevel2012ferguson}. Defined similarly is the tempered multistable process, for which the CGMY process \cite{carr2002fine} with L\'evy measure \eqref{CGMY Levy measure} is an example of. Shot noise representations for tempered multistable processes are provided in \cite{lecourtois2018some}.
Another example of representing L\'evy processes via their shot noise series is the case of $t$-distributed increments \cite{massing2018simulation}.

\subsubsection{Discussion}
We discuss some aspects of simulation by Poisson truncation of shot noise representations.
To begin with, we remark that simulating L\'evy processes based on the thinning method (Proposition \ref{prop series representations} (iii)) may be computationally taxing for obtaining a large number of jumps. This is due to the decreasing acceptance probability as the summation index $k$ increases. For example, in the case of the thinning method for the gamma process (Example \ref{example gamma law}), the acceptance probability for taking the $k$-th summand as a jump is given by $\mathbb P(\Gamma_k V_k < \alpha)$.
As a result, significantly more computations are required to obtain a large quantity of jumps compared to, say, the inverse L\'evy measure method. This explains the absence of exponential convergence rates of the thinning method \eqref{thinning method for gamma law} in the case of the gamma process (Example \ref{example gamma truncation error}), that are enjoyed by the other methods \eqref{inverse Levy measure method for gamma law}, \eqref{rejection method for gamma law} and \eqref{Bondesson's method for gamma law}.
Nevertheless, the thinning method plays a crucial role from a theoretical point of view. For instance, it is used in \cite{talagrand1993regularity} to discuss the boundedness of infinitely divisible processes. As such, different shot noise representations have different potential and use.
For example, while the truncation of the shot noise representation by the inverse L\'evy measure method simulates a L\'evy process by discarding its smallest jumps, the truncation of a representation by the rejection method illustrates the relationship between a L\'evy process with another.

The Poisson truncation approximation of a L\'evy process is able to preserve various key properties from the original process.
For example, the discontinuity of sample paths clearly hold for both the original L\'evy process and its Poisson truncation approximation.
As the largest jumps of a L\'evy process without Gaussian components account for the majority of the variation, some key moment properties are retained by the Poisson truncation approximation.
In the case of the tempered stable and gamma processes, the sample paths resulting from Poisson truncation retain the finiteness of moments of all polynomial orders. Meanwhile, for the stable process, marginal moments of order $[\alpha,+\infty)$ remain infinite even after Poisson truncation. Similarly, the extremal behaviour remains unchanged after Poisson truncation as it is attributed to the largest jumps \cite{hult2007extremal}.

It is well-known that the finiteness of the total variation of a L\'evy process without Gaussian components with L\'evy measure $\nu(d\mathbf z)$ depends on the finiteness of the integral $\int_{\mathbb R_0^d}(\|\mathbf z\|\land1)\,\nu(d\mathbf z)$.
In particular, a sufficiently intense activity of small jumps is the only possible source of infinite variation.
Truncation of shot noise representation cuts off all small jumps and thus necessarily leads to sample paths of finite variation, so this property is preserved for finite variation L\'evy processes. However, sample paths of infinite variation, such as that of stable and tempered stable processes with stability $\alpha\in[1,2)$, will result in finite variation after applying Poisson truncation. We mention that in this case, the total variation of the Poisson truncation approximation diverges fastest in the case of the inverse L\'evy measure method, echoing the dominance we saw in Theorem \ref{theorem finite truncation comparison} (iii) but with the integrand $\|\mathbf z\|\land1$ instead. The lack of preservation of infinite variation means that even though the truncation of shot noise representation leads to sample paths based on individual jumps, the resulting paths cannot be employed in full for investigating sample path properties. From the perspective of simulation, this is not an issue, as the compromise of generating finite variation approximations of the L\'evy process is implied by the very nature of numerical investigation.

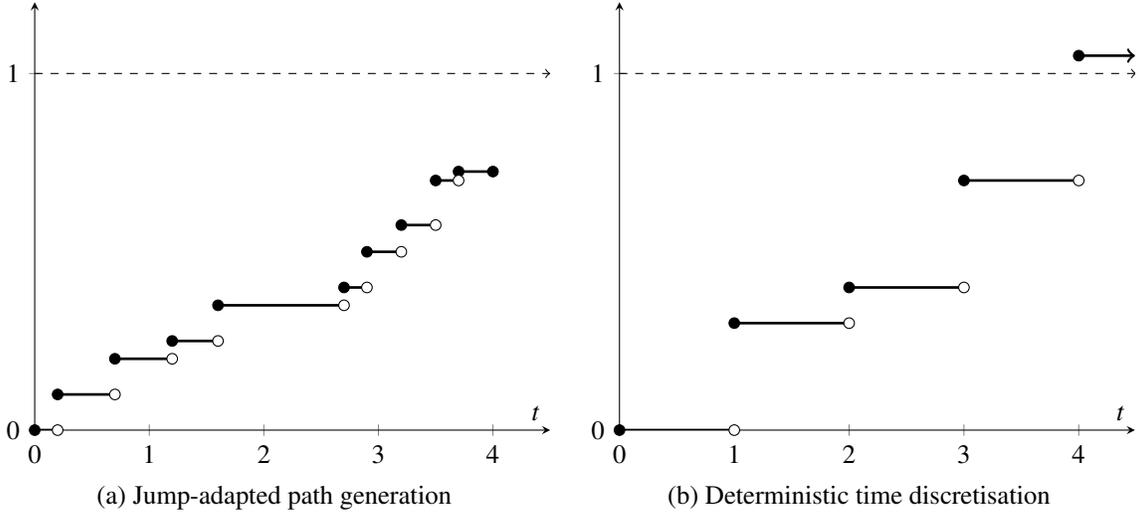
\begin{figure}[h]
	\centering
	\begin{tabular}{cc}
		\begin{tikzpicture}
		\begin{axis}[%
		,xlabel=$t$
		,axis x line = bottom,axis y line = left
		,xtick={0,1,2,3,4}
		,ytick={0,1}
		,ymax=1.2 
		,xmax=4.5
		]
		\addplot[cmhplot,-,domain=0:0.2]{0};
		\addplot[cmhplot,-,domain=0.2:0.7]{0.1};
		\addplot[cmhplot,-,domain=0.7:1.2]{0.2};
		\addplot[cmhplot,-,domain=1.2:1.6]{0.25};
		\addplot[cmhplot,-,domain=1.6:2.7]{0.35};
		\addplot[cmhplot,-,domain=2.7:2.9]{0.40};
		\addplot[cmhplot,-,domain=2.9:3.2]{0.50};
		\addplot[cmhplot,-,domain=3.2:3.5]{0.575};
		\addplot[cmhplot,-,domain=3.5:3.7]{0.70};
		\addplot[cmhplot,-,domain=3.7:4]{0.725};
		\addplot[soldot]coordinates{(0,0) (0.2,0.1)(0.7,0.2) (1.2,0.25) (1.6,0.35) (2.7,0.40) (2.9,0.50) (3.2,0.575) (3.5,0.70) (3.7,0.725) (4,0.725)};
		\addplot[holdot]coordinates{(0.2,0) (0.7,0.1) (1.2,0.2) (1.6,0.25) (2.7,0.35) (2.9,0.40) (3.2,0.50) (3.5,0.575) (3.7,0.70)};
		\draw [dashed,->] (0,1) -- (4.5,1);
		\end{axis}
		\end{tikzpicture}
		&
		\begin{tikzpicture}
		\begin{axis}[%
		,xlabel=$t$
		,axis x line = bottom,axis y line = left
		,xtick={0,1,2,3,4}
		,ytick={0,1}
		,ymax=1.2 
		]
		\addplot[cmhplot,-,domain=0:1]{0};
		\addplot[cmhplot,-,domain=1:2]{0.3};
		\addplot[cmhplot,-,domain=2:3]{0.4};
		\addplot[cmhplot,-,domain=3:4]{0.7};
		\addplot[cmhplot,->,domain=4:4.49]{1.05};
		\addplot[soldot]coordinates{(0,0) (1,0.3) (2,0.4) (3,0.7) (4,1.05)};
		\addplot[holdot]coordinates{(1,0) (2,0.3) (3,0.4) (4,0.7)};
		\draw [dashed,->] (0,1) -- (4.5,1);
		\end{axis}
		\end{tikzpicture}\\
		(a) Jump-adapted path generation & (b) Deterministic time discretisation
	\end{tabular}
	\caption{Comparison between generating sample paths to estimate the first crossing time $\inf\{t>0 : X_t>1\}$ by simulation via (a) jumps and (b) increments. The Poisson truncation method corresponds to (a), with $T=4$. The particular realisation does not reach the unit threshold, thus one would need to rerun the simulation with larger $T$. For simulation via unit increments (b), one could simply keep piling on increments until the unit threshold is reached.}
	\label{threshold hitting}
\end{figure}

Before closing this subsection, we mention a possible disadvantage of the truncation method for sample path generation of L\'evy processes. When sample path generation via increments is possible, one important distinction between such a scheme and that of a Poisson truncation method is that in the latter, a terminal time $T$ is fixed beforehand and cannot be extended during the simulation, whereas in the former, one can keep piling on as many increments as desired, possibly until some condition is fulfilled. An example of when this distinction is important is through a comparison between the settings of \cite{yuan2020asymptotic} and \cite{carnaffan2017solving}. In the former, Poisson truncation of shot noise representation is used to generate sample paths of an inverse subordinator. In that setting, the goal of evaluating immobility periods requires the observation of jumps, thus rendering the truncation method as ideal and the incremental method as inappropriate. In \cite{carnaffan2017solving}, however, evaluations of the inverse subordinator at any time $t>0$ is required, so if the subordinator (before inversion) has yet to reach $t$, one must keep adding increments until that level is reached. This is illustrated in Figure \ref{threshold hitting}.
Thus, in such a scenario, sample path generation by increments should be preferred.

\subsection{Simulating infinitely divisible processes}\label{section Levy-driven stochastic integrals}
So far, we have studied approximations of L\'evy processes via truncation of shot noise representations. We will now generalise the technique further to approximate stochastic processes governed by L\'evy-driven stochastic integrals of the form \eqref{stochastic integral process}, which is a large class of infinitely divisible processes. We emphasise that our focus is on the case where the driving L\'evy process is of infinite jump activity, as otherwise, exact simulation methods are readily available. The technique revolves around the partitioning the stochastic process based on the the sizes and timings of the underlying jumps, as we will present shortly in \eqref{stochastic integral decomposition}. Shot noise representation is perhaps even more pertinent to the construction and approximation of infinitely divisible processes, which often requires the consideration of individual jumps more so than L\'evy processes for both theoretical and numerical purposes.

Suppose we want to simulate the stochastic process $\{X_t:t\in[0,T]\}$ in $\mathbb R^d$ such that the marginal is described by the stochastic integral
$X_t = \int_{\mathcal T}f(t,s)\,dL_s$,
where $\mathcal T\subseteq \mathbb R$ and $\{L_s:s\in\mathcal T\}$ is a L\'evy process with an infinite L\'evy measure $\nu(d\mathbf z)$ with a decomposition \eqref{generalised shot noise method decomposition}. This corresponds to the stochastic integral form \eqref{stochastic integral process} where the integrator is a L\'evy process. Naturally, each jump of the underlying L\'evy process at time $s$ is modulated by $s\mapsto f(t,s)$, so we have a L\'evy-It\^o decomposition of the form
\begin{equation}\label{Levy-Ito decomposition of ID process}
	X_t = \int_{\mathcal T}\int_{\mathbb R_0^d} f(t,s)\mathbf z(\mu(d\mathbf z,ds) - \nu(d\mathbf z)ds\mathbbm1_{(0,1]}(\|\mathbf z\|)),
\end{equation}
where $\mu(d\mathbf z,ds)$ is the Poisson random measure on $\mathbb R_0^d\times\mathcal T$ associated with $\{L_s:s\in\mathcal T\}$ and $\nu(d\mathbf z)ds$ is the corresponding compensator.
Various theoretical developments for infinitely divisible processes and their series representations have been established, such as
their spectral representations \cite{rajput1989spectral} and path properties \cite{rosinski1989path,talagrand1993regularity}. Much like in the case of L\'evy processes, shot noise representations for infinitely divisible processes converge almost surely uniformly under suitable technical conditions \cite{basse2013uniform}. Moreover, if the probability space is rich enough, then one can choose the random sequences such that the shot noise representation is almost surely equal to the infinitely divisible process \cite{rosinski2018representations}.

We look to truncate the infinitely divisible process $\{X_t:t\in[0,T]\}$\eqref{Levy-Ito decomposition of ID process} based on jump timings and sizes of the underlying L\'evy process $\{L_s:s\in\mathcal T\}$ via shot noise representation in order to obtain independent simulatable and error components for analysis. If $\text{Leb}(\mathcal T) < +\infty$, that is, when the domain of the jump timings is bounded, we require no truncation on the index set. Otherwise, we introduce a nondecreasing sequence $\{\mathcal T_n\}_{n\in\mathbb N}$ of connected subintervals of $\mathcal T$ such that $\cup_{n\in\mathbb N}\mathcal T_n = \mathcal T$ and $\text{Leb}(\mathcal T_n) < +\infty$. That is, the parameter $n$ represents a truncation based on time timings.
For the truncation over jump sizes, we note that the finite measure $\nu_m(d\mathbf z) \coloneqq \int_{0+}^{m}\mathbb P(H(r,U)\in d\mathbf z)\,dr$ is the L\'evy measure of the Poisson truncation approximation of the L\'evy process $\{L_s:s\in\mathcal T_n\}$ with respect to the index set $\{k\in\mathbb N : \Gamma_k \le m\text{Leb}(\mathcal T_n)\}$, for every $n\in\mathbb N$. That is, $m\in\mathbb N$ fulfils the role of a finite truncation parameter on the L\'evy measure $\nu(d\mathbf z)$. We remark that for any fixed $m$, the impact of the truncation on the sizes of the surviving jumps depends on the shot noise representation used, in the same vein as our discussion in Section \ref{section analysis of finite truncation}.
We decompose the stochastic process as \cite{idsim}
\begin{equation}\label{stochastic integral decomposition}
	X_t = X_t(m,n) + Q_t(m) + R_t(m,n), \quad t\in[0,T],
\end{equation}
where
\begin{align}
	\notag X_t(m,n) &\coloneqq \int_{\mathcal T_n}\int_{\mathbb R_0^{d}}f(t,s)\mathbf z\,(\mu_m - \nu_m\mathbbm1_{(0,1]}(\|\mathbf z\|))(d\mathbf z,ds),\\
	\notag Q_t(m) &\coloneqq \int_{\mathcal T}\int_{\mathbb R_0^{d}}f(t,s)\mathbf z\,((\mu - \mu_m) - (\nu - \nu_m)\mathbbm1_{(0,1]}(\|\mathbf z\|))(d\mathbf z,ds),\\
	\notag R_t(m,n) &\coloneqq \int_{\mathcal T\backslash\mathcal T_n}\int_{\mathbb R_0^{d}}f(t,s)\mathbf z\,(\mu_m - \nu_m\mathbbm1_{(0,1]}(\|\mathbf z\|))(d\mathbf z,ds),
\end{align}
Intuitively, the component $X_t(m,n)$ corresponds to large jumps over $\mathcal T_n$, $Q_t(m)$ corresponds to small jumps over $\mathcal T$, and $R_t(m,n)$ corresponds to large jumps larger over $\mathcal T\backslash\mathcal T_n$. Since the regions of jump sizes and timings simulated by each component form a disjoint union, the stochastic processes on the right hand side of \eqref{stochastic integral decomposition} can be treated independently thanks to the independent scattering of the Poisson random measure $\mu(d\mathbf z,ds)$.

As $\nu_m(\mathbb R_0^d)\text{Leb}(\mathcal T_n)=m\text{Leb}(\mathcal T_n)<+\infty$ for every $(m,n)\in\mathbb N^2$, the former component $\{X_t(m,n):t\in[0,T]\}$ can be treated as the principal simulatable component of the approximation by the Poisson truncation with the shot noise representation
\begin{equation}\label{principal truncation}
	\{X_t(m,n):t\in[0,T]\} \stackrel{\mathscr L}{=} \left\{\sum_{\{k\in\mathbb N:\Gamma_k\le \text{Leb}(\mathcal T_n)m\}}\left[f(t,T_k)H\left(\frac{\Gamma_k}{\text{Leb}(\mathcal T_n)},U_k\right) - c_k\int_{\mathcal T_n}f(t,s)\,ds\right] : t\in[0,T]\right\},
\end{equation}
where $\{T_k\}_{k\in\mathbb N}$ is a sequence of iid uniform random variables on $\mathcal T_n$, and the other random sequences are as in \eqref{generalised shot noise series}.
For generating sample paths of the infinitely divisible process based on the Poisson truncation approximation \eqref{principal truncation} over sample points $0=t_0<t_1<\cdots<t_{J-1}<t_J=T$, we provide the numerical recipe as follows:

\begin{enumerate}[\bf Step 1.]
	\setlength{\parskip}{0cm}
	\setlength{\itemsep}{0cm}
	\item Generate a standard exponential random variable $E_1$. If $E_1 \le m\text{Leb}(\mathcal T_n)$, then assign $\Gamma_1 \leftarrow E_1$. Otherwise, return the degenerate zero process as the approximate sample path and terminate the algorithm.
	
	\item While $\Gamma_k \le m\text{Leb}(\mathcal T_n)$, generate a standard exponential random variable $E_{k+1}$ and assign $\Gamma_{k+1} \leftarrow \Gamma_k + E_{k+1}$. Denote this as $\{\Gamma_k\}_{k\in\{1,\cdots,N\}}$, where $N$ satisfies $\Gamma_N \le m\text{Leb}(\mathcal T_n) < \Gamma_{N+1}$.
	
	\item Generate a sequence $\{T_k\}_{k\in\{1,\cdots,N\}}$ of iid uniform random variables on $\mathcal T_n$.
	
	\item Generate a sequence $\{U_k\}_{k\in\{1,\cdots,N\}}$ of suitable iid random variables.
	
	\item For every $j\in\{1,\cdots,J\}$, assign $L_{t_j}^{H,\alpha,n} \leftarrow \sum_{k=1}^{N} (f_n(t_j,T_k)H(\Gamma_k/\text{Leb}(\mathcal T_n),U_k) - c_k\int_{\mathcal T_n}f(t_j,s)\,ds)$.
	
	\item Return $\{0,L_{t_1}^{H,\alpha,n},\cdots,L_{t_J}^{H,\alpha,n}\}$ as the positions of the approximate sample path at the sample times $\{0,t_1,\cdots,t_J\}$.
\end{enumerate}

We can consider the residual components $\{Q_t(m):t\in[0,T]\}$ and $\{R_t(m,n):t\in[0,T]\}$ as the error processes to analyse. In particular, if the kernel $f(t,\cdot)$ is square-integrable and essentially uniformly bounded for every $t\in[0,T]$, and the kernel of the shot noise representation $H(\cdot,\bm\xi)$ satisfies some suitable technical conditions, then the stochastic process $\{Q_t(m):t\in[0,T]\}$ comprising of small jumps can be approximated by a Gaussian process in the same vein as Section \ref{section Gaussian approximation}. For more details, we refer the reader to \cite{idsim}.

In what follows, we demonstrate this approximate sample path generation method along with error analysis using the examples of higher order fractional stable motion (Section \ref{section nfsm}) and L\'evy-driven continuous-time autoregressive moving average (CARMA) processes (Section \ref{section CARMA process}). In both of those cases, for simplicity we restrict ourselves to the univariate setting with $\mathcal T = \mathbb R$, and typically consider the Poisson truncation based on the inverse L\'evy measure method, so $\nu_m(d\mathbf z) = \mathbbm1_{(\eta(m),+\infty)}(\|\mathbf z\|)\nu(d\mathbf z)$, where $\eta(m) = \sup\{r>0 : \int_{\|\mathbf z\|>r}\nu(d\mathbf z) > m\}$.
We also discuss L\'evy-driven Ornstein-Uhlenbeck processes as interesting cases in which simulation by increments may be preferred over truncation of shot noise representation (Section \ref{section OU process}).

\subsubsection{Fractional stable motion}\label{section nfsm}
The higher order fractional stable motion, which generalises the higher order fractional Brownian motion, was defined in \cite{kawai2016higher}. Let $\{L_t^{\alpha,+} : t\ge0\}$ be a stable process in $\mathbb R$ with stability parameter $\alpha\in(0,2)$ and a symmetric L\'evy measure 
\[
	\nu(dz) \coloneqq\frac{\alpha c_\alpha}{2}\frac{1}{|z|^{\alpha+1}}\,dz,
\]
where $c_\alpha\coloneqq(\Gamma(1-\alpha)\cos(\pi\alpha/2))^{-1}$ for $\alpha\in(0,2)\backslash\{1\}$ and $c_1 \coloneqq 2/\pi$. We extend the stable process temporally by defining another stable process $L^{\alpha,-}$ similarly such that $L^{\alpha,-}\stackrel{\mathscr L}{=}L^{\alpha,+}$, and defining $\{L_t^\alpha:t\in\mathbb R\}$, where $L_t^\alpha\coloneqq L_t^{\alpha,+}\mathbbm1(t\ge0) + L_{-t-}^{\alpha,-}\mathbbm1(t<0)$.
Let $n\in\mathbb N$ and $H\in(n-1,n)\backslash\{1/\alpha\}$ such that $H\notin\mathbb Z$. We define the \textit{$n$-th order fractional stable motion with Hurst parameter $H$} by $\{L_t^{H,\alpha,n}:t\in\mathbb R\}$, where
\begin{equation}\label{nfsm}
	L_t^{H,\alpha,n}\coloneqq \int_{\mathbb R}f_n(t,s;H,\alpha)\,dL_s^\alpha,\quad t\in\mathbb R,
\end{equation}
with the higher order moving average kernel
\[
	f_n(t,s;H,\alpha)\coloneqq\frac{1}{\Gamma(H-1/\alpha + 1)}\left((t-s)_+^{H-1/\alpha} - \sum_{k=0}^{n-1}\binom{H-1/\alpha}{k}(-s)_+^{H-1/\alpha-k}t^k\right).
\]
As $\alpha\to2$, $L^{H,\alpha,n}$ converges towards an $n$-th order fractional Brownian motion. We call $L^{H,\alpha,1}$ a \textit{linear fractional stable motion}.

Using the framework of the decomposition \eqref{stochastic integral decomposition} based on the inverse L\'evy measure method with the truncation on jump timings $\mathcal T_\kappa \coloneqq (\kappa,T]$ (where we denote the truncation parameter as $\kappa$ instead of $n$ to avoid notational conflict), we decompose the higher order fractional stable motion as
\begin{equation}\label{higher order fBm decomposition}
	L_t^{H,\alpha,n} = L_t^{H,\alpha,n}(m,\kappa) + Q_t(m) + R_t(m,\kappa).
\end{equation}
The following shot noise representation of $\{L_t^{H,\alpha,n}(m,\kappa):t\in[0,T]\}$ offers a clear solution to its simulation. Let $\{\Gamma_k\}_{k\in\mathbb N}$ be a sequence of standard Poisson arrival times, $\{T_k\}_{k\in\mathbb N}$ a sequence of iid uniform random variables on $(-\kappa,T)$, $\{U_k\}_{k\in\mathbb N}$ a sequence of iid random variables uniformly distributed on $S^0$, such that all random sequences are mutually independent. Then, we have that \cite[Lemma 6.1]{kawai2016higher}
\begin{equation}\label{nfsm truncation}
	L^{H,\alpha,n}(m,\kappa)\stackrel{\mathscr L}{=} \left\lbrace\sum_{\{k\in\mathbb N: \Gamma_k \le m(T+\kappa)\}} \left(\frac{\Gamma_k}{(T+\kappa)c_\alpha}\right)^{-1/\alpha}U_k f_n(t,T_k;H,\alpha):t\in[0,T]\right\rbrace.
\end{equation}

For generating sample paths of the higher order fractional stable motion based on the Poisson truncation approximation \eqref{nfsm truncation}, we follow the numerical recipe described previously in Section \ref{section Levy-driven stochastic integrals} but with the shot noise representation \eqref{nfsm truncation} instead. Examples of this sample path generation scheme for the higher order fractional motion are provided in Figure \ref{fig nfsm} below. We observe the versatility of the higher order fractional stable motion for generating rough paths (Figure \ref{fig nfsm} (a)) as well as aggregated smooth paths (Figure \ref{fig nfsm} (c)).

\begin{figure}[h]\centering
	\begin{tabular}{ccc}
		\includegraphics[width=61mm]{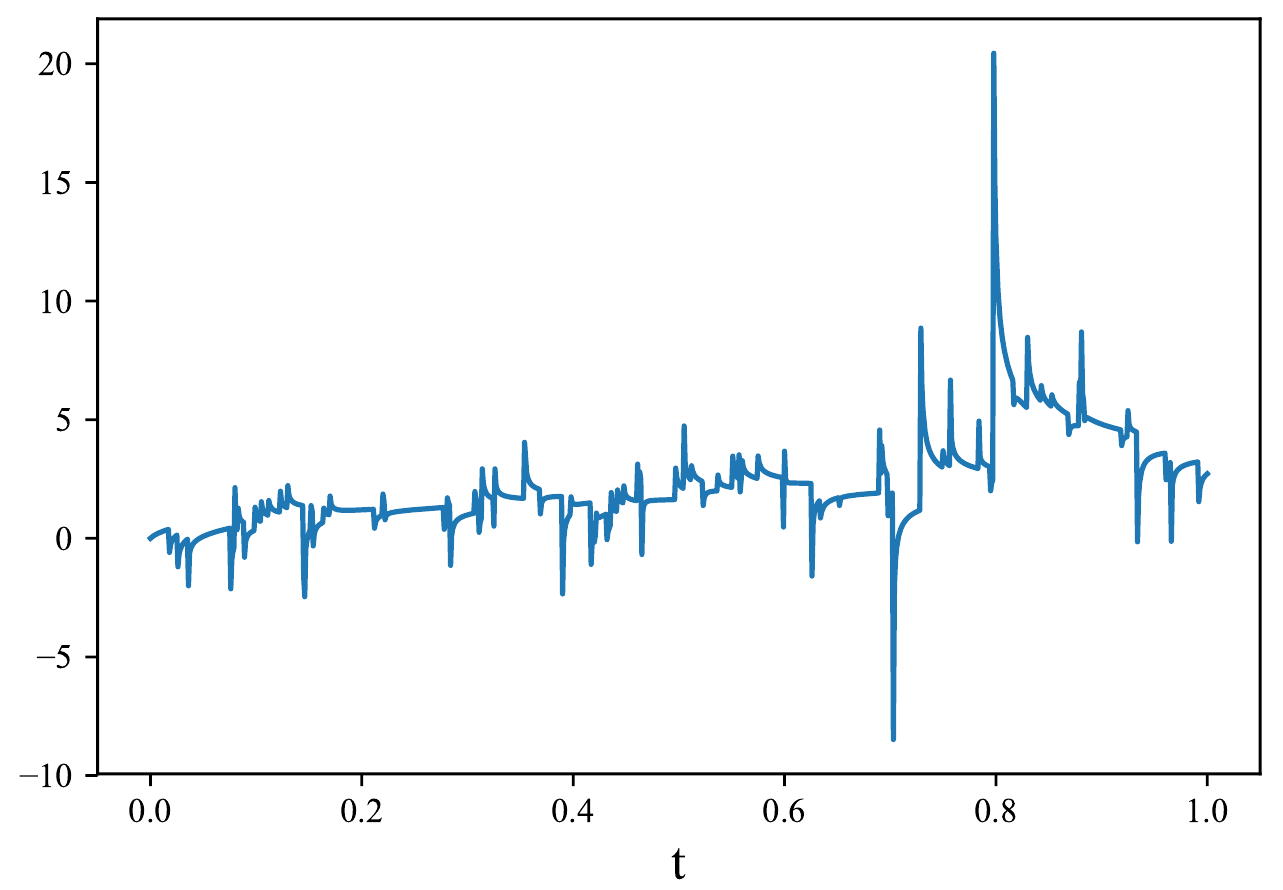} & \includegraphics[width=61mm]{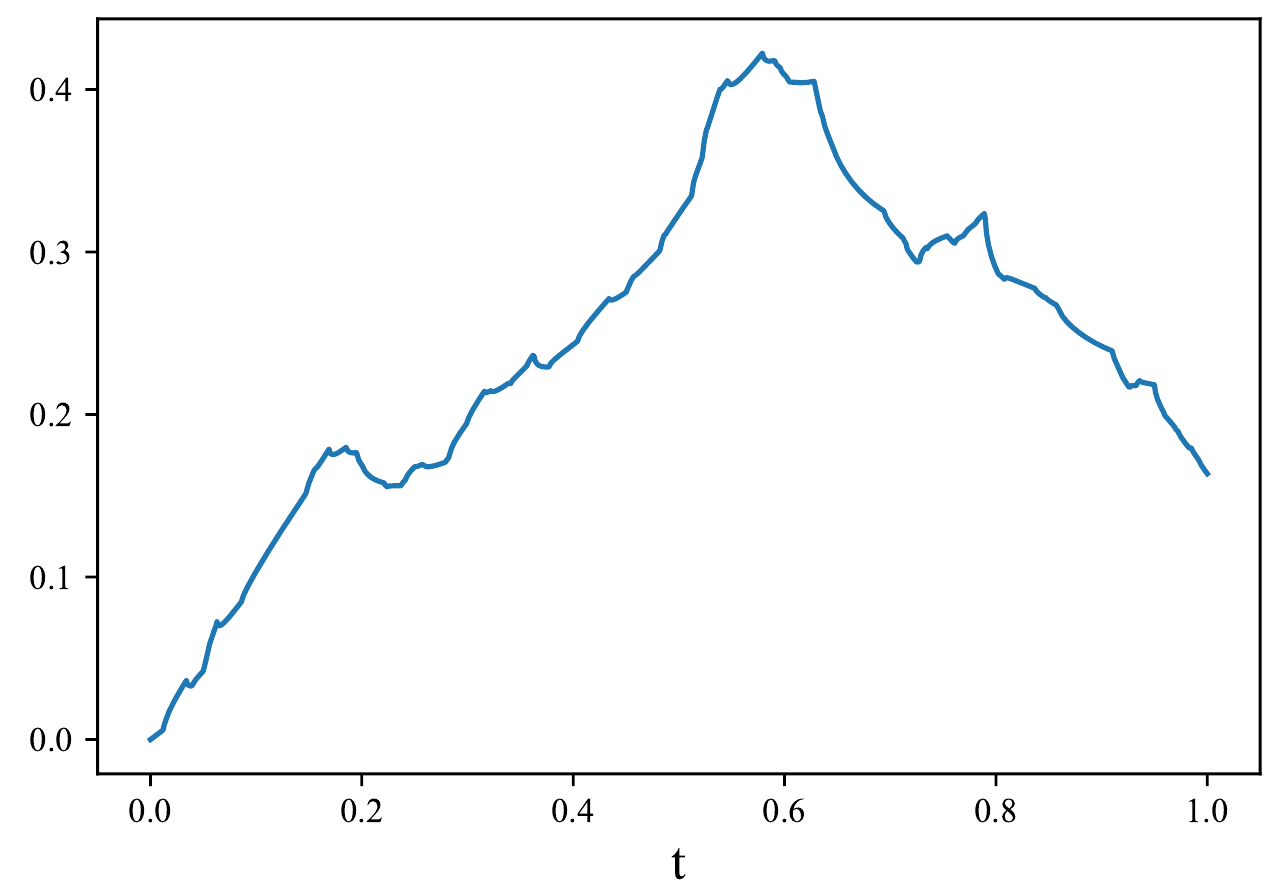} & \includegraphics[width=61mm]{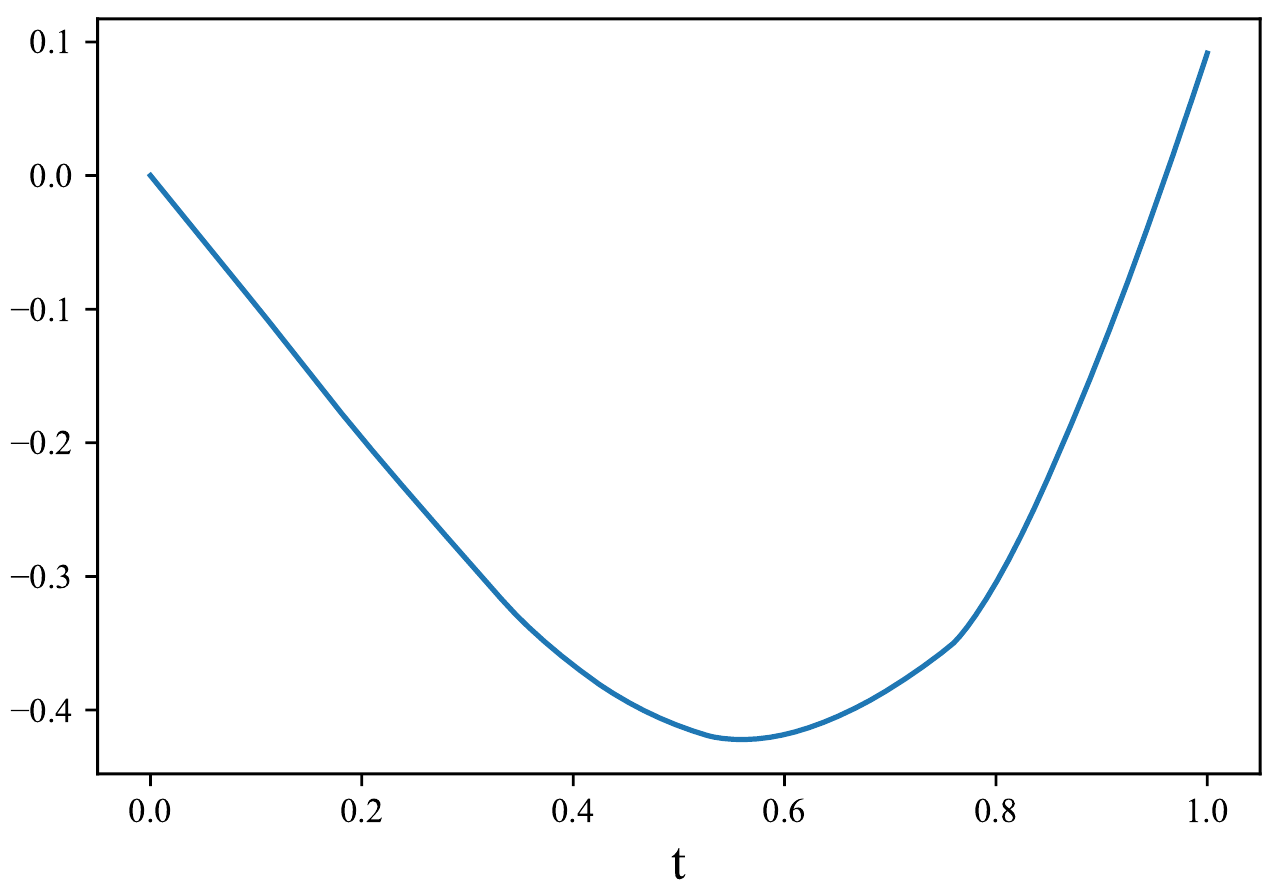} \\
		(a) First-order, $H=0.1$ & (b) Second-order, $H=1.1$ & (c) Third-order, $H=2.1$ 
	\end{tabular}
	\caption{Sample paths of first, second and third-order fractional stable motions via the truncation approximation \eqref{nfsm truncation}. The other parameters are $\alpha=1.7$, $\kappa=5$ and $m=100$. Time steps are of size $10^{-3}$.}
	\label{fig nfsm}
\end{figure}

We remark on the scenario where the sample paths of the stochastic integral process are almost surely unbounded on every finite interval of positive length, which happens whenever the integrand of the stochastic integral is explosive. For instance, in the case of the higher order fractional stable motion \eqref{nfsm}, this occurs when $\alpha\in(0,2\wedge (1/(n-1)))$ and $H\in(n-1,n\wedge(1/\alpha))$ \cite[Theorem 4.1]{kawai2016higher}. In such a case, it is nonsensical to use truncation of shot noise representation to simulate the stochastic integral process, and even misleading if done so without the knowledge of sample path unboundedness. This is because truncation to a finite L\'evy measure will almost surely produce a bounded sample path, which would otherwise be unbounded in the absence of truncation. For this reason, we advise to check that the stochastic process is almost surely bounded over $[0,T]$ prior to generating its sample paths.

We now turn our attention to the stochastic process $\{Q_t(m) : t\in[0,T]\}$ consisting of the small jumps over $(+\infty,T]$. Let us further decompose this stochastic process as
\[
	Q_t(m) = Y_t(m) + Z_t(m),
\]
where
\[
	Y_t(m) \coloneqq \int_0^t\int_{|z|\le\eta(m)}\frac{(t-s)^{H-1/\alpha}}{\Gamma(H - 1/\alpha + 1)}z\, (\mu(dz,ds) - \nu(dz)ds\mathbbm1_{(0,1]}(|z|)),
\]
and
\[
	Z_t(m) \coloneqq \int_{-\infty}^0\int_{|z|\le\eta(m)}f_n(t,s;H,\alpha)z\, (\mu(dz,ds) - \nu(dz)ds\mathbbm1_{(0,1]}(|z|)).
\]
Since the kernel $f(t,s;H,\alpha)$ in $Z_t(m)$ can be written without $(\cdot)_+$, it is continuously differentiable with respect to $t$. The following result gives us the asymptotic behaviour of $Q_t(m)$, and validates its approximation via a fractional Brownian motion \cite[Theorem 6.3]{kawai2016higher}.
Let $\{B_t : t\in\mathbb R\}$ be a temporally extended standard Brownian motion in $\mathbb R$ and denote $\sigma(m)\coloneqq(\int_{|z|\le\eta(m)}|z|^2\,\nu(dz))^{1/2}$.
\begin{enumerate}[(i)]
	\setlength{\parskip}{0cm}
	\setlength{\itemsep}{0cm}
	\item If $\alpha\in(1/n,2)$ and $H\in((n-1)\vee(1/\alpha),n)$, then as $m\to+\infty$, it holds that
	\[
		\left\{\frac{Y_t(m)}{\sigma(m)} : t\in[0,T]\right\} \to \left\{\int_0^T \frac{(t-s)_+^{H-1/\alpha}}{\Gamma(H - 1/\alpha +1)}\,dB_s : t\in[0,T]\right\},
	\]
	where the convergence is weak convergence in $\mathscr C([0,T],\mathbb R)$, the space of continuous functions from the compact interval $[0,T]$ to $\mathbb R$, endowed with the uniform metric.
	\item If $\alpha\in(0,2)$ and $H\in(n-1,n)$, then $\{Z_t(m) : t\in[0,T]\}$ converges in probability to the degenerate zero process uniformly on $[0,T]$ as $m\to+\infty$.
	\item If $\alpha\in(2/3,2)$ and $H\in(n+1/\alpha - 3/2,n)$, then as $m\to+\infty$, it holds that
	\[
		\left\{\frac{Z_t(m)}{\sigma(m)} : t\in[0,T]\right\}\stackrel{\mathscr L}{\to} \left\{\int_{-\infty}^0 f_n(t,s;H,\alpha)\,dB_s : t\in[0,T]\right\}.
	\]
	Moreover, if $n\ge2$, $\alpha\in(2/3,2)$ and $H\in(n+1/\alpha-3/2,n)$, then the weak convergence holds true in $\mathscr C([0,T],\mathbb R)$.
\end{enumerate}

Additionally, by noting that we can write $f_n(t,s;H,\alpha) = f_n(t,s;H - 1/\alpha + 1/2,2)$, the previous result gives us that if $\alpha\in((1/n)\vee(2/3),2)$ and $H\in((n + 1/\alpha - 3/2)\vee(1/\alpha),n)$, then
\[
	\left\lbrace\frac{Q_t(m)}{\sigma(m)} : t\in[0,T]\right\rbrace\stackrel{\mathscr L}{\to}\left\lbrace\int_{-\infty}^T f_n(t,s;H-1/\alpha + 1/2,2)\,dB_s : t\in[0,T]\right\rbrace,
\]
as $m\to+\infty$. That is, we are able to approximate the stochastic process $\{Q_t(m) : t\in[0,T]\}$ with an $n$-th order fractional Brownian motion.
Finally, we state an asymptotic analysis of the limiting error process $\{R_t(\kappa) : t\in[0,T]\}$, where $R_t(\kappa)\coloneqq\lim_{m\to+\infty}R_t(m,\kappa)$.
The following statements hold \cite[Theorem 6.4]{kawai2016higher}.
\begin{enumerate}[(i)]
	\setlength{\parskip}{0cm}
	\setlength{\itemsep}{0cm}
	\item The finite dimensional distribution of $\{R_t(\kappa) : t\in[0,T]\}$ converges in probability to zero as $\kappa\to+\infty$. Moreover, if $\alpha\in(1,2)$, then the convergence can be strengthened to convergence in probability to the degenerate zero process uniformly on $[0,T]$.
		
	\item For every $\kappa>0$, as $\lambda\to+\infty$, it holds that
	\[
		\lambda^\alpha\mathbb P\left(\sup_{t\in[0,T]} |R_t(\kappa)|>\lambda\right) \to c_\alpha\int_{-\infty}^{-\kappa}|f_n(T,s;H,\alpha)|^\alpha\,ds.
	\]
	Moreover, the supremum over $[0,T]$ can be replaced by a maximum over a finite number of observation times in $[0,T]$.
\end{enumerate}

Thus, we see the error process can be made arbitrarily close to the degenerate zero process by increasing $\kappa$, and that the tail of $\sup_{t\in[0,T]}|R_t(\kappa)|$ resembles that of a Pareto distribution. We have therefore successfully approximated the higher order fractional stable motion.

This method of approximation for infinitely divisible processes is rather systematic. To summarise, we do so by firstly considering a domain of the L\'evy measure for which it is finite, which is analogous to an approximation via Poisson truncation of shot noise representation. If possible, Gaussian approximation can be additionally used on remaining components to obtain greater accuracy. The remaining stochastic process $\{R_t(m,\kappa):t\in[0,T]\}$ is treated as the error of the approximation, for which we analyse its properties.

We briefly mention the tempered stable generalisation of the higher order fractional stable motion. By replacing the L\'evy measure of the driving process in the stochastic integral \eqref{nfsm} with that of a symmetric tempered stable L\'evy measure
\begin{equation}\label{symmetric tempered stable Levy measure}
	\nu(dz) = \frac{\alpha c_\alpha}{2}\frac{e^{-\beta|z|}}{|z|^{\alpha+1}}dz,\quad z\in\mathbb R,
\end{equation}
where $\beta>0$ is a tempering parameter, we obtain the higher order fractional tempered stable motion. Its properties, such as persistent autocorrelations and behaviour in short and long time regimes, are investigated in \cite{carnaffan2019analytic}. A shot noise representation in the vein of \eqref{nfsm truncation} is given by \cite[Section 6]{carnaffan2019analytic}
\begin{equation}\label{nftsm truncation}
	\left\{\sum_{\{k\in\mathbb N: \Gamma_k\le m(T+\kappa)\}}\left(\left(\frac{\Gamma_k}{(T+\kappa)c_\alpha}\right)^{-1/\alpha} \wedge \left(\frac{E_k R_k^{1/\alpha}}{\beta}\right)\right)U_kf_n(t,T_k;H,\alpha) : t\in[0,T]\right\},
\end{equation}
where $\{E_k\}_{k\in\mathbb N}$ is a sequence of iid standard exponential random variables and $\{R_k\}_{k\in\mathbb N}$ is a sequence of iid standard uniform random variables, such that all random sequences are mutually independent. This shot noise truncation is not surprising at all, given Rosi\'nski's series representation \eqref{eq tempered stable law series representation} for the tempered stable law.

An alternative fractional tempered stable motion is studied in \cite{houdre2006fractional}, where the L\'evy measure of the driving process in the stochastic integral \eqref{nfsm} is once again replaced by a tempered stable L\'evy measure \eqref{symmetric tempered stable Levy measure}, but the first-order moving average kernel $f_1$ is also replaced, by a Volterra kernel
\begin{equation}\label{Volterra kernel}
	K(t,s;H,\alpha) \coloneqq c_{H,\alpha}\left(\left(\frac ts\right)^{H-1/\alpha}(t-s)^{H-1/\alpha} - \left(H - \frac1\alpha\right)s^{1/\alpha - H}\int_s^t u^{H-1/\alpha - 1}(u-s)^{H-1/\alpha}\,du\right)\mathbbm1_{[0,t]}(s),
\end{equation}
where $H\in(1/\alpha - 1/2, 1/\alpha + 1/2)$, $\alpha\in(0,2)$ and $c_{H,\alpha}$ is a constant depending only on $H$ and $\alpha$. Note that as the kernel in this fractional tempered stable motion is distinct from the higher order fractional tempered stable motion discussed in \cite[Section 6]{carnaffan2019analytic}, the former is not simply a lower-order version of the latter. The fractional tempered stable motion based on the Volterra kernel \eqref{Volterra kernel} does not integrate over negative time, so within the context of the decomposition \ref{stochastic integral decomposition}, the truncation parameter on time is degenerate ($\kappa=0$). Consequently, the error analysis is simplified as the stochastic process $\{Q_t(m):t\in[0,T]\}$ of the small jumps also has no negative time component, and the component $\{R_t(m,\kappa):t\in[0,T]\}$ is irrelevant. Similarly to the moving average kernel, the Volterra kernel can capture selfsimilar dynamics. However, unlike the moving average kernel, there is no obvious generalisation of the Volterra kernel that can preserve its key properties.

\subsubsection{L\'evy-driven CARMA processes}\label{section CARMA process}
We now turn to a numerical scheme for generating approximate sample paths of L\'evy-driven CARMA processes via Poisson truncation of shot noise representation \cite{kawai2017sample}. L\'evy-driven CARMA processes naturally generalise Gaussian CARMA processes so as to capture asymmetry and heavy tails in a variety of physical and social science settings.
We begin defining the L\'evy-driven CARMA process as follows.
Fix $a_1,\cdots,a_p,b_0,\cdots,b_{p-1}\in\mathbb R$ such that $b_q = 1$, $q\le p-1$ and $b_k = 0$ for $k>q$, and define the polynomials $a(z) \coloneqq z^p + a_1z^{p-1} + \cdots + a_p$ and $b(z)\coloneqq b_0 + b_1 z + \cdots + b_qz^q$ in such a way that $a(z)$ and $b(z)$ have no common roots. Define $A\in\mathbb R^{p\times p}$ by
\[
	A \coloneqq \begin{bmatrix}
		0 & 1 & 0 & \cdots & 0\\
		0 & 0 & 1 & \cdots & 0\\
		\vdots & \vdots & \vdots & \cdots & \vdots\\
		0 & 0 & 0 & \cdots & 1\\
		-a_p & -a_{p-1} & -a_{p-2} & \cdots & -a_1
	\end{bmatrix}.
\]
Denote the eigenvalues of $A$ as $\lambda_1,\cdots,\lambda_p\in\mathbb R$, that is, $a(z) = \prod_{k=1}^p(z-\lambda_k)$. Let $\{L_t:t\in\mathbb R\}$ be a temporally extended univariate L\'evy process in same sense as in Section \ref{section nfsm}.
Define $e_p\in\mathbb R^p$ as the unit vector in the $p$-th direction and $\mathbf b\coloneqq[b_1, b_1,\cdots,b_{p-1}]^\intercal\in\mathbb R^p$.
A \textit{L\'evy-driven CARMA process in $\mathbb R$ of order $(p,q)$ with $p>q$} is defined as $\{Y_t:t\in\mathbb R\}$, where $Y_t\coloneqq\langle\mathbf b,X_t\rangle$ and $\{X_t:t\in\mathbb R\}$ is a stochastic process in $\mathbb R^p$ satisfying
\[
	X_{t_2} = e^{A(t_2-t_1)}X_{t_1} + \int_{t_1}^{t_2} e^{A(t_2 - s)}e_p\,dL_s,\quad t_1\le t_2.
\]
Under suitable technical conditions, the L\'evy-driven CARMA process $\{Y_t:t\in\mathbb R\}$ is strictly stationary, and can be expressed as a linear combination of the real and dependent L\'evy-driven Ornstein--Uhlenbeck processes as follows:
\[
	Y_t = \sum_{k=1}^p\frac{b(\lambda_k)}{a'(\lambda_k)}\int_{-\infty}^t e^{\lambda_k(t-s)}\,dL_s,\quad t\in\mathbb R.
\]
Note that the CARMA$(1,0)$ process corresponds to the L\'evy-driven Ornstein-Uhlenbeck process, which we will discuss exclusively in Section \ref{section OU process} due to their special feature with respect to sample path generation.

We first consider the stable CARMA process. Specifically, the driver $\{L_t:t\in\mathbb R\}$ is a temporally extended stable process with L\'evy measure
\[
	\nu(dz) = \alpha c_\alpha\left(\frac{1+\beta}{2}\mathbbm1_{(0,+\infty)}(z) + \frac{1-\beta}{2}\mathbbm1_{(-\infty,0)}(z)\right)\frac{dz}{|z|^{\alpha+1}},\quad z\in\mathbb R_0,
\]
where $\beta\in[-1,1]$ if $\alpha\in(0,1)\cup(1,2)$ and zero if $\alpha=1$. We define a constant $\lambda_{\alpha,\beta,m} \coloneqq \beta c_\alpha\eta(m)^{1-\alpha}\alpha/(\alpha-1)$, which is used in a correction term to centre the error term $Q_t(m)$ when $\alpha\in(0,1)$. Approximation and error analysis for this case is provided as follows.
Define the function $g(u) \coloneqq \sum_{k=1}^p \mathbbm1(u\ge0)e^{\lambda_k u}b(\lambda_k)/a'(\lambda_k)$.
In the context of the decomposition \eqref{stochastic integral decomposition} for the stable CARMA process based on the inverse L\'evy measure method and the truncation on the jump timings $\mathcal T_n\coloneqq (n,T]$, the following statements hold \cite[Section 4]{kawai2017sample}:
\begin{enumerate}[(i)]
	\setlength{\parskip}{0cm}
	\setlength{\itemsep}{0cm}
	\item It holds that
	\begin{equation}\label{stable CARMA truncation}
		\{Y_t(m,n) : t\in[0,T]\} \stackrel{\mathscr L}{=}\left\{\sum_{\{k\in\mathbb N: \Gamma_k \le m(T+n)\}}\left(\frac{\Gamma_k}{(T+n)c_\alpha}\right)^{-1/\alpha}g(t-U_k)r_k - \mathbbm1_{(1,2)}(\alpha)\lambda_{\alpha,\beta,m}\int_{-n}^T g(t-s)\,ds : t\in[0,T]\right\},
	\end{equation}
	where $\{\Gamma_k\}_{k\in\mathbb N}$ is a sequence of standard Poisson arrival times, $\{U_k\}_{k\in\mathbb N}$ is a sequence of iid uniform random variables on $(-n,T)$, $\{r_k\}_{k\in\mathbb N}$ a sequence of iid random variables such that $\mathbb P(V_1 = 1) = (1+\beta)/2$ and $\mathbb P(V_1 = -1) = (1-\beta)/2$, and that all random sequences are mutually independent.
		
	\item As $m\to+\infty$, it holds that
	\[
		\left\{\frac{1}{\sigma(m)}\left(Q_t(m) + \mathbbm1_{(0,1)}(\alpha)\lambda_{\alpha,\beta,m}\int_{-\infty}^T g(t,s)\,ds\right) : t\in[0,T]\right\} \stackrel{\mathscr L}{\to} \left\{\int_{-\infty}^T g(t-s)\,dB_s:t\in[0,T]\right\}.
	\]
	If $q<p-1$, then the convergence can be strengthened to weak convergence in $\mathscr C([0,T];\mathbb R)$.
		
	\item The finite dimensional distributions of $\{R_t(n) : t\in[0,T]\}$ converges in probability to zero as $n\to+\infty$. If $\alpha\in(1,2)$, then the convergence can be strengthened to the convergence in probability uniformly on $[0,T]$.
		
	\item It holds that for every $n>0$ as $\lambda\to+\infty$,
	\[
		\lambda^\alpha \mathbb P\left(\sup_{t\in[0,T]}|R_t(n)|>\lambda\right) \to c_\alpha\int_{-\infty}^{-n}\sup_{t\in[0,T]}|g(t-s)|^\alpha\,ds.
	\]
\end{enumerate}

\begin{figure}[h]\centering
	\begin{tabular}{cc}
		\includegraphics[width=90mm]{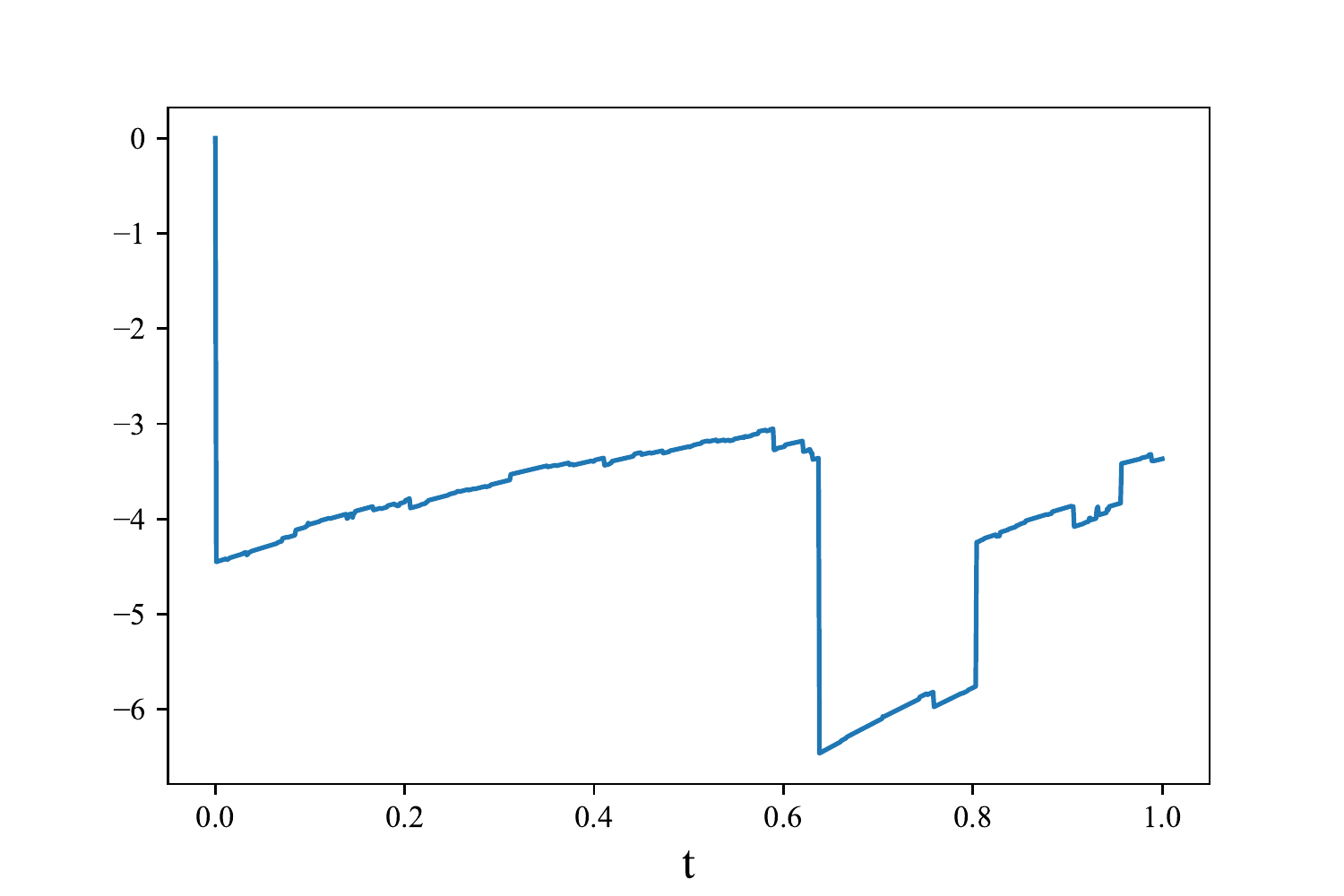} & \includegraphics[width=90mm]{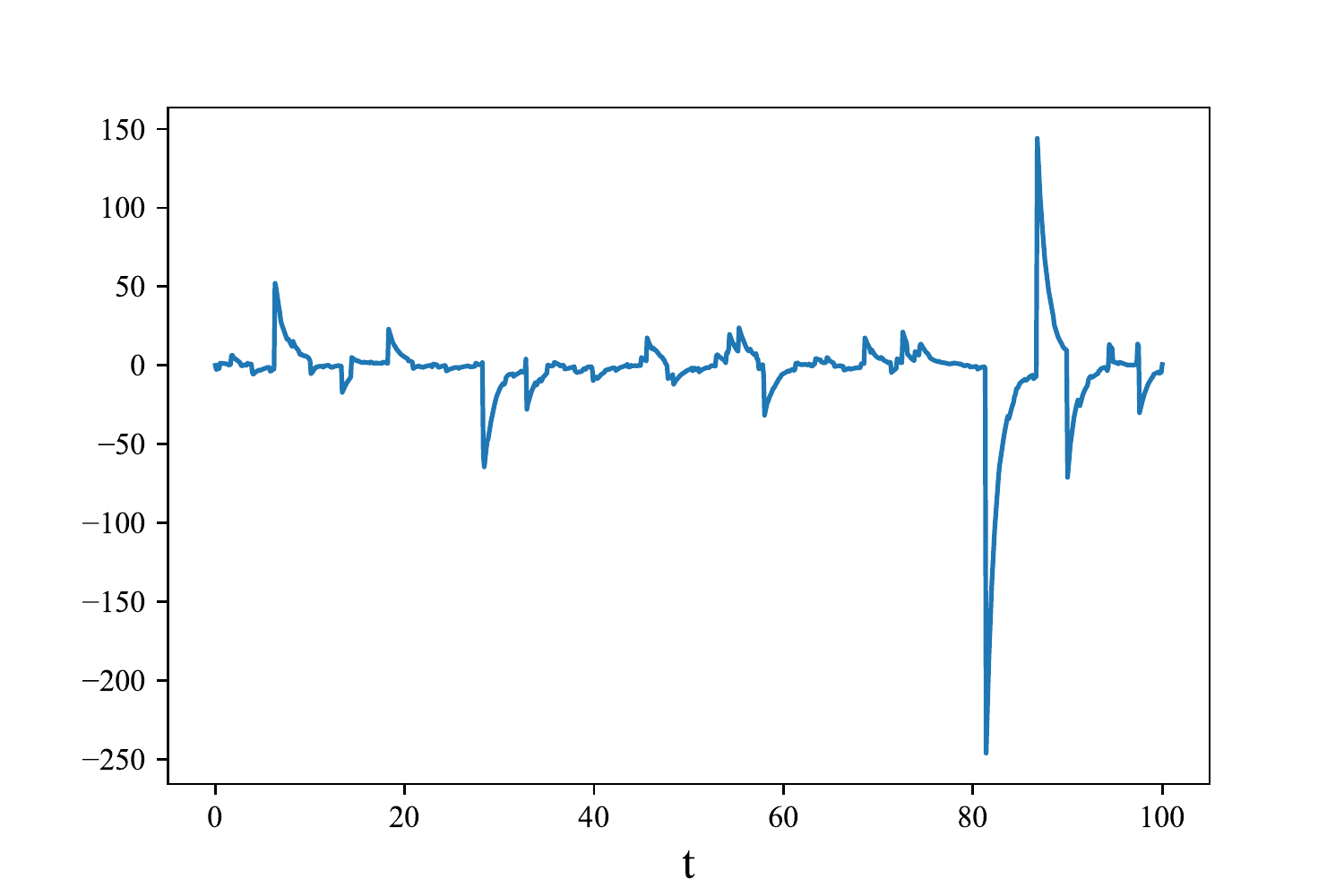} \\
		(a) $T=1$ & (b) $T=100$
	\end{tabular}
	\caption{Examples of sample paths of the stable CARMA$(2,1)$ process with $\lambda_1 = -1.2$, $\lambda_2 = -0.5$, $\alpha=0.8$ and $\beta=0$, based on the truncation approximation \eqref{stable CARMA truncation}. Time steps are of size $10^{-3}$.}
	\label{fig stable CARMA}
\end{figure}

So clearly, the error component $\{Q_t(m) : t\in[0,T]\}$ is asymptotically Gaussian while the error component $\{R_t(n) : t\in[0,T]\}$ asymptotically has a Pareto-tailed distribution, just as for the higher order fractional stable motion in Section \ref{section nfsm}. The latter is unsurprising due to the presence of the stable driver. Following the simulation recipe outlined in Section \ref{section Levy-driven stochastic integrals} but based on the shot noise representation \eqref{stable CARMA truncation} instead, we provide approximate sample paths in Figure \ref{fig stable CARMA}. From the generated sample paths, we observe the mean-reverting property of the CARMA process in diminishing the impact of jumps over time.

Next, we consider the case where the driving L\'evy process $\{L_t:t\in\mathbb R\}$ without Gaussian components is centred and has finite second-order moments. Specifically, we demand that the characteristic function follows the form
\[
	\mathbb E\left[e^{i\theta L_t}\right] = \exp\left[t \int_{\mathbb R_0}\left(e^{-i\theta z} - 1 - i\theta z\right)\,\nu(dz)\right],\quad t\ge0,
\]
such that the infinite L\'evy measure satisfies $\int_{\mathbb R_0}z^2\,\nu(dz) < +\infty$.
We now deviate from the inverse L\'evy measure method and leave the underlying shot noise representation general. We state the results on the error processes from the decomposition \eqref{stochastic integral decomposition} as follows \cite[Proposition 5.3]{kawai2017sample}:
\begin{enumerate}[(i)]
	\setlength{\parskip}{0cm}
	\setlength{\itemsep}{0cm}
	\item The error process $\{Q_t(m):t\in[0,T]\}$ converges in probability to the zero process uniformly on $[0,T]$ as $m\to+\infty$. Moreover, where $\sigma^2(m) \coloneqq\int_{\mathbb R_0}z^2\,\nu_m(dz)$, if for every $c>0$
		\[
			\lim_{m\to+\infty}\frac{1}{\sigma^2(m)}\int_{\{|z|>c\sigma(m)\}}z^2\,(\nu-\nu_m)(dz) = 0,
		\]
		then it holds that
		\[
			\left\{\frac{Q_t(m)}{\sigma(m)} : t\in[0,T]\right\} \stackrel{\mathscr L}{=} \left\{\int_{-\infty}^T g(t-s)\,dB_s : t\in[0,T]\right\}.
		\]
		Additionally, if $q<p-1$, then the convergence can be strengthened to weak convergence in the space $\mathscr C([0,T];\mathbb R)$.
		
	\item It holds that for every $m>0$, the error process $\{R_t(m,n) : t\in[0,T]\}$ converges in probability to the zero process uniformly on $[0,T]$ as $n\to+\infty$.
\end{enumerate}

By (ii), we see that the error component $\{Q_t(m) : t\in[0,T]\}$ can be approximated by a Gaussian process. Moreover, it has been shown that when possible, including the Gaussian approximation in the simulation for $\{Y_t:t\in[0,T]\}$ not only improves precision, but is also necessary to preserve the second-order structure \cite[Proposition 5.3]{kawai2017sample}.

In conclusion, we see that the general idea of decomposing a L\'evy measure in terms of its jump magnitudes and timings is a powerful tool that can be applied systematically for approximating sample paths of infinitely divisible processes.

\subsubsection{L\'evy-driven Ornstein--Uhlenbeck process}\label{section OU process}
We now consider the L\'evy-driven Ornstein--Uhlenbeck (OU) process, which can be thought of as a special case of the L\'evy-driven CARMA process with $(p,q)=(1,0)$ but with an additional centring parameter. This class of infinitely divisible processes deserves special attention in the context of numerical methods, as exact simulation methods by increments are available in some cases, which may be preferred over the Poisson truncation of shot noise representation.

Let $\{L_t:t\ge0\}$ be a L\'evy process without Gaussian components and with L\'evy measure $\nu(d\mathbf z)$. The L\'evy-driven OU process $\{X_t:t\ge0\}$ is described by the stochastic differential equation
\begin{equation}\label{OU SDE}
	dX_t = \lambda(\mu-X_t)\,dt + dL_{t},
\end{equation}
where $\lambda>0$ and $\mu\in\mathbb R^d$. The explicit solution is given by
\begin{equation}\label{OU explicit solution}
	X_t = e^{-\lambda t}X_0 + \mu\left(1-e^{-\lambda t}\right) + \int_0^t e^{-\lambda(t-s)}\,dL_{s}.
\end{equation}
In light of this stochastic integral representation, the L\'evy-driven OU process is an infinitely divisible process.

The OU process is typically defined by its invariant law $\lim_{t\to+\infty}\mathscr L(X_t)$. In the case where the invariant law is Gaussian, the driver of the OU process \eqref{OU SDE} is a Brownian motion, which is outside the scope of the present survey.
It is known that \cite[Theorem 17.5]{sato1999levy} if the L\'evy measure $\rho(d\mathbf z)$ of the driving process $\{L_t:t\ge0\}$ satisfies the integrability condition $\int_{\|\mathbf z\|>2} \ln\|\mathbf z\|\,\rho(d\mathbf z) < +\infty$, then there exists a L\'evy-driven OU process \eqref{OU explicit solution} such that the invariance law $\lim_{t\to+\infty}\mathscr L(X_t)$ exists and is selfdecomposable and infinitely divisible with the L\'evy measure $\nu(B) = \lambda^{-1} \int_{\mathbb R^d}\int_0^{+\infty}\mathbbm1_B(e^{-s}\mathbf z)\,ds\,\nu(d\mathbf z)$, for $B\in\mathcal B(\mathbb R^d_0)$.
Conversely, every selfdecomposable law admits the unique existence of a L\'evy-driven OU process such that the L\'evy measure of its L\'evy driver satisfies the aforementioned integrability condition \cite{sato1999levy}. For simplicity, we focus on the univariate setting. Where $w(z)$ is the L\'evy density of the unit-time marginal $L_1$ of the driving L\'evy process, it is related to $u(z)$ by the equation
\begin{equation}\label{driver decomposition}
	w(z) = -\lambda\left(u(z) + z\frac{\partial}{\partial z}u(z)\right).
\end{equation}

A (non-Gaussian) stable OU process is defined as the L\'evy-driven OU process such that its invariant law $\lim_{t\to+\infty}\mathscr L(X_t)$ is a (non-Gaussian) stable law.
If the invariant law is a stable law with the one-sided L\'evy density $u(z) = az^{-\alpha-1}$ on $(0,+\infty)$ with $\alpha\in(0,1)$, then the driving process admits the L\'evy density $w(z) = \lambda a\alpha z^{-\alpha-1}$ on $(0,+\infty)$, according to the relation \eqref{driver decomposition}, that is, the driving process is a stable subordinator with stability $\alpha$, but with scale $\lambda a\alpha$. Combining results from Example \ref{example stable series representation} and Section \ref{section Levy-driven stochastic integrals}, we have the shot noise representation
\begin{equation}\label{SOU series representation}
	\{X_t:t\in[0,T]\} \stackrel{\mathscr L}{=} \left\{e^{-\lambda t}X_0 + \mu\left(1-e^{-\lambda t}\right) + \sum_{k=1}^{+\infty}e^{-\lambda(t-T_k)}\left(\frac{\Gamma_k}{\lambda aT}\right)^{-1/\alpha}\mathbbm1_{[0,t]}(T_k) : t\in[0,T]\right\},
\end{equation}
where $\{\Gamma_k\}_{k\in\mathbb N}$ is a sequence of standard Poisson arrival times independent of $\{T_k\}_{k\in\mathbb N}$, a sequence of iid uniform random variables on $(0,T)$. Sample paths based on the Poisson truncation of the shot noise representation \eqref{SOU series representation} are provided in Figure \ref{fig SOU} below. Similar to our numerical illustrations of the stable CARMA process in Figure \ref{fig stable CARMA}, we also observe mean-reverting behaviour in the sample paths of the stable OU process. In contrast, the driving L\'evy process in this case is a stable subordinator, thus we see that all jumps in the sample paths of Figure \ref{fig SOU} are in the positive direction. Similarly to the case of the stable process in Figure \ref{fig 2D stable paths}, the jumps associated with a lower stability parameter (Figure \ref{fig SOU} (a)) are significantly larger.

\begin{figure}[h]\centering
	\begin{tabular}{cc}
		\includegraphics[width=90mm]{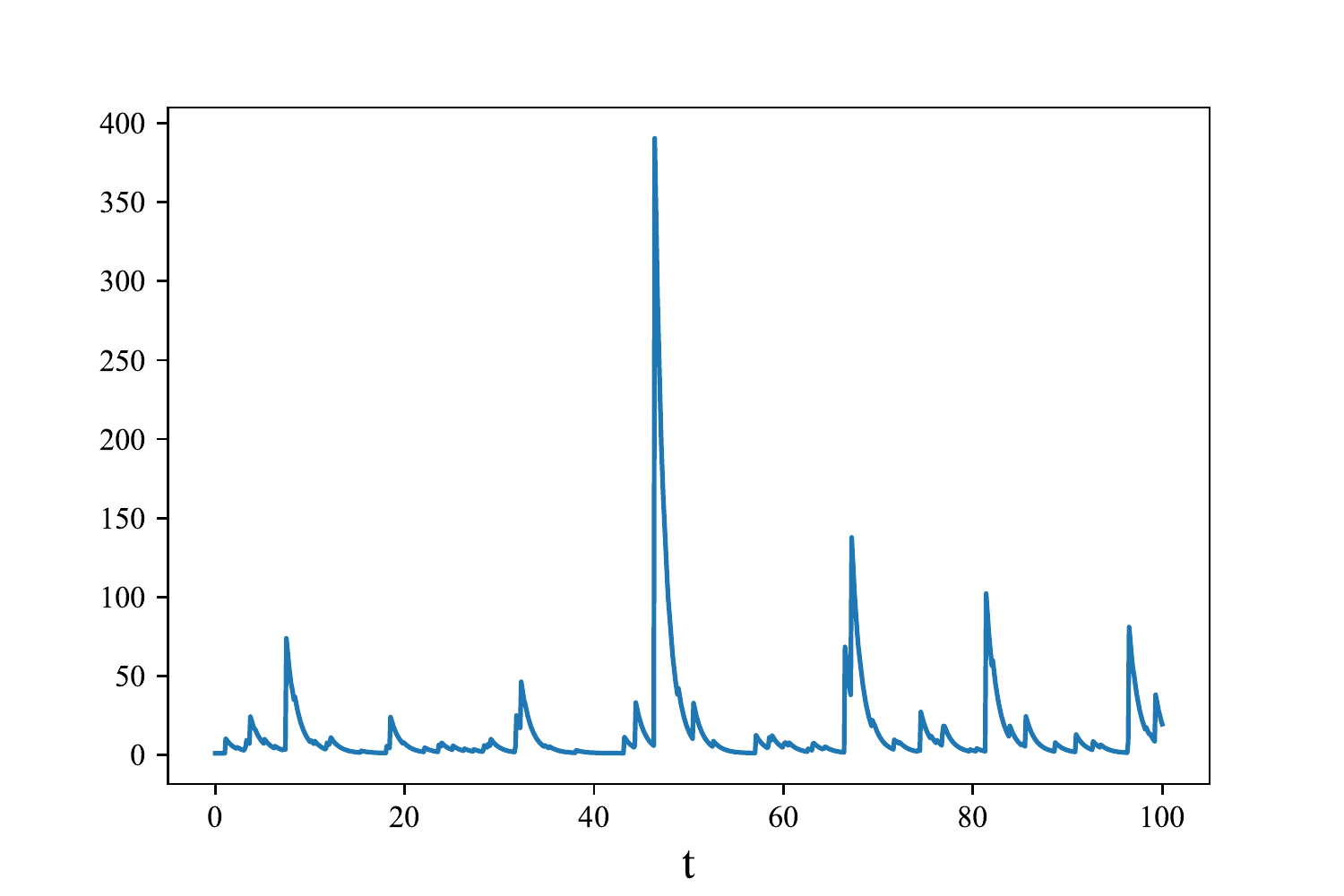} & \includegraphics[width=90mm]{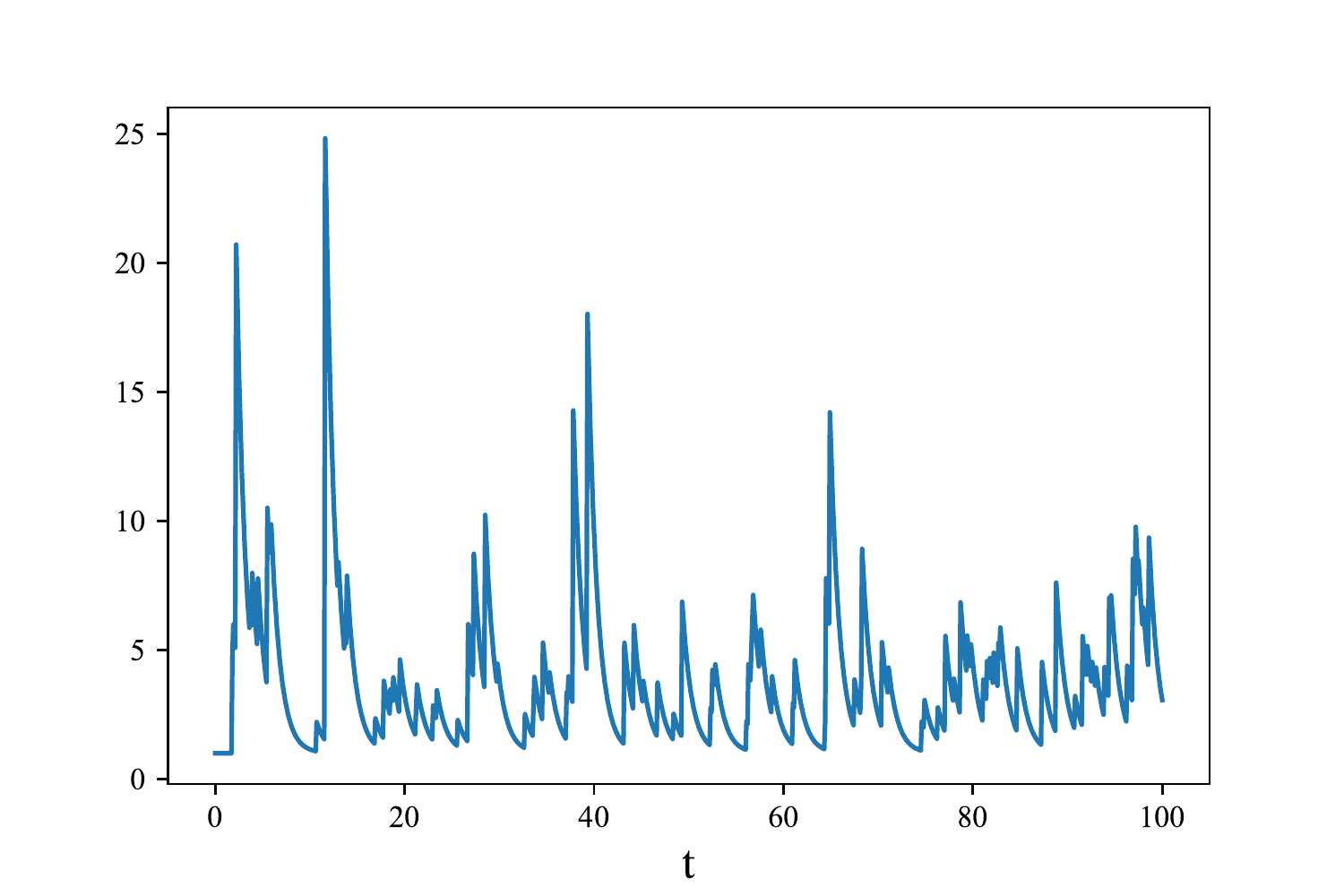} \\
		(a) $\alpha=0.7$ & (b) $\alpha=1.5$
	\end{tabular}
	\caption{Examples of approximate sample paths of the stable OU process \eqref{OU explicit solution} based on truncation of the shot noise representation \eqref{SOU series representation}. Other parameters are fixed at $a=1$, $\mu=1$, $\lambda=1$, $X_0=1$ and $T=100$. The truncation parameter is fixed at $n=100$ and time steps are of size $10^{-3}$.}
	\label{fig SOU}
\end{figure}

An interesting property of the stable OU process is that the L\'evy measures of the invariant law and the driving process both correspond to a stable law with stability $\alpha$, though with different scales. This invariance of the stability parameter between the L\'evy measures does not hold in general, as we see in the following.

\begin{figure}[h]\centering
	\begin{tabular}{cc}
		\includegraphics[width=90mm]{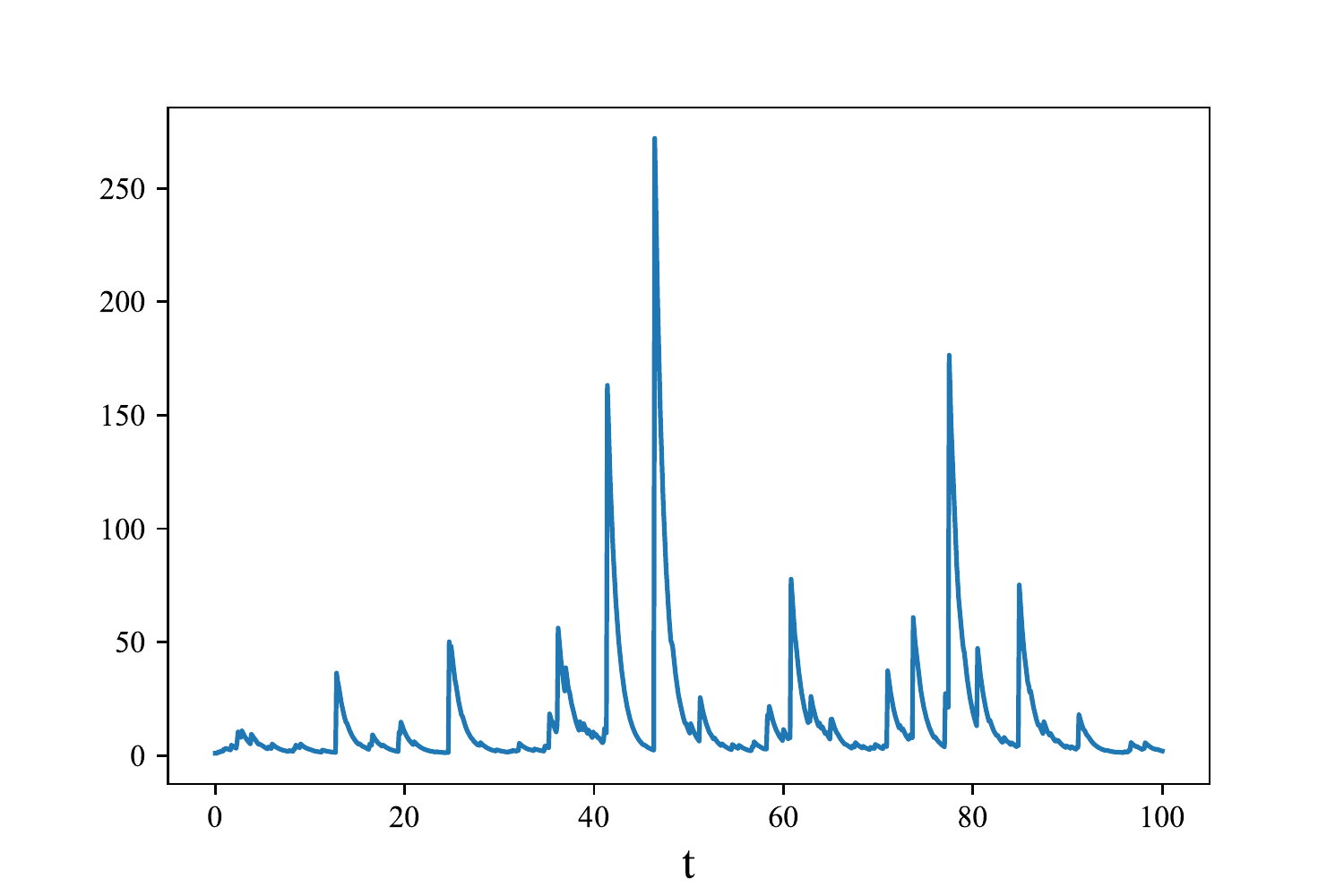} & \includegraphics[width=90mm]{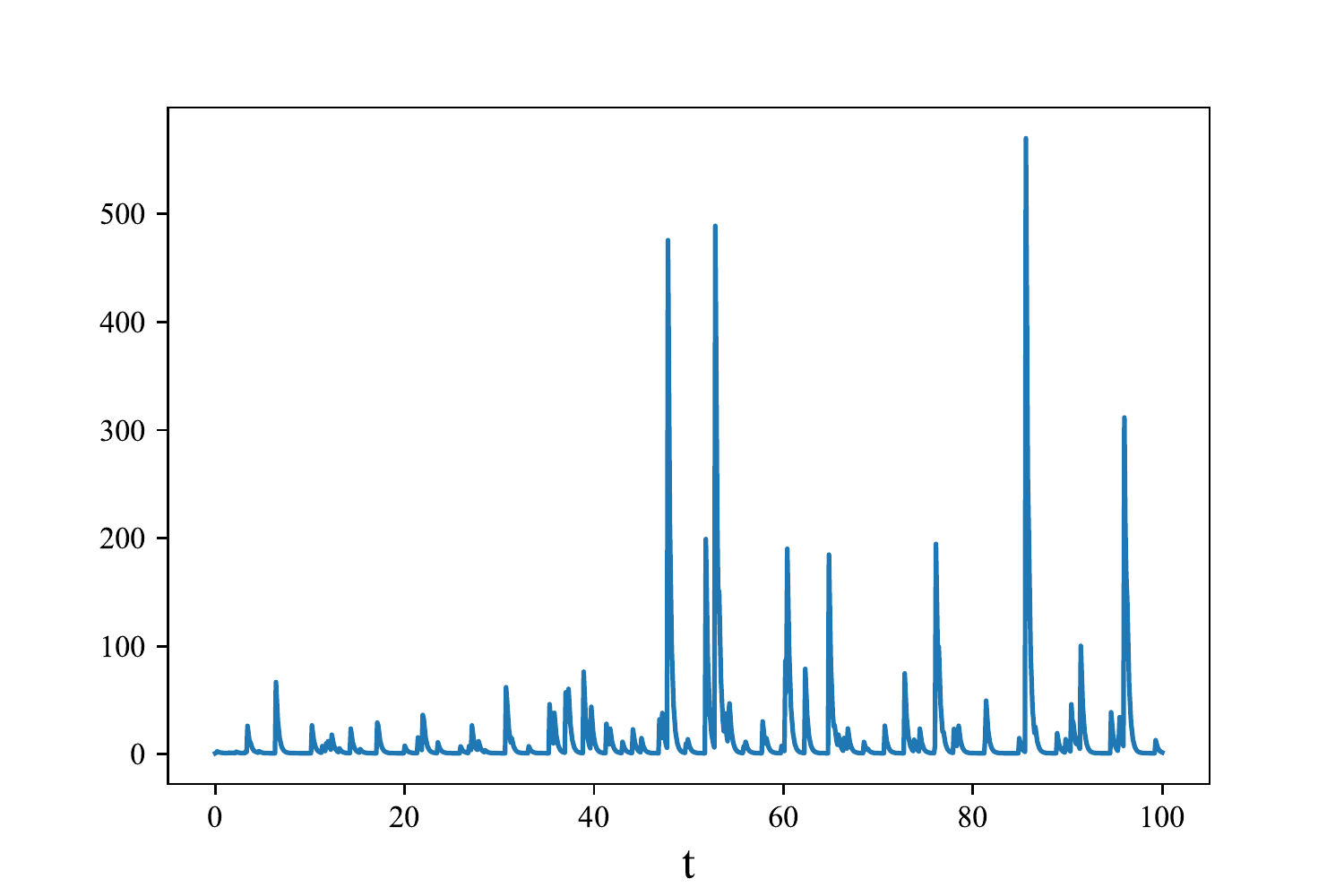} \\
		(a) $\lambda=1$ & (b) $\lambda=4$
	\end{tabular}
	\caption{Examples of approximate sample paths of the tempered stable OU process \eqref{OU explicit solution} based on truncation of the shot noise representation \eqref{TSOU series representation}. Other parameters are fixed at $\alpha=0.7$, $a=1$, $\beta=1$, $\mu=1$, $X_0=1$ and $T=100$. The truncation parameter is fixed at $n=100$ and time steps are of size $10^{-3}$.}
	\label{fig TSOU}
\end{figure}

Consider a tempered stable OU process such that its invariant law $\lim_{t\to+\infty}\mathscr L(X_t)$ is a tempered stable law with L\'evy density $u(z) = ae^{-\beta z}z^{-\alpha-1}$ on $(0,+\infty)$ with $\alpha\in(0,1)$. From the relation \eqref{driver decomposition}, the L\'evy density of the unit-time marginal $L_1$ of the driving process is given by
\begin{equation}
	w(z) = \lambda a\alpha\frac{e^{-\beta z}}{z^{\alpha+1}} + \lambda a \beta\frac{e^{-\beta z}}{z^{(\alpha-1)+1}},\quad z\ge0.
\end{equation}
Thus, for $\alpha\in(0,1)$, the L\'evy driver $\{L_t:t\ge0\}$ is the superposition of a tempered stable subordinator with stability $\alpha$ and scale $\lambda a\alpha$, and a compound Poisson process with L\'evy density $\lambda a\beta z^{-\alpha}e^{-\beta z}$, which is indeed a gamma density scaled by a factor of $\lambda a$.
This leads to the shot noise representation \cite{kawai2012infinite,rosinski2007tempering}
\begin{multline}\label{TSOU series representation}
	\{X_t:t\in[0,T]\} \stackrel{\mathscr L}{=} \Bigg\lbrace e^{-\lambda t}X_0 + \mu\left(1-e^{-\lambda t}\right) + \sum_{k=1}^{+\infty}e^{-\lambda(t-T_k)}\left(\left(\frac{\Gamma_k}{\lambda aT}\right)^{-1/\alpha}\wedge \left(\frac{W_k U_k^{1/\alpha}}{\beta}\right)\right)\mathbbm1_{[0,t]}(T_k) \\+ \sum_{k=1}^{+\infty}e^{\widetilde \Gamma_k - \lambda t}G_k\mathbbm1_{[0,t]}(\widetilde \Gamma_k) : t\in[0,T]\Bigg\rbrace,
\end{multline}
where, in addition to the random sequences appearing in Rosi\'nski's series representation \eqref{eq tempered stable law series representation}, $\{\widetilde \Gamma_k\}_{k\in\mathbb N}$ is a sequence of Poisson arrival times with intensity $\lambda a\beta^\alpha\Gamma(1-\alpha)$ independent to $\{\Gamma_k\}_{k\in\mathbb N}$, and $\{G_k\}_{k\in\mathbb N}$ is a sequence of iid gamma random variables with shape $1-\alpha$ and rate $\beta$. Examples of sample paths based on the Poisson truncation of the shot noise representation \eqref{TSOU series representation} are provided in Figure \ref{fig TSOU} as follows. As expected, the jumps of the tempered stable OU process tend to be smaller than those of its stable counterpart in Figure \ref{fig SOU} (a). We see from the plots in Figure \ref{fig TSOU} that the parameter $\lambda$ controls the intensity of the mean-reverting property.

For $\alpha\in(1,2)$, by the relation \eqref{driver decomposition}, the driving process $\{L_t:t\ge0\}$ is instead a superposition of two independent tempered stable processes, one with stability $\alpha$ and scale $\lambda a\alpha$, and another with stability $\alpha-1$ and scale $\lambda a\beta$. A shot noise representation similar to \eqref{TSOU series representation} is available, however intricate centring terms must appear. For the sake of simplicity, consider the case where the invariant law of the OU process $\lim_{t\to+\infty}\mathscr L(X_t)$ is a symmetric tempered stable law, so that the independent tempered stable processes forming the driving L\'evy process are also symmetric. As almost sure convergence without centres is guaranteed, the shot noise series for the symmetric tempered stable OU process with stability $\alpha\in(1,2)$ is given by
\begin{multline}\label{TSOU series representation 2}
	\{X_t:t\in[0,T]\} \stackrel{\mathscr L}{=} \Bigg\lbrace e^{-\lambda t}X_0 + \mu\left(1-e^{-\lambda t}\right) + \sum_{k=1}^{+\infty}e^{-\lambda(t-T_k)}\left(\left(\frac{\Gamma_k}{\lambda aT}\right)^{-1/\alpha}\wedge \left(\frac{W_k U_k^{1/\alpha}}{\beta}\right)\right)V_k\mathbbm1_{[0,t]}(T_k) \\+ \sum_{k=1}^{+\infty}e^{-\lambda\left(t-\widetilde T_k\right)}\left(\left(\frac{(\alpha-1)\widetilde\Gamma_k}{\lambda a\beta T}\right)^{-1/(\alpha-1)}\wedge \left(\frac{\widetilde W_k \widetilde U_k^{1/(\alpha-1)}}{\beta}\right)\right)\widetilde V_k\mathbbm1_{[0,t]}\left(\widetilde T_k\right) : t\in[0,T]\Bigg\rbrace,
\end{multline}
where the familiar random sequences are the same as those in the representation \eqref{TSOU series representation}, and the random sequences distinguished by the tilde are independent copies of the corresponding sequences.

With the availability of shot noise series representations for the stable \eqref{SOU series representation} and tempered stable OU processes \eqref{TSOU series representation}, we immediately have an easy method for sample path generation by Poisson truncation. As usual, this entails truncation errors due to the discarded jumps, and a decomposition of the stochastic integral process in the same vein as \eqref{stochastic integral decomposition} can be performed to analyse the error. The kernel $s\mapsto e^{-\lambda(t-s)}\mathbbm1_{[0,t]}(s)$ is square-integrable and uniformly bounded for every $t$, which conveniently permits Gaussian approximation of the discarded jumps.
Alternatively, as the domain of integration over time is bounded, it may be easier to investigate the mean-squared error. Consider an infinitely divisible process $\{\int_0^t f(t,s)\,dL_s:t\in[0,T]\}$, where $\{L_t:t\in[0,T]\}$ is a L\'evy process with an isotropic L\'evy measure $\nu(dz)$ such that it has the shot noise representation \eqref{generalised shot noise series for Levy process} without the centring terms. Then, a shot noise representation for the infinitely divisible process can be obtained by simply including $f(t,T_k)$ as a factor in each of the summands.
For the Poisson truncation approximation, the It\^o-Wiener isometry can be applied in the same vein as \eqref{mse} to express the mean-squared error as $\mathbb E[|Q_t(m)|^2] = \int_0^t |f(t,s)|^2\,ds\int_{\mathbb R_0}|z|^2\,(\nu-\nu_n)(dz)$, where $Q_t(m)$ previously defined in \eqref{stochastic integral decomposition} integrates over $[0,T]$ instead of $(-\infty,T]$.
This result can be straightforwardly generalised to the case of a multivariate L\'evy integrator similarly to \eqref{mse}. We see that as long as the integrand $f(t,s)$ is square-integrable, the quality of the Poisson truncation is characterised by the L\'evy integrator.

Yet another representation is avaliable via the Markovian property, for which exact sample path generation of the L\'evy-driven OU process is possible in some cases. Specifically, where $\Delta$ denotes the size of the time step, it holds that
\begin{equation}\label{OU Markovian formula}
	X_{(n+1)\Delta} = e^{-\lambda\Delta}X_{n\Delta} + \mu\left(1-e^{-\lambda\Delta}\right) + \int_{n\Delta}^{(n+1)\Delta} e^{-\lambda((n+1)\Delta - s)}dL_s \stackrel{\mathscr L}{=} e^{-\lambda\Delta}X_{n\Delta} + \mu\left(1-e^{-\lambda\Delta}\right) + \int_{0}^{\Delta} e^{-\lambda(\Delta - s)}dL_s,
\end{equation}
where the equality in law holds by the independence and stationarity of the L\'evy driver. In the one-dimensional setting, exact time-discretisation schemes based on \eqref{OU Markovian formula} for the stable OU process with stability $\alpha\in(0,2)$ and the tempered stable OU process with stability $\alpha\in(0,1)$ are possible, as a stable random variable $S_{\alpha,a}$ in $\mathbb R$ with stability $\alpha$ and scale $a$ can be exactly simulated through the well-known representation
\begin{equation}\label{stable representation}
	S_{\alpha,a} \stackrel{\mathscr L}{=} \left(\frac{a\Gamma(1-\alpha)}{\alpha\cos(V)}\right)^{1/\alpha}\sin(\alpha(V+\pi/2))\left(\frac{\cos(V-\alpha(V+\pi/2))}{E}\right)^{(1-\alpha)/\alpha},
\end{equation}
where $V$ is a uniform random variable on $(-\pi/2,\pi/2)$ independent of $E$, a standard uniform random variable. For $\alpha\in(0,1)$, the transition law of the tempered stable OU process consists of a compound Poisson random variable and a tempered stable random variable with stability $\alpha$. The latter of which can be generated exactly by an acceptance-rejection algorithm on $S_{\alpha,a}$ \eqref{stable representation}, so an exact time-discretisation scheme for sample path generation is readily available \cite{kawai2011exact}.

Therefore, we see that there is a reasonable basis for turning to the classical recursive sampling by increments as opposed to by Poisson truncation of shot noise representation. However, there are two bases for which one may prefer using Poisson truncation of shot noise representation. Firstly, if the application context requires the observation of jumps, such as in insurance mathematics, then simulation via increments cannot be applied. Secondly, the Markovian successive representation \eqref{OU Markovian formula} may not offer an exact sampling scheme for all OU processes. For instance, there is yet an exact sampling method for tempered stable random variables with stability $\alpha\in(1,2)$, which appears in the transition law for the tempered stable OU process with $\alpha\in(1,2)$. As a result, the simulation method by increments for this case does not carry the advantage of being exact \cite{kawai2012infinite}. Moreover, the problem of sampling of multivariate stable random vectors hinders multidimensional generalisations of the recursive scheme based on the Markovian representation \eqref{OU Markovian formula} in practice.

Before moving on, we briefly mention that stochastic integral processes with respect to L\'evy processes of type G also admit shot noise representations based on \eqref{type G series}. Thus, similar Poisson truncation schemes can be used for their sample path generation. Mean-squared error bounds for the truncation of shot noise representation for such stochastic integral processes are provided in \cite{rosinski1991class,wiktorsson2002simulation}, including for schemes based upon subordinated Gaussian representations of L\'evy processes of type G.

\subsection{Simulating infinitely divisible fields}\label{section Levy-driven random fields}

We consider the real harmonisable multifractional L\'evy motion (RHMLM) \cite{dejean2005fracsim,lacaux2004series} as an example of the simulation of an infinitely divisible field via a shot noise representation. The RHMLM is defined by
\[
	X(\mathbf x) = \int_{\mathbb R^d} \frac{e^{-i\langle \mathbf x,\bm\xi\rangle} - 1}{\|\bm\xi\|^{h(\mathbf x) + d/2}} \,L(d\bm\xi),\qquad \mathbf x\in\mathbb R^d,
\]
where $L(d\bm\xi)$ is a L\'evy random measure without Gaussian components with control measure $\nu(d\mathbf z)$ on $\mathbb C$ and the function $h:\mathbb R^d\to(0,1)$ varies in the place of the H\"older exponent.

If the control measure $\nu(d\mathbf z)$ is finite, then a shot noise representation for the RHMLM is given by \cite{dejean2005fracsim}
\begin{equation}\label{RHMLM series}
	\{X(\mathbf x):\mathbf x\in\mathbb R^d\} \stackrel{\mathscr L}{=} \left\{\sum_{k=1}^{+\infty}f\left(\mathbf x,\left(\frac{\Gamma_k}{c_d\nu(\mathbb 
	C)}\right)^{1/d}U_k\right)Z_k : \mathbf x\in\mathbb R^d\right\},
\end{equation}
where $\{\Gamma_k\}_{k\in\mathbb N}$ is a sequence of arrival times of the standard Poisson process, $\{U_k\}_{k\in\mathbb N}$ is a sequence of iid uniform random vectors on $S^{d-1}$, and $\{Z_k\}_{k\in\mathbb N}$ is a sequence of iid random vectors distributed according to $\nu(d\mathbf z)/\nu(\mathbb C)$ such that all random sequences are mutually independent, and $c_d \coloneqq 2\pi^{d/2}/(\Gamma(d/2)d)$. Note that the shot noise representation \eqref{RHMLM series} is distinct from the shot noise representations we have seen thus far. In particular, for the generalised shot noise representation \eqref{generalised shot noise series}, the Poisson arrival times simulate the L\'evy measure through a random walk on the space of jump sizes, while in \eqref{RHMLM series}, the Poisson arrival times simulate a random walk over the space of jump positions.
Armed with a shot noise representation \eqref{RHMLM series} for the RHMLM, a straightforward sample path generation method is to simulate the first $n$ summands. Denoting $\{X^{(n)}(\mathbf x):\mathbf x\in\mathbb R^d\}$ as this approximation, it was found in \cite{lacaux2004series} that for $q\ge p$ and every $n\ge q/2 + q(\max_K h)/d + 1$, it holds that
\[
	\mathbb E\left[\left\|X - X^{(n)}\right\|_{p,K}^q\right] \le C_q\frac{D_{n,q}(\min_K h)}{n^{q(\min_K h)/d}},
\]
where $\|\cdot\|_{p,K}$ denotes $L^p$ norm over a compact set $K\subset\mathbb R^d$, $C_q$ is a constant independent of $n$ and $D_{n,q}(y) = \Gamma(n+1 - q/y - qy/d)^{q/2 + qy/d}/\Gamma(n+1)$.

When the control measure is infinite, it was suggested in \cite{lacaux2004series} to separate the RHMLM in terms of large jumps and small jumps, where the former component can be simulated via the shot noise representation \eqref{RHMLM series}, while under certain technical conditions, the latter can be approximated by a Gaussian field in the same vein as Sections \ref{section Gaussian approximation} and \ref{section Levy-driven stochastic integrals}. We remark that the Poisson truncation scheme described in Section \ref{section truncation} can be applied if the removal of the conditioning on the number of jumps is desired.

\subsection{Simulating L\'evy-driven stochastic differential equations}\label{section Levy-driven SDEs}
As we have applied the theory of shot noise representations to approximations of a large class of infinitely divisible processes (Section \ref{section Levy-driven stochastic integrals}), the next natural step is in generalising the method to approximating L\'evy-driven stochastic differential equations (SDEs).
Working towards a Poisson truncation method of approximating L\'evy-driven SDEs, we first establish the setting of multivariate L\'evy-driven SDEs of the form
\begin{equation}\label{jump-diffusion SDE}
	dX_t = \mu(t,X_t)\,dt + \sigma(t,X_t)\,dB_t + \theta(t,X_{t-})\,dL_t,\quad t\ge0,
\end{equation}
where $\{B_t:t\ge0\}$ is an $l$-dimensional Brownian motion, $\{L_t : t\ge0\}$ is an $l$-dimensional L\'evy process with L\'evy measure $\nu(d\mathbf z)$ and the coefficients $\mu:[0,+\infty)\times\mathbb R^d\to \mathbb R^d$ and $\sigma,\theta:[0,+\infty)\times\mathbb R^d\to \mathbb R^{d\times l}$ are continuous. A sufficient condition for the unique existence of the solution to the SDE \eqref{jump-diffusion SDE} is for the coefficients to satisfy linear growth and Lipschitz conditions \cite[Section 6.2]{applebaum2009levy}.

Sample path generation for solutions to the general SDE is more difficult than for stochastic integral processes.  In short, this is because stochastic integral processes discussed in Section \ref{section Levy-driven stochastic integrals} are described explicitly, whether Markovian or not, over the time interval of interest, leading to straightforward implementation and error analysis. In contrast, the solution to the general SDE can only be described implicitly, that is, the state $X_t$ appears in both sides of \eqref{jump-diffusion SDE}, thus requiring a recursive scheme.
The exception is the very rare case when a closed-form solution for the SDE is available, such as the examples of the L\'evy-driven OU process and Dol\'eans-Dade stochastic exponential. We have seen the explicitness of the OU process in Section \ref{section OU process}, which leads to its simulation in the form of a stochastic integral process \eqref{stochastic integral process}. For the latter, consider the L\'evy-driven Dol\'eans-Dade stochastic exponential as described by the SDE $dX_t = X_{t-}\,dL_t$. By It\^o's lemma for discontinuous semimartingales, the explicit solution is available as
\[
	X_t = X_0\exp\left[-t\int_{-1}^{+\infty}z\,\nu(dz) + \int_0^t\int_{-1}^{+\infty}\ln(1+z)\mu(dz,ds)\right],
\]
if the support of the L\'evy measure $\nu(dz)$ is the half-line $(-1,+\infty)$.
The stochastic exponential can be simulated through the shot noise representation of the integral on the right hand side (see \cite[Section 4.1]{kashima2011optimization}).
More recently, exact methods based on rejection sampling have been developed \cite{giesecke2013exact,pollock2016exact}, which can exactly simulate a class of univariate jump-diffusion processes with finite jump intensity without the need for recursion.

As usual with the study of differential equations, the immediate approach when met with the problem of investigating solutions to the SDE \eqref{jump-diffusion SDE} is via the deterministic time-discretisation paradigm.
The shortfall here is that simulation of the underlying L\'evy integrator via increments does not observe jumps, which renders such a method inappropriate for certain practical scenarios, such as in the case of insurance models where individual claims often need to be observed. In order to incorporate the information of individual jumps for approximations of better quality, jump-adapted time-discretisation, where the discretisation of time includes jump timings, is preferred \cite{bruti2007strong,mordecki2008adaptive}.

With jump-adapted methods, the necessity for computing a finite number of time steps requires the underlying jump component to correspond to a finite L\'evy measure. In the case of the infinite L\'evy measure, the only solution is via truncation to obtain a compound Poisson approximation of the jump component \cite{kohatsu2010jump}. When possible, the accuracy of the numerical scheme can be improved via a Gaussian approximation of the discarded jumps \cite{asmussen2001approximations,cohen2007gaussian}. In the case where the discretisation of time includes deterministic and random jump times, careful balance between the order of the scheme for the Gaussian component and the truncation of the L\'evy measure for the jump component is desired, as the overall rate of convergence is only as fast as the slowest component \cite{kohatsu2014optimal}. In the case where the driver is a subordinated L\'evy process, the Euler method in which the subordinator is approximated by truncation of its shot noise representation is studied in \cite{rubenthaler2003improved}. We also mention here the recent emergence of an alternative approach to approximating SDEs driven by Marcus-type L\'evy noise via homogenisation of deterministic maps \cite{chevyrev2020superdiffusive}.

The basic framework for simulating SDEs via Poisson truncation is described as follows \cite{jum2015numerical}. 
Denote $\nu_n(d\mathbf z)$ as the L\'evy measure corresponding to the Poisson truncation of the shot noise series \eqref{generalised shot noise series for Levy process}. For example, if the shot noise series we truncate corresponds to the inverse L\'evy measure method in the isotropic case, then $\nu_n(d\mathbf z) = \mathbbm1_{[\epsilon(n),+\infty)}(\|\mathbf z\|)\nu(d\mathbf z)$, where the cutoff threshold $\epsilon(n)$ for the magnitude of jumps is decreasing towards zero in $n$ by the definition of the kernel \eqref{inverse Levy measure kernel}. Thus, the corresponding approximation for the solution to the L\'evy-driven SDE \eqref{jump-diffusion SDE} is described by the following SDE as a mixture of integrals and a summation
\begin{equation}\label{truncation for jump-diffusion SDE}
	\begin{split}
	X_{t}^{(n)} &= X_0 + \int_0^t\mu\left(s,X_{s}^{(n)}\right)\,ds + \int_0^t\sigma\left(s,X_{s}^{(n)}\right)\,dB_s - \int_0^t\theta\left(s,X_{s}^{(n)}\right)\,ds \sum_{\{k\in\mathbb N:\Gamma_k\le nT\}}c_k\\
	&\quad + \sum_{\{k\in\mathbb N:\Gamma_k\le nT\}}\theta\left(T_k,X_{T_k-}^{(n)}\right)H\left(\frac{\Gamma_k}{T},U_k\right)\mathbbm1_{[0,t]}(T_k).
	\end{split}
\end{equation}
For the numerical scheme to produce sample paths of the solution process to the SDE \eqref{truncation for jump-diffusion SDE}, one may combine the almost surely finite number of jump timings of the Poisson truncation for the jump component and strong Taylor approximation for the drift and diffusion components (see \cite[Section 8.2]{Platen2010Numerical}).

Numerical methods and error analysis for L\'evy-driven SDEs via the truncation of its L\'evy measure in full generality are still on-going research. For other promising directions of future research, we conjecture that the shot noise representation has the potential to offer effective numerical methods for stochastic delay differential equations with jumps \cite{dareiotis2016tamed}, backward stochastic differential equations with jumps  \cite{fujii2019asymptotic,geiss2016simulation} and stochastic partial differential equations with jumps \cite{dunst2012approximate,hausenblas2006numerical}.

\section{Computation of expectations via shot noise representation}\label{section computing expectations}
For the practical use of models involving infinitely divisible laws and related processes, it is often crucial to be able to approximate relevant expectations, such as estimating probabilities and moments.
As mentioned earlier, theoretical discussions have taken advantage of the importance of the initial arrival times \cite{breton2010regularity}. In our context, this fact can be effectively exploited for numerically approximating expectations \cite{imai2010quasi,kawai2012monte}. Moreover, techniques involving changing the underlying probability measure \cite{kawai2006importance,kawai2012monte} are also useful.
In what follows, we discuss some topics relevant to shot noise representations for approximating expectations.

Throughout, we assume that an infinitely divisible random vector of our interest admits the shot noise representation \eqref{generalised shot noise series}. Crucially for the succeeding discussion, recall that the arrival times $\{\Gamma_k\}_{k\in\mathbb N}$ of the standard Poisson process equal in law to successive summations of iid standard exponential random variables, that is,
\begin{equation}\label{interarrival times}
	\{\Gamma_1,\Gamma_2,\Gamma_3,\cdots\} \stackrel{\mathscr L}{=} \left\{\sum\nolimits_{k=1}^1E_k,\sum\nolimits_{k=1}^2E_k,\sum\nolimits_{k=1}^3E_k,\cdots\right\},
\end{equation}
where $\{E_k\}_{k\in\mathbb N}$ is a sequence of iid standard exponential random variables. As the first (inter-)arrival time $E_1$ appears in all summands $\{H(\Gamma_k,U_k)\}_{k\in\mathbb N}$, it accounts for a large portion of the variation of the underlying randomness. Similar statements can be said of the first few exponential interarrival times. We remark that L\'evy processes and more general infinitely divisible processes are within our scope through the arrival times $\{\Gamma_k\}_{k\in\mathbb N}$ in their shot noise representations, such as \eqref{generalised shot noise series for Levy process}, \eqref{nfsm truncation} and \eqref{stable CARMA truncation}. Moreover, as information on individual jumps are perhaps even more crucial in the case of infinitely divisible processes, the following discussions are even more pertinent for the aforementioned more general setting.

\subsection{Effective dimension on interarrival exponentials}\label{section qmc1}

Intuitively, in the case of a shot noise representation \eqref{generalised shot noise series}, the faster the kernel $\|H(\cdot,\bm\xi)\|$ decays, the greater the proportion of variation explained by the first few exponential random variables $\{E_1,\cdots,E_n\}$ \eqref{interarrival times}. This leads to the notion that in many cases, the variation can be captured by only considering a lower dimension of the exponential interarrival times $\{E_k\}_{k\in\mathbb N}$. This idea is formalised in \cite{imai2010quasi} through the cumulative explanatory ratio (CER), which is interpreted as the proportion of the variance explained by the first few dimensions. Where $X$ is an infinitely divisible random vector in $\mathbb R^d$ without a Gaussian component and $f:\mathbb R^d\to\mathbb R$ is a continuous function, suppose we want to compute $\mathbb E[f(X)]$ provided that $\text{Var}(f(X))$ is finite. The CER associated with the first $n$ interarrival times is defined by
\begin{equation}\label{CER}
	\text{CER}_n \coloneqq \frac{\text{Cov}(X_f(n),X_f(0))}{\text{Var}(X_f(0))},
\end{equation}
where, for natural numbers $n\in\{1,\cdots,N\}$,
\[
	X_f(n)\coloneqq f\left(\sum\nolimits_{k=1}^N\left(H\left(\sum\nolimits_{l=1}^{k\wedge n} E_l + \sum\nolimits_{l=k\wedge n+1}^k E_l',U_k\right) - c_k\right)\right) \stackrel{\mathscr L}{=} f\left(\sum\nolimits_{k=1}^N\left(H(\Gamma_k,U_k) - c_k\right)\right)\eqqcolon X_f(0),
\]
where $\{E_k'\}_{k\in\mathbb N}$ is an iid copy of the sequence of exponential interarrival times $\{E_k\}_{k\in\mathbb N}$.

A high CER \eqref{CER} corresponds to a low effective dimension structure of the shot noise representation. This is a desirable property, as rather than using Monte Carlo methods with the convergence rate $\mathcal O(n^{-1/2})$, quasi-Monte Carlo methods can be reliably used on the first few interarrival times to achieve the faster convergence rate $\mathcal O(n^{-1}(\ln n)^d)$ in practice. This is in contrast to problems with high-dimensionality, for which the speed improvement may often be too marginal to justify the use of quasi-Monte Carlo methods for such problems.
The CER \eqref{CER} for stable random variables and a modified CER for the randomised quasi-Monte Carlo method for the gamma random variable were investigated in \cite{imai2010quasi}, which were found to be remarkably high with only the first few terms of their shot noise series. Thus, quasi-Monte Carlo methods are applicable for greater accuracy in expectation computations.

\subsection{Stratification on interarrival exponentials}\label{section stratification}
We describe the variance reduction method of stratifying the exponential interarrival times $\{E_k\}_{k\in\mathbb N}$ \eqref{interarrival times}, as investigated in \cite{kawai2012monte}. We present the technique of stratified sampling when computing a random variable $F$ involving a shot noise representation, say \eqref{generalised shot noise series}, which is built upon the sequence $\{E_k\}_{k\in\mathbb N}$.
For simplicity, we only consider the stratification of the first interarrival time $E_1$. Fix $M\in\mathbb N$ and partition $(0,+\infty) = \cup_{m=1}^M B_m$ so $B_1 = (0,b_1],B_2 = (b_1,b_2],\cdots,B_{M} = (b_{M-1},+\infty)$ with $0<b_1<b_2<\cdots< b_{M-1}<+\infty$, such that $p_m\coloneqq \mathbb P(E_1\in B_m) = 1/M$ for every $m\in\{1,\cdots,M\}$. For every $m$, let $\{F_{m,k}\}_{k\in\mathbb N}$ be a sequence of iid random variables with the distribution $\mathscr{L}(F | E_1\in B_m)$. Then, where $n_m$ denotes the number of samples from the $m$-th stratum, the random variable
\[
	\sum_{m=1}^M p_m\frac{1}{n_m}\sum_{k=1}^{n_m} F_{m,k}
\]
is an unbiased estimator for the expectation $\mathbb E[F]$, with estimator variance
\[
	\text{Var}_{\mathbb P}\left(\sum_{m=1}^M p_m\frac{1}{n_m}\sum_{k=1}^{n_m}F_{m,k}\right) = \frac1n \sum_{m=1}^M p_m^2\frac{\alpha_m^2}{q_m},
\]
where $q_m\coloneqq n_m/n$. The minimal variance $(\sum_{m=1}^M p_m\sigma_m)^2$ is achieved with the allocation $q_m = p_m\sigma_m/(\sum_{m=1}^M p_m\sigma_m)$, where $\sigma^2_m \coloneqq \text{Var}_{\mathbb P}(F | E_1\in B_m)$ is the stratum variance. It should be noted that the stratified sampling technique can be extended to further interarrival times. However, the extension may be computationally taxing with only marginal returns due to the low effective dimension of the shot noise representation (Section \ref{section qmc1}) and the rapid growth of strata with the Monte Carlo dimension \cite{imai2010quasi}.

\subsection{Importance sampling and control variates on interarrival exponentials}\label{section qmc1 2}
Next, we outline the variance reduction methods of control variates and importance sampling \cite{kawai2012monte}. We assume the same simplified scenario considered in Section \ref{section stratification}. Fix $n\in\mathbb N$ and define $E^{(n)} \coloneqq [E_1,\cdots,E_n]^\intercal$ of \eqref{interarrival times} and $\bm\lambda \coloneqq [\lambda_1,\cdots,\lambda_n]^\intercal\in(-\infty,1)^n$. We parametrise the first $n$ interarrival times as $\{E_k/(1-\lambda_k)\}_{k\in\{1,\cdots,n\}}$, thus in turn we parametrise the random variable of interest $F$ as $F(\bm\lambda)$. Define a family of probability measures $\{\mathbb Q_{\bm\lambda}\}_{\bm\lambda\in(-\infty,1)^n}$ by
\[
	\frac{d\mathbb Q_{\bm\lambda}}{d\mathbb P}\big|_{\sigma(E^{(n)})} \coloneqq \frac{e^{\langle\bm\lambda,E^{(n)}\rangle}}{\mathbb E_{\mathbb P}[e^{\langle\bm\lambda,E^{(n)}\rangle}]} = \prod_{k=1}^n (1-\lambda_k)e^{\lambda_k E_k},\quad\mathbb P\text{-a.s.}
\]
Note that $\mathbb Q_{\bm\lambda}(E_k\in B) = \mathbb P(E_k/(1-\lambda_k)\in B)$ and $\mathbb Q_{\bm\lambda}(F(0)\in B) = \mathbb P(F(\bm\lambda)\in B)$ for every $B\in\mathcal B(\mathbb R)$. It holds that for every $\bm\theta\coloneqq[\theta_1,\cdots,\theta_n]^\intercal\in\mathbb R^n$,
\begin{align}
	\mathbb E_{\mathbb P}[F(0)] &= \mathbb E_{\mathbb P}\left[F(0) - \left\langle\bm\theta, E^{(n)} - \mathbb E_{\mathbb P}\left[E^{(n)}\right]\right\rangle\right] \label{control variates}\\
	&= \mathbb E_{\mathbb Q_{\bm\lambda}}\left[\frac{d\mathbb P}{d\mathbb Q_{\bm\lambda}}\big|_{\sigma(E^{(n)})}\left(F(0) - \left\langle\bm\theta,E^{(n)} - \mathbb E_{\mathbb P}\left[E^{(n)}\right]\right\rangle\right)\right] \label{importance sampling}\\
	\notag &= \mathbb E_{\mathbb P}\left[\left(\prod_{k=1}^n \frac{\exp\left[-\frac{\lambda_k}{1-\lambda_k}E_k\right]}{1-\lambda_k}\right)\left(F(\bm\lambda) - \sum_{k=1}^n\theta_k\left(\frac{E_k}{1-\lambda_k} - 1\right)\right)\right].
\end{align}
Introducing the subtraction term in the expectation in \eqref{control variates} with equality is referred to as the method of control variates. The change of measure in \eqref{importance sampling} is referred to as the method of importance sampling. The estimator variance is given by
\begin{align*}
	V(\bm\lambda,\bm\theta) &\coloneqq \text{Var}_{\mathbb Q_{\bm\lambda}}\left(\frac{d\mathbb P}{d\mathbb Q_{\bm\lambda}}\big|_{\sigma(E^{(n)})}\left(F(0) - \left\langle\bm\theta,E^{(n)} - \mathbb E_{\mathbb P}\left[E^{(n)}\right]\right\rangle\right)\right)\\
	&= \mathbb E_{\mathbb P}\left[\left(\prod_{k=1}^n\frac{e^{-\lambda_k E_k}}{1-\lambda_k}\right)\left(F(0) - \left\langle F(0),E^{(n)} - 	\mathbf1_n\right\rangle\right)^2\right] - (\mathbb E_{\mathbb P}[F(0)])^2,
\end{align*}
where $\mathbf1_n\coloneqq [1,\cdots,1]^\intercal\in\mathbb R^n$.
For example, for the method of control variates ($\bm\lambda=0$) applied only to the first interarrival time $E_1$ ($n=1$), the optimal parameter can be derived as $\theta^* = \text{Cov}_{\mathbb P}(F(0),E_1)$. For the method of importance sampling alone ($\bm\theta=0$), stochastic approximation techniques may be employed for searching an optimal parameter.

\subsection{Importance sampling on all individual jumps}
Yet another importance sampling method, different from the one considered in Section \ref{section qmc1 2}, can be constructed using more information of the sample path via density transformations between individual jumps \cite{kawai2006importance}. Let $(\{X_t:t\ge0\},\mathbb P)$ and $(\{X_t:t\ge0\},\mathbb Q)$ be L\'evy processes in $\mathbb R^d$ characterised by triples $(\mathbf\gamma_{\mathbb P},A,\nu_{\mathbb P})$ and $(\mathbf\gamma_{\mathbb Q},A,\nu_{\mathbb Q})$, respectively. Given some technical conditions on the distance between the two laws, it holds that the probability measures $\mathbb P$ and $\mathbb Q$ are equivalent and
\[
	\frac{d\mathbb P}{d\mathbb Q}\Big|_{\mathscr F_t} = e^{U_t}, \quad \mathbb Q\text{-a.s.},
\]
where the stochastic process $\{U_t:t\ge0\}$ in $\mathbb R$ satisfies
\begin{equation}\label{vred U}
	U_t = \langle \bm\eta,\widetilde X_t\rangle - \frac t2 \langle \bm\eta,A\bm\eta\rangle - t\langle \mathbf\gamma_{\mathbb Q},\bm\eta\rangle + \lim_{\epsilon\to0}\left[\sum_{(s,\|X_s - X_{s-}\|)\in(0,t]\times(\epsilon,+\infty)}\varphi(X_s - X_{s-}) - t\int_{\|\mathbf z\|>\epsilon}\left(e^{\varphi(\mathbf z)}-1\right)\,\nu_{\mathbb Q}(d\mathbf z)\right],\quad \mathbb Q\text{-a.s.},
\end{equation}
where $\widetilde X_t = X_t - \sum_{s\in(0,t]}(X_s - X_{s-})$, $\varphi \coloneqq \ln d\nu_{\mathbb P}/d\nu_{\mathbb Q}$ and the vector $\bm\eta\in\mathbb R^d$ satisfies $\mathbf\gamma_{\mathbb Q} - \mathbf\gamma_{\mathbb P} - \int_{\|\mathbf z\|\le 1}\mathbf z\,(\nu_{\mathbb Q} - \nu_{\mathbb P})(d\mathbf z) = A\bm\eta$.
Moreover, the stochastic process $\{U_t:t\ge0\}$ is uniformly convergent in $t$ on any bounded interval $\mathbb Q\text{-a.s.}$ and $\mathbb E_{\mathbb Q}[e^{U_t}] = \mathbb E_{\mathbb P}[e^{-U_t}] = 1$ for every $t\ge0$.

For simplicity, we present the one-dimensional case. Suppose we want to evaluate $\mathbb E_{\mathbb P}[F]$, where $F$ is a random variable involving the sample path of the L\'evy process $\{X_t:t\in[0,T]\}$. As it holds that $\mathbb E_{\mathbb P}[F] = \mathbb E_{\mathbb Q}[e^{U_T}F]$, we can estimate our quantity of interest via a Monte Carlo iteration of the latter, that is,
\[
	\lim_{n\to+\infty}\frac1n \sum_{k=1}^n e^{U_{k,T}}F_{k} = \mathbb E_{\mathbb P}[F],\quad\mathbb Q\text{-a.s.},
\]
where $\{F_k\}_{k\in\mathbb N}$ is a sequence of iid copies of the random variable $F$ and $\{U_{k,t}:t\in[0,T]\}_{k\in\mathbb N}$ is a sequence of iid copies of $\{U_t:t\in[0,T]\}$. Observe from \eqref{vred U} that generating the sample path $\{U_t:t\ge0\}$ requires the jumps of the L\'evy process, so generally the implementation of the computation significantly benefits from, or more precisely, often requires shot noise representation, for instance, infinitely divisible processes driven by the L\'evy process $\{X_t:t\ge0\}$.

\section{Concluding remarks}\label{section conclusion}
In this survey, we have summarised shot noise representation with a view towards sampling infinitely divisible laws and generating sample paths of related processes. In particular, we reviewed the important aspects of shot noise representation through the rather scenic route of the L\'evy-It\^o decomposition approach. We provided shot noise representations of various popular laws and stochastic processes in the literature. Through our description of the truncation of shot noise representation, a general and systematic method for simulation and computation of expectations was discussed. Examples of simulation recipes were provided, and the key results for error analysis and our numerical demonstrations should provide the confidence that truncation of shot noise representation does not simply satisfy the need to be accurate, but also the desire for a straightforward and computationally feasible approach for simulation in practice.

We hope that the present survey makes clear the practicality of numerical methods for simulating infinitely divisible laws and related processes based on shot noise representations, and encourages future development into expanding the technique. We reiterate that the approximation of L\'evy-driven SDEs and more general shot noise processes via truncation of shot noise representation is still the subject of further research. As mentioned previously, future directions for this area of research include numerical schemes for stochastic delay differential equations, backward stochastic differential equations and stochastic partial differential equations with jumps based on shot noise representations. As shot noise representation not only yields a viable method of sample path generation for a wide-range of stochastic processes, but can also provide insights into their properties, we expect further studies of such stochastic processes will continue to invoke shot noise representation techniques. For this reason, the investigation of shot noise representations and their truncation will remain in priority for the foreseeable future.


{\small
	\bibliographystyle{abbrv}
	\bibliography{master}

\begin{thebibliography}{100}

\bibitem{almasary2017approximate}
Z.~Al~Masry, S.~Mercier, and G.~Verdier.
\newblock Approximate simulation techniques and distribution of an extended
  gamma process.
\newblock {\em Methodology and Computing in Applied Probability}, 19:213--235,
  2017.

\bibitem{applebaum2009levy}
D.~Applebaum.
\newblock {\em {L{\'e}vy Processes and Stochastic Calculus}}.
\newblock Cambridge University Press, United Kingdom, 2009.

\bibitem{asmussen2001approximations}
S.~Asmussen and J.~Rosi{\'n}ski.
\newblock {Approximations of small jumps of L{\'e}vy processes with a view
  towards simulation}.
\newblock {\em Journal of Applied Probability}, 38(2):482--493, 2001.

\bibitem{basse2013uniform}
A.~Basse-O'Connor and J.~Rosi{\'n}ski.
\newblock {On the uniform convergence of random series in Skorohod space and
  representations of c{\`a}dl{\`a}g infinitely divisible processes}.
\newblock {\em The Annals of Probability}, 41(6):4317--4341, 2013.

\bibitem{beenakker1999photon}
C.~Beenakker and M.~Patra.
\newblock Photon shot noise.
\newblock {\em Modern Physics Letters B}, 13(11):337--347, 1999.

\bibitem{beenakker2003quantum}
C.~Beenakker and C.~Sch{\"o}nenberger.
\newblock Quantum shot noise.
\newblock {\em Physics Today}, 56(5):37, 2003.

\bibitem{beenakker1992suppression}
C.~W.~J. Beenakker and M.~B\"uttiker.
\newblock Suppression of shot noise in metallic diffusive conductors.
\newblock {\em Physics Review B: Condensed Matter and Materials Physics},
  46(3):1889--1892, 1992.

\bibitem{bentkus1996bounds}
V.~Bentkus, F.~G{\"o}tze, and V.~Paulauskas.
\newblock {Bounds for the accuracy of Poissonian approximations of stable
  laws}.
\newblock {\em Stochastic Processes and their Applications}, 65(1):55--68,
  1996.

\bibitem{bentkus2001levy}
V.~Bentkus, A.~Juozulynas, and V.~Paulauskas.
\newblock {L{\'e}vy--LePage series representation of stable vectors:
  convergence in variation}.
\newblock {\em Journal of Theoretical Probability}, 14(4):949--978, 2001.

\bibitem{bertoin1998levy}
J.~Bertoin.
\newblock {\em {L{\'e}vy Processes}}.
\newblock {Cambridge University Press}, 1998.

\bibitem{bianchi2011tempered}
M.~L. Bianchi, S.~T. Rachev, Y.~S. Kim, and F.~J. Fabozzi.
\newblock {Tempered infinitely divisible distributions and processes}.
\newblock {\em Theory of Probability \& Its Applications}, 55(1):2--26, 2011.

\bibitem{bondesson1982simulation}
L.~Bondesson.
\newblock On simulation from infinitely divisible distributions.
\newblock {\em Advances in Applied Probability}, 14(4):855--869, 1982.

\bibitem{bremaud2002power}
P.~Br{\'e}maud and L.~Massouli{\'e}.
\newblock {Power spectra of general shot noises and Hawkes point processes with
  a random excitation}.
\newblock {\em Advances in Applied Probability}, 34(1):205--222, 2002.

\bibitem{breton2010regularity}
J.-C. Breton.
\newblock Regularity of the laws of shot noise series and of related processes.
\newblock {\em Journal of Theoretical Probability}, 23:21--38, 2010.

\bibitem{brix1999generalized}
A.~Brix.
\newblock {Generalized gamma measures and shot-noise Cox processes}.
\newblock {\em Advances in Applied Probability}, 31(4):929--953, 1999.

\bibitem{bruti2007strong}
N.~Bruti-Liberati and E.~Platen.
\newblock Strong approximations of stochastic differential equations with
  jumps.
\newblock {\em Journal of Computational and Applied Mathematics},
  205(2):982--1001, 2007.

\bibitem{carnaffan2017cusping}
S.~Carnaffan and R.~Kawai.
\newblock {Cusping, transport and variance of solutions to generalized
  Fokker--Planck equations}.
\newblock {\em Journal of Physics A: Mathematical and Theoretical},
  50(24):245001, 2017.

\bibitem{carnaffan2017solving}
S.~Carnaffan and R.~Kawai.
\newblock {Solving multidimensional fractional Fokker--Planck equations via
  unbiased density formulas for anomalous diffusion processes}.
\newblock {\em SIAM Journal on Scientific Computing}, 39(5):B886--B915, 2017.

\bibitem{carnaffan2019analytic}
S.~Carnaffan and R.~Kawai.
\newblock Analytic model for transient anomalous diffusion with highly
  persistent correlations.
\newblock {\em Physics Review E}, 99(6):062120, 2019.

\bibitem{carr2002fine}
P.~Carr, H.~Geman, D.~B. Madan, and M.~Yor.
\newblock {The fine structure of asset returns: An empirical investigation}.
\newblock {\em The Journal of Business}, 75(2):305--332, 2002.

\bibitem{chevyrev2020superdiffusive}
I.~Chevyrev, P.~K. Friz, A.~Korepanov, and I.~Melbourne.
\newblock Superdiffusive limits for deterministic fast-slow dynamical systems.
\newblock {\em Probability Theory and Related Fields}, 178(3):735--770, 2020.

\bibitem{cohen2007gaussian}
S.~Cohen and J.~Rosi{\'n}ski.
\newblock {Gaussian approximation of multivariate L{\'e}vy processes with
  applications to simulation of tempered stable processes}.
\newblock {\em Bernoulli}, 13(1):195--210, 2007.

\bibitem{tankov2003financial}
R.~Cont and P.~Tankov.
\newblock {\em {Financial Modelling with Jump Processes}}.
\newblock Chapman \& Hall / CRC Press, 2003.

\bibitem{daley1971definition}
D.~Daley.
\newblock The definition of a multi-dimensional generalization of shot noise.
\newblock {\em Journal of Applied Probability}, 8(1):128--135, 1971.

\bibitem{dareiotis2016tamed}
K.~Dareiotis, C.~Kumar, and S.~Sabanis.
\newblock {On tamed Euler approximations of SDEs driven by L{\'e}vy noise with
  applications to delay equations}.
\newblock {\em SIAM Journal on Numerical Analysis}, 54(3):1840--1872, 2016.

\bibitem{davydov2012convergence}
Y.~Davydov and C.~Dombry.
\newblock {On the convergence of LePage series in Skorokhod space}.
\newblock {\em Statistics \& Probability Letters}, 82(1):145--150, 2012.

\bibitem{dejean2005fracsim}
S.~D{\'e}jean and S.~Cohen.
\newblock {FracSim: An R package to simulate multifractional L{\'e}vy motions}.
\newblock {\em Journal of Statistical Software}, 14(1):1--19, 2005.

\bibitem{dunst2012approximate}
T.~Dunst, E.~Hausenblas, and A.~Prohl.
\newblock {Approximate Euler method for parabolic stochastic partial
  differential equations driven by space-time L{\'e}vy noise}.
\newblock {\em SIAM Journal on Numerical Analysis}, 50(6):2873--2896, 2012.

\bibitem{echternach2013photon}
P.~Echternach, K.~Stone, C.~Bradford, P.~Day, D.~Wilson, K.~Megerian,
  N.~Llombart, and J.~Bueno.
\newblock Photon shot noise limited detection of terahertz radiation using a
  quantum capacitance detector.
\newblock {\em Applied Physics Letters}, 103(5):053510, 2013.

\bibitem{eliazar2005nonlinear}
I.~Eliazar and J.~Klafter.
\newblock On the nonlinear modeling of shot noise.
\newblock {\em Proceedings of the National Academy of Sciences},
  102(39):13779--13782, 2005.

\bibitem{ferguson1972representation}
T.~S. Ferguson and M.~J. Klass.
\newblock {A representation of independent increment processes without Gaussian
  components}.
\newblock {\em The Annals of Mathematical Statistics}, 43(5):1634--1643, 1972.

\bibitem{fujii2019asymptotic}
M.~Fujii and A.~Takahashi.
\newblock {Asymptotic expansion for forward-backward SDEs with jumps}.
\newblock {\em Stochastics}, 91(2):175--214, 2019.

\bibitem{geiss2016simulation}
C.~Geiss and C.~Labart.
\newblock {Simulation of BSDEs with jumps by Wiener chaos expansion}.
\newblock {\em Stochastic Processes and their Applications}, 126(7):2123--2162,
  2016.

\bibitem{giesecke2013exact}
K.~Giesecke and D.~Smelov.
\newblock Exact sampling of jump diffusions.
\newblock {\em Operations Research}, 61(4):894--907, 2013.

\bibitem{gilchrist1975shot}
J.~H. Gilchrist and J.~B. Thomas.
\newblock A shot process with burst properties.
\newblock {\em Advances in Applied Probability}, 7(3):527--541, 1975.

\bibitem{giraitis1991shot}
L.~Giraitis and D.~Surgailis.
\newblock {On shot noise processes attracted to fractional L{\'e}vy motion}.
\newblock In {\em Stable Processes and Related Topics}, pages 261--273.
  Springer, 1991.

\bibitem{grothe2013vine}
O.~Grothe and S.~Nicklas.
\newblock {Vine constructions of L{\'e}vy copulas}.
\newblock {\em Journal of Multivariate Analysis}, 119:1--15, 2013.

\bibitem{hackmann2018karhunen}
D.~Hackmann.
\newblock {Karhunen--Lo{\`e}ve expansions of L{\'e}vy processes}.
\newblock {\em Communications in Statistics-Theory and Methods},
  47(23):5675--5687, 2018.

\bibitem{hausenblas2006numerical}
E.~Hausenblas and I.~Marchis.
\newblock {A numerical approximation of parabolic stochastic partial
  differential equations driven by a Poisson random measure}.
\newblock {\em BIT Numerical Mathematics}, 46(4):773, 2006.

\bibitem{hawkes1974cluster}
A.~G. Hawkes and D.~Oakes.
\newblock A cluster process representation of a self-exciting process.
\newblock {\em Journal of Applied Probability}, 11(3):493--503, 1974.

\bibitem{houdre2006fractional}
C.~Houdr{\'e} and R.~Kawai.
\newblock On fractional tempered stable motion.
\newblock {\em {Stochastic Processes and their Applications}},
  116(8):1161--1184, 2006.

\bibitem{houdre2007layered}
C.~Houdr{\'e} and R.~Kawai.
\newblock On layered stable processes.
\newblock {\em Bernoulli}, 13(1):252--278, 2007.

\bibitem{hult2007extremal}
H.~Hult and F.~Lindskog.
\newblock {Extremal behavior of stochastic integrals driven by regularly
  varying L{\'e}vy processes}.
\newblock {\em The Annals of Probability}, 35(1):309--339, 2007.

\bibitem{imai2010quasi}
J.~Imai and R.~Kawai.
\newblock {Quasi-Monte Carlo method for infinitely divisible random vectors via
  series representations}.
\newblock {\em SIAM Journal on Scientific Computing}, 32(4):1879--1897, 2010.

\bibitem{imai2011finite}
J.~Imai and R.~Kawai.
\newblock On finite truncation of infinite shot noise series representation of
  tempered stable laws.
\newblock {\em {Physica A: Statistical Mechanics and its Applications}},
  390(23-24):4411--4425, 2011.

\bibitem{imai2013numerical}
{J. Imai and R. Kawai}.
\newblock {Numerical inverse L{\'e}vy measure method for infinite shot noise
  series representation}.
\newblock {\em {Journal of Computational and Applied Mathematics}},
  253:264--283, 2013.

\bibitem{janicki1993simulation}
A.~Janicki and A.~Weron.
\newblock {\em {Simulation and Chaotic Behavior of {$\alpha$}-Stable Stochastic
  Processes}}.
\newblock CRC Press, 1993.

\bibitem{jum2015numerical}
E.~Jum.
\newblock {\em {Numerical Approximation of Stochastic Differential Equations
  Driven by L{\'e}vy Motion with Infinitely Many Jumps}}.
\newblock PhD thesis, University of Tennessee, 2015.

\bibitem{kahle2016degradation}
W.~Kahle, S.~Mercier, and C.~Paroissin.
\newblock {\em {Degradation Processes in Reliability}}.
\newblock John Wiley \& Sons, 2016.

\bibitem{kallenberg1974series}
O.~Kallenberg.
\newblock {Series of random processes without discontinuities of the second
  kind}.
\newblock {\em The Annals of Probability}, 2(9):729--737, 1974.

\bibitem{kashima2011optimization}
K.~Kashima and R.~Kawai.
\newblock An optimization approach to weak approximation of stochastic
  differential equations with jumps.
\newblock {\em Applied Numerical Methods}, 61(5):641--650, 2011.

\bibitem{kawai2006importance}
R.~Kawai.
\newblock {An importance sampling method based on the density transformation of
  L{\'e}vy processes}.
\newblock {\em Monte Carlo Methods and Applications}, 12(2):171 -- 186, 2006.

\bibitem{kawai2016higher}
R.~Kawai.
\newblock {Higher order fractional stable motion: hyperdiffusion with heavy
  tails}.
\newblock {\em {Journal of Statistical Physics}}, 165(1):126--152, 2016.

\bibitem{kawai2017sample}
R.~Kawai.
\newblock {Sample path generation of L{\'e}vy-driven continuous-time
  autoregressive moving average processes}.
\newblock {\em Methodology and Computing in Applied Probability},
  19(1):175--211, 2017.

\bibitem{idsim}
R.~Kawai.
\newblock A general approach to sample path generation of infinitely divisible
  processes via shot noise representation.
\newblock {\em Statistics \& Probability Letters}, 174:109091, 2021.

\bibitem{kawai2012monte}
R.~Kawai and J.~Imai.
\newblock {On Monte Carlo and Quasi-Monte Carlo methods for series
  representation of infinitely divisible laws}.
\newblock In {\em Monte Carlo and Quasi-Monte Carlo Methods 2010}, pages
  471--486. Springer, 2012.

\bibitem{kawai2011exact}
R.~Kawai and H.~Masuda.
\newblock {Exact discrete sampling of finite variation tempered stable
  Ornstein--Uhlenbeck processes}.
\newblock {\em {Monte Carlo Methods and Applications}}, 17(3):279--300, 2011.

\bibitem{kawai2012infinite}
R.~Kawai and H.~Masuda.
\newblock {Infinite variation tempered stable Ornstein--Uhlenbeck processes
  with discrete observations}.
\newblock {\em Communications in Statistics-Simulation and Computation},
  41(1):125--139, 2012.

\bibitem{kluppelberg2004fractional}
C.~Kl{\"u}ppelberg and C.~K{\"u}hn.
\newblock {Fractional Brownian motion as a weak limit of Poisson shot noise
  processes{--}with applications to finance}.
\newblock {\em Stochastic Processes and their Applications}, 113(2):333--351,
  2004.

\bibitem{kohatsu2014optimal}
A.~Kohatsu-Higa, S.~Ortiz-Latorre, and P.~Tankov.
\newblock {Optimal simulation schemes for L{\'e}vy driven stochastic
  differential equations}.
\newblock {\em Mathematics of Computation}, 83(289):2293--2324, 2014.

\bibitem{kohatsu2010jump}
A.~Kohatsu-Higa and P.~Tankov.
\newblock {Jump-adapted discretization schemes for L{\'e}vy-driven SDEs}.
\newblock {\em Stochastic Processes and their Applications},
  120(11):2258--2285, 2010.

\bibitem{koponen1995analytic}
I.~Koponen.
\newblock {Analytic approach to the problem of convergence of truncated
  L{\'e}vy flights towards the Gaussian stochastic process}.
\newblock {\em Physical Review E}, 52(1):1197--1199, 1995.

\bibitem{kotz2001laplace}
S.~Kotz, T.~J. Kozubowski, and K.~Podg^^c3^^b3rski.
\newblock {\em {The Laplace Distribution and Generalizations: A Revisit with
  Applications to Communications, Economics, Engineering, and Finance}}.
\newblock {Birkh{\"a}user}, Boston, 2001.

\bibitem{lacaux2004series}
C.~Lacaux.
\newblock {Series representation and simulation of multifractional L{\'e}vy
  motions}.
\newblock {\em Advances in Applied Probability}, 36(1):171--197, 2004.

\bibitem{lane1984central}
J.~A. Lane.
\newblock {The central limit theorem for the Poisson shot-noise process}.
\newblock {\em Journal of Applied Probability}, 21(2):287--301, 1984.

\bibitem{lecourtois2018some}
O.~Le~Courtois.
\newblock Some further results on the tempered multistable approach.
\newblock {\em Asia-Pacific Financial Markets}, 25:87--109, 2018.

\bibitem{leguevel2012ferguson}
R.~Le~Gu{\'e}vel and J.~L. V{\'e}hel.
\newblock {A Ferguson--Klass--LePage series representation of multistable
  multifractional motions and related processes}.
\newblock {\em Bernoulli}, 18(4):1099--1127, 2012.

\bibitem{lemke2015fully}
T.~Lemke, M.~Riabiz, and S.~J. Godsill.
\newblock {Fully Bayesian inference for {$\alpha$}-stable distributions using a
  Poisson series representation}.
\newblock {\em Digital Signal Processing}, 47:96 --115, 2015.

\bibitem{lepage1980multidimensional}
R.~LePage.
\newblock {Multidimensional infinitely divisidle variables and processes Part
  II}.
\newblock In {\em Probability in Banach Spaces III}, pages 279--284.
  Springer-Verlag, 1980.

\bibitem{lepage1981convergence}
R.~LePage, M.~Woodroofe, and J.~Zinn.
\newblock Convergence to a stable distribution via order statistics.
\newblock {\em The Annals of Probability}, 9(4):624--632, 1981.

\bibitem{mantegna1994stochastic}
R.~N. Mantegna and H.~E. Stanley.
\newblock {Stochastic process with ultraslow convergence to a Gaussian: the
  truncated L{\'e}vy flight}.
\newblock {\em Physical Review Letters}, 73(22):2946--2949, 1994.

\bibitem{massing2018simulation}
T.~Massing.
\newblock {Simulation of Student--L{\'e}vy processes using series
  representations}.
\newblock {\em Computational Statistics}, 33:1649--1685, 2018.

\bibitem{moller2006approximate}
J.~M{\o}ller and J.~G. Rasmussen.
\newblock {Approximate simulation of Hawkes processes}.
\newblock {\em Methodology and Computing in Applied Probability}, 8(1):53--64,
  2006.

\bibitem{moller2005generalised}
J.~M{\o}ller and G.~L. Torrisi.
\newblock {Generalised shot noise Cox processes}.
\newblock {\em Advances in Applied Probability}, 37(1):48--74, 2005.

\bibitem{mordecki2008adaptive}
E.~Mordecki, A.~Szepessy, R.~Tempone, and G.~E. Zouraris.
\newblock {Adaptive weak approximation of diffusions with jumps}.
\newblock {\em SIAM Journal on Numerical Analysis}, 46(4):1732--1768, 2008.

\bibitem{pang2019nonstationary}
G.~Pang and M.~S. Taqqu.
\newblock {Nonstationary self-similar Gaussian processes as scaling limits of
  power-law shot noise processes and generalizations of fractional Brownian
  motion}.
\newblock {\em High Frequency}, 2(2):95--112, 2019.

\bibitem{peters2015advances}
G.~W. Peters and P.~V. Shevchenko.
\newblock {\em {Advances in Heavy Tailed Risk Modeling: A Handbook of
  Operational Risk}}.
\newblock John Wiley \& Sons, 2015.

\bibitem{Platen2010Numerical}
E.~Platen and N.~Bruti-Liberti.
\newblock {\em Numerical Solution of Stochastic Differential Equations with
  Jumps in Finance}.
\newblock Springer-Verlag, 2010.

\bibitem{pollock2016exact}
M.~Pollock, A.~M. Johansen, and G.~O. Roberts.
\newblock {On the exact and {$\varepsilon$}-strong simulation of (jump)
  diffusions}.
\newblock {\em Bernoulli}, 22(2):794--856, 2016.

\bibitem{rajput1989spectral}
B.~S. Rajput and J.~Rosi{\'n}ski.
\newblock Spectral representations of infinitely divisible processes.
\newblock {\em Probability Theory and Related Fields}, 82(3):451--487, 1989.

\bibitem{resnick1976extremal}
S.~I. Resnick.
\newblock An extremal decomposition of a process with independent, stationary
  increments.
\newblock {\em Technical Report 79, Department of Statistics, Stanford
  University}, 1976.

\bibitem{reznikov1998quantum}
M.~Reznikov, R.~De~Picciotto, M.~Heiblum, D.~Glattli, A.~Kumar, and
  L.~Saminadayar.
\newblock Quantum shot noise.
\newblock {\em Superlattices and Microstructures}, 23(3):901--915, 1998.

\bibitem{rice1977generalized}
J.~Rice.
\newblock On generalized shot noise.
\newblock {\em Advances in Applied Probability}, 9(3):553--565, 1977.

\bibitem{rosinski1989path}
J.~Rosi{\'n}ski.
\newblock On path properties of certain infinitely divisible processes.
\newblock {\em Stochastic Processes and their Applications}, 33(1):73--87,
  1989.

\bibitem{rosinski1990series}
J.~Rosi{\'n}ski.
\newblock {On series representations of infinitely divisible random vectors}.
\newblock {\em The Annals of Probability}, 18(1):405--430, 1990.

\bibitem{rosinski1991class}
J.~Rosi{\'n}ski.
\newblock {On a class of infinitely divisible processes represented as mixtures
  of Gaussian processes}.
\newblock In {\em Stable Processes and Related Topics}, pages 27--41. Springer,
  1991.

\bibitem{rosinski2001series}
J.~Rosi{\'n}ski.
\newblock {Series representations of L{\'e}vy processes from the perspective of
  point processes}.
\newblock In {\em L{\'e}vy Processes}, pages 401--415. Springer, 2001.

\bibitem{rosinski2007tempering}
J.~Rosi{\'n}ski.
\newblock Tempering stable processes.
\newblock {\em {Stochastic Processes and their Applications}}, 117(6):677--707,
  2007.

\bibitem{rosinski2008simulation}
J.~Rosi{\'n}ski.
\newblock {Simulation of L{\'e}vy Processes}.
\newblock In F.~Ruggeri, R.~S. Kenett, and F.~W. Faltin, editors, {\em
  Encyclopedia of Statistics in Quality and Reliability: Computationally
  Intensive Methods and Simulation}. John Wiley \& Sons, 2008.

\bibitem{rosinski2018representations}
J.~Rosi{\'n}ski.
\newblock Representations and isomorphism identities for infinitely divisible
  processes.
\newblock {\em The Annals of Probability}, 46(6):3229--3274, 2018.

\bibitem{rubenthaler2003improved}
S.~Rubenthaler and M.~Wiktorsson.
\newblock {Improved convergence rate for the simulation of stochastic
  differential equations driven by subordinated L{\'e}vy processes}.
\newblock {\em Stochastic Processes and their Applications}, 108(1):1--26,
  2003.

\bibitem{samorodnitsky1996class}
G.~Samorodnitsky.
\newblock A class of shot noise models for financial applications.
\newblock In {\em Athens Conference on Applied Probability and Time Series
  Analysis}, pages 332--353. Springer, 1996.

\bibitem{samorodnitsky1991construction}
G.~Samorodnitsky and M.~S. Taqqu.
\newblock {Construction of multiple stable measures and integrals using LePage
  representation}.
\newblock In {\em Stable Processes and Related Topics}, pages 121--141.
  Springer, 1991.

\bibitem{samorodnitsky1994stable}
G.~Samorodnitsky and M.~S. Taqqu.
\newblock {\em {Stable Non-Gaussian Random Processes}}.
\newblock Chapman \& Hall, New York, 1994.

\bibitem{sato1999levy}
K.~Sato.
\newblock {\em {L{\'e}vy Processes and Infinitely Divisible Distributions}}.
\newblock {Cambridge University Press}, United Kingdom, 1999.

\bibitem{scherer2012shot}
M.~Scherer, L.~Schmid, and T.~Schmidt.
\newblock Shot-noise driven multivariate default models.
\newblock {\em European Actuarial Journal}, 2(2):161--186, 2012.

\bibitem{sears2012photon}
A.~Sears, A.~Petrenko, G.~Catelani, L.~Sun, H.~Paik, G.~Kirchmair, L.~Frunzio,
  L.~Glazman, S.~Girvin, and R.~Schoelkopf.
\newblock {Photon shot noise dephasing in the strong-dispersive limit of
  circuit QED}.
\newblock {\em Physical Review B}, 86(18):180504, 2012.

\bibitem{stanislavsky2009transport}
A.~Stanislavsky and K.~Weron.
\newblock {Transport of magnetic bright points on the Sun. Analysis of
  subdiffusion scenarios}.
\newblock {\em Astrophysics and Space Science}, 323:351--355, 2009.

\bibitem{talagrand1993regularity}
M.~Talagrand.
\newblock Regularity of infinitely divisible processes.
\newblock {\em The Annals of Probability}, 21(1):362--432, 1993.

\bibitem{tankov2016levy}
P.~Tankov.
\newblock L{\'e}vy copulas: review of recent results.
\newblock In {\em The Fascination of Probability, Statistics and their
  Applications}, pages 127--151. Springer, 2016.

\bibitem{vervaat1979stochastic}
W.~Vervaat.
\newblock On a stochastic difference equation and a representation of
  non-negative infinitely divisible random variables.
\newblock {\em Advances in Applied Probability}, 11(4):750--783, 1979.

\bibitem{walker2000miscellanea}
S.~Walker and P.~Damien.
\newblock {Miscellanea. Representations of L{\'e}vy processes without Gaussian
  components}.
\newblock {\em Biometrika}, 87(2):477--483, 2000.

\bibitem{westcott1976existence}
M.~Westcott.
\newblock On the existence of a generalized shot-noise process.
\newblock In E.~J. Williams, editor, {\em Studies in Probability and
  Statistics. Papers in Honour of Edwin JG Pitman, North-Holland, Amsterdam},
  pages 73--88. North-Holland, Amsterdam, 1976.

\bibitem{wiktorsson2002simulation}
M.~Wiktorsson.
\newblock {Simulation of stochastic integrals with respect to L{\'e}vy
  processes of type G}.
\newblock {\em Stochastic Processes and their Applications}, 101(1):113--125,
  2002.

\bibitem{wilt2013photon}
B.~A. Wilt, J.~E. Fitzgerald, and M.~J. Schnitzer.
\newblock Photon shot noise limits on optical detection of neuronal spikes and
  estimation of spike timing.
\newblock {\em Biophysical Journal}, 104(1):51--62, 2013.

\bibitem{yuan2020asymptotic}
S.~Yuan and R.~Kawai.
\newblock Asymptotic degeneracy and subdiffusivity.
\newblock {\em Journal of Physics A: Mathematical and Theoretical},
  53(9):095002, 2020.

\end{thebibliography}
}

\end{document}